\documentclass[a4paper,fleqn]{cas-sc}

\usepackage{placeins}

\usepackage{graphicx}
\usepackage{xcolor}

\usepackage[numbers]{natbib}

\usepackage{amsmath}
\usepackage{amssymb}
\usepackage{amsthm}
\usepackage{stmaryrd}
\SetSymbolFont{stmry}{bold}{U}{stmry}{m}{n} 

\numberwithin{equation}{section}

\newcommand{\R}{\mathbb{R}}

\newcommand{\Ee}{\mathcal{E}}

\newcommand{\Tt}{\mathcal{T}}
\newcommand{\Vv}{\mathcal{V}}

\newtheorem{remark}{Remark}

\begin{document}
\let\WriteBookmarks\relax
\def\floatpagepagefraction{1}
\def\textpagefraction{.001}

\shorttitle{Slip Wall Boundary Conditions in Curved Domains}
\shortauthors{N. Atallah et~al.}

\title[mode = title]{Slip Wall Boundary Conditions in Curved Domains for
                     Finite Element ALE Hydrodynamics}

\author[4]{Nabil Atallah}
\author[2]{Kevin Sweet}
\author[3]{Paul Moujaes}
\author[1]{Jan Nikl}
\author[1]{Ketan Mittal}
\author[1]{Vladimir Z. Tomov}
\cormark[1]
\ead{tomov2@llnl.gov}

\affiliation[1]{organization={Lawrence Livermore National Laboratory}}
\affiliation[2]{organization={Portland State University}}
\affiliation[3]{organization={Dortmund Technical University}}
\affiliation[4]{organization={Divergent 3D}}

\cortext[1]{Corresponding author}

\nonumnote{Performed under the auspices of the U.S. Department of Energy under
Contract DE-AC52-07NA27344 (LLNL-JRNL-2024639).}

\begin{abstract}
We present a high-order finite element arbitrary
Lagrangian-Eulerian (ALE) method for hydrodynamics with slip-wall
boundary conditions on analytically curved domains.
The method combines distinct Lagrangian, remesh, and remap phases
while preserving the geometric features of curved boundaries
and allowing tangential motion of boundary nodes.
In the Lagrangian phase, we use weak boundary condition enforcement to
impose the slip-wall condition without freezing tangential motion.
In the mesh optimization phase, we adapt TMOP
(target-matrix optimization paradigm) to analytic walls by replacing
selected boundary coordinate degrees of freedom with curve or surface
parameters.
The Newton solve then moves boundary nodes tangentially with
respect to the prescribed geometry while carrying the drift introduced by the
Lagrangian phase.
In the remap phase, we transfer the solution from the old mesh to the
optimized mesh by a conservative advection-based remap over pseudo-time.
We present numerical results in two and three dimensions demonstrating the
method on high-order curved meshes and showing that it preserves boundary
fidelity while producing the expected hydrodynamic solution features.
\end{abstract}


\maketitle


\section{Introduction}
\label{sec:introduction}

ALE hydrodynamics has a long history in computational hydrocodes
\cite{Hirt1974, Benson1992}.
We build on the high-order finite element methods
for Lagrangian and ALE hydrodynamics developed in \cite{Dobrev2012, ALE2018}.
In their current form, these methods have been demonstrated on
domains with straight, axis-aligned boundaries.
In this work, we extend them to analytically curved domains with
slip-wall boundary conditions.
This extension is nontrivial because the condition
$\boldsymbol{v}\cdot \boldsymbol{n}=0$ must be enforced while still permitting
tangential boundary motion, the high-order curved geometry must be preserved
during mesh optimization, and the remap must remain compatible with the
resulting tangential mesh displacement.

Several related boundary treatments have been developed for finite element
and ALE flow discretizations.
Early finite element work addressed the algebraic imposition of
normal and tangential velocity constraints, including slip conditions on
curved boundaries \cite{Engelman1982, Behr2004}. For high-order compressible
flow discretizations, wall boundary conditions have been studied primarily
from the viewpoint of stability, entropy stability, and implementation in
Euler and Navier--Stokes solvers \cite{Parsani2015, Kopriva2020, Chan2022}.
Slip-wall boundary conditions have also been considered in high-order methods
for compressible flows on curved domains, including shifted-boundary
corrections that retain a piecewise linear representation of the curved
geometry \cite{Ciallella2023}. This approach has recently been extended to
moving curved boundaries in an ALE framework \cite{Boscheri2025}.
Moving slip-wall boundary conditions have also been used in high-order ALE
methods for compressible flows on moving meshes \cite{Boscheri2017}.
More generally, the freedom to move ALE boundary nodes
tangentially while constraining their normal motion has been exploited as a
mechanism for controlling mesh deformation and improving mesh quality
\cite{Gadala2002, RiveroRodriguez2021}, and moving slip-wall conditions also
arise in ALE formulations for incompressible and weakly compressible flows
\cite{Busto2023}.
The present work addresses a different but related problem: the construction
of a unified wall-boundary treatment for high-order finite element
Lagrangian and ALE hydrodynamics in general curved geometries, in which weak
slip-wall enforcement, geometry-constrained tangential mesh optimization,
and conservative remap remain mutually compatible.
These issues are especially important when the boundary itself represents a
physical feature that should remain sharp and geometrically faithful
throughout the simulation.

Our approach combines ideas from several components of the existing
high-order methodology. In the Lagrangian phase, we adopt the weak boundary
condition enforcement developed in \cite{Nabil2024}, which was designed
precisely to handle slip-wall conditions on curved boundaries without
suppressing tangential motion. In the remesh phase, we use the TMOP
mesh-optimization framework \cite{TMOP2019, Knupp2012}, but replace selected
boundary coordinate degrees of freedom by analytic curve or surface
parameters.
This keeps boundary nodes on the prescribed geometry during the Newton
solve and preserves geometric features of the curved domain.
In the remap phase, we transfer the solution from the old mesh to the optimized
mesh by conservative advection over a pseudo-time interval.

The resulting method provides a high-order ALE treatment of slip walls on
curved domains in both two and three dimensions. We present numerical
examples that illustrate the behavior of the three phases, demonstrate the
preservation of curved boundaries and their geometric features, and show the
expected hydrodynamic solution structure on the resulting meshes.

The implementation uses open-source components. The Lagrangian hydrodynamics
driver is based on the Laghos miniapp \cite{laghos}, which is built on top of
the finite element library MFEM \cite{MFEM2024}.

Section \ref{sec:lagrange} summarizes the Lagrange phase, including the
high-order finite element discretization, the weak boundary treatment, and the
time integration. Sections \ref{sec:mesh_opt} and
\ref{sec:remap-phase} describe the mesh optimization and remap phases.
Finally, in Section \ref{sec:results}, we present numerical results.


\section{Lagrange Phase}
\label{sec:lagrange}

In this section, we summarize the Lagrangian phase used in our ALE method.
The formulation follows \cite{Nabil2024}; our goal here is to record the
ingredients that are needed later for the mesh optimization and remap
phases, with particular emphasis on the weak slip-wall treatment.

We begin with a reference configuration $\Omega_0 = \Omega(0) \subset \R^d$
and a time-dependent physical domain $\Omega(t) \subset \R^d$. The fluid
motion is described by the map
$\boldsymbol{\varphi}_t : \Omega_0 \to \Omega(t)$, with
$\boldsymbol{x} = \boldsymbol{\varphi}_t(\boldsymbol{x}_0)$.
We denote the boundary of $\Omega_0$ by $\Gamma_0$ with outward unit normal
$\boldsymbol{n}_0$, and the boundary of $\Omega(t)$ by $\Gamma(t)$
with outward unit normal $\boldsymbol{n}$.


\subsection{Governing equations and boundary conditions}

We use the non-conservative form of the Lagrangian Euler equations commonly
used in shock hydrodynamics \cite{Benson1992, Dobrev2012}:
\begin{subequations}
\begin{align}
  \frac{d\boldsymbol{x}}{dt} &= \boldsymbol{v}, \\
  \rho J &= \rho_0, \label{eq:strong-mass-conservation} \\
  \rho \frac{d \boldsymbol{v}}{dt} &= \nabla_{\boldsymbol{x}} \cdot \boldsymbol{\sigma}, \\
  \rho \frac{de}{dt} &= \nabla_{\boldsymbol{x}}\boldsymbol{v} : \boldsymbol{\sigma}.
\end{align}
\end{subequations}
Here, $\boldsymbol{v}$ is the velocity, $\rho$ is the density,
$\rho_0$ is the initial density, and
$J = \det(\nabla_{\boldsymbol{x}_0}\boldsymbol{\varphi}_t)$ is the
deformation Jacobian. We omit body forces, source terms, and heat fluxes,
since they are not used in the tests considered in this paper.
In our formulation, the Cauchy stress tensor is written as
\begin{equation}
  \boldsymbol{\sigma} = -pI + \boldsymbol{\sigma}_a,
\end{equation}
where $\boldsymbol{\sigma}_a$ denotes the artificial-viscosity contribution
and the thermodynamic pressure is determined by an equation of state
$p = EOS(\rho,e)$. For an ideal gas,
\begin{equation}
  p = (\gamma - 1)\rho e, \qquad
  c_s = \sqrt{\gamma(\gamma - 1)e},
\end{equation}
with $\gamma$ the adiabatic index and $c_s$ the sound speed.
On the boundary we impose slip-wall conditions,
\begin{subequations}
\begin{align}
  \boldsymbol{v} \cdot \boldsymbol{n} &= 0, \\
  \boldsymbol{\tau}_i \cdot (\boldsymbol{\sigma}\boldsymbol{n}) &= 0,
\end{align}
\end{subequations}
where $\boldsymbol{\tau}_i$ denotes a tangential direction. The key feature
of the weak formulation below is that it enforces the normal-velocity
condition without freezing tangential boundary motion.


\subsection{High-order finite element discretization}

The reference domain $\Omega_0$ is partitioned into quadrilateral elements
in 2D or hexahedral elements in 3D. We use the standard $Q_kQ_{k-1}$
discretization from \cite{Nabil2024}:
\begin{itemize}
  \item The kinematic variables $\boldsymbol{x}$ and $\boldsymbol{v}$ are
        approximated in the continuous Gauss-Lobatto space
        $\Vv = (Q_k)^d \subset (H^1(\Omega_0))^d$,
        with basis $\{w_i\}_{i=1}^{N_\Vv}$.
  \item The specific internal energy is approximated in the discontinuous
        space $\Ee = \widehat{Q}_{k-1} \subset L^2(\Omega_0)$ using
        Gauss-Legendre nodes and Bernstein basis functions
        $\{\phi_i\}_{i=1}^{N_\Ee}$.
\end{itemize}
The basis functions are transported by the Lagrangian map:
\begin{equation}
  w_i(\boldsymbol{x},t) = w_i \circ \boldsymbol{\varphi}_t^{-1}(\boldsymbol{x}),
  \qquad
  \phi_j(\boldsymbol{x},t) = \phi_j \circ \boldsymbol{\varphi}_t^{-1}(\boldsymbol{x}).
\end{equation}
Using Einstein notation, the discrete unknowns are written as
\begin{align}
  \boldsymbol{x}(\boldsymbol{x}_0,t) &=
    \text{x}_{i,\ell}(t) w_i(\boldsymbol{x}_0)\boldsymbol{\zeta}_\ell, \\
  \boldsymbol{v}(\boldsymbol{x}_0,t) &=
    \text{v}_{i,\ell}(t) w_i(\boldsymbol{x}_0)\boldsymbol{\zeta}_\ell, \\
  e(\boldsymbol{x}_0,t) &= e_j(t)\phi_j(\boldsymbol{x}_0),
\end{align}
where $\boldsymbol{\zeta}_\ell$ denotes the $\ell$th Cartesian basis
vector,
with $\text{v}_{i,\ell} = d\text{x}_{i,\ell}/dt$. Density is recovered from
the geometry using strong mass conservation:
\begin{equation}
  \rho(\boldsymbol{x},t) =
    \frac{\rho_0(\boldsymbol{x}_0)}{J(\boldsymbol{x}_0,t)}.
\end{equation}


\subsection{Weak form and boundary terms}

The weak form below follows \cite{Nabil2024}; the detailed motivation,
analysis, and numerical study of the boundary terms can be found there.
The semi-discrete system can be written in terms of the mass matrices
$\mathbf{M}_{\Vv}$ and $\mathbf{M}_{\Ee}$,
\begin{subequations}
\begin{align}
  (\mathbf{M}_\Vv)_{i_{\ell}j_{m}} &= \int_\Omega \rho w_i w_j \delta_{\ell m} +
                                      \int_{\Gamma_0} \alpha_0 \rho_{\max} L
                                      w_i n_{0,\ell} w_j n_{0,m}, \\
  (\mathbf{M}_\Ee)_{ij} &= \int_\Omega \rho \phi_i \phi_j,
\end{align}
\end{subequations}
and the operator
\begin{equation}
  (\mathbf{F})_{i_\ell j} =
    \int_\Omega \sigma_{\ell m}\phi_j \nabla_{x_m} w_i -
    \int_\Gamma (n_m \sigma_{mk} n_k) n_\ell w_i \phi_j +
    \int_\Gamma \beta \rho c_s (\boldsymbol{v}\cdot\boldsymbol{n}) w_i n_\ell \phi_j.
\end{equation}
Here $\rho_{\max} = \sup_{\boldsymbol{x}_0 \in \Omega_0}\rho_0$,
$L$ is the perimeter of the bounding box of the domain, and
$\alpha_0 = \beta L / (\det(\nabla_{\hat{\boldsymbol{x}}}
\boldsymbol{x}_0))^{1/d}$ with $\beta = \lambda C_I$ and
\begin{equation}
C_I = \begin{cases}
        \dfrac{(k+1)(k+d)}{d},
          & \text{for simplices in } d \text{ dimensions}, \\
        (k+1)^2,
          & \text{for quadrilaterals and hexahedra}.
\end{cases}
\end{equation}

The resulting semi-discrete momentum and energy equations are
\begin{equation}
  \mathbf{M}_\Vv \frac{d\mathbf{v}}{dt} = -\mathbf{F}\mathbf{1}, \qquad
  \mathbf{M}_\Ee \frac{d\mathbf{e}}{dt} = \mathbf{F}^T\mathbf{v}.
\end{equation}
The weak enforcement of the slip-wall condition enters through three
boundary contributions.
The projected traction term in $\mathbf{F}$ is motivated directly by the
slip-wall condition $\boldsymbol{\tau}_i\cdot(\boldsymbol{\sigma}\boldsymbol{n}) = 0$:
after integration by parts, only the normal component of the traction should remain
on the boundary, so the full boundary traction is replaced by its normal projection.
The boundary term in $\mathbf{M}_\Vv$ then adds control of the normal component
of the velocity, which is the quantity that must vanish at the wall.
More precisely, this term enforces $\dot{\boldsymbol{v}}\cdot\boldsymbol{n}_0 = 0$
in the initial configuration, which is equivalent to
$d(\boldsymbol{v}\cdot\boldsymbol{n})/dt = 0$
for the fixed wall boundaries considered here.
Finally, the impedance-like term $\beta \rho c_s (\boldsymbol{v}\cdot\boldsymbol{n})$
penalizes any residual normal motion, with the factor $\rho c_s$ providing the
characteristic scaling used in \cite{Nabil2024}.
The scaling of $\alpha_0 = \beta L/J_{\hat{\boldsymbol{x}}\to\boldsymbol{x}_0}^{1/d}$
is chosen for dimensional consistency and time-independence; in the referenced
paper it is also identified as critical for obtaining accurate tangential
motion under mesh refinement.
Together, these terms weakly enforce the wall
condition while preserving tangential freedom of motion.


\subsection{Viscosity and time integration}

For shock calculations, we retain the artificial viscosity treatment from
\cite{Dobrev2012, Nabil2024}.
Since our present work does not modify that component, we
focus here on the boundary terms and the ALE coupling that follows.

Time integration in the Lagrange phase uses a modified midpoint
second-order Runge-Kutta method from compatible Lagrangian hydrodynamics
\cite{Caramana1998, Barlow2008, Scovazzi2008, Dobrev2012, Sandu2021},
which conserves discrete total energy exactly and
discrete total linear momentum up to $O(h^{m+1})$.
Denoting a quantity at time $t_n$ by a superscript $n$, the update is
\begin{align*}
  \mathbf{v}^{n + \frac{1}{2}} &= \mathbf{v}^n - \frac{\Delta t}{2}
        \mathbf{M}_{\Vv}^{-1}\mathbf{F}^n \mathbf{1}, &
  \mathbf{v}^{n+1} &= \mathbf{v}^n - \Delta t \mathbf{M}_{\Vv}^{-1}
        \mathbf{F}^{n + \frac{1}{2}}\mathbf{1}, \\
  \mathbf{e}^{n + \frac{1}{2}} &= \mathbf{e}^n + \frac{\Delta t}{2}
        \mathbf{M}_{\Ee}^{-1}\left(\left(\mathbf{F}^n\right)^T
        \mathbf{v}^{n + \frac{1}{2}}\right), &
  \mathbf{e}^{n+1} &= \mathbf{e}^n + \Delta t \mathbf{M}_{\Ee}^{-1}
        \left(\left(\mathbf{F}^{n + \frac{1}{2}}\right)^T
        \left(\frac{\mathbf{v}^{n} + \mathbf{v}^{n+1}}{2}\right)\right), \\
  \mathbf{x}^{n + \frac{1}{2}} &= \mathbf{x}^n + \frac{\Delta t}{2}
        \mathbf{v}^{n+\frac{1}{2}}, &
  \mathbf{x}^{n+1} &= \mathbf{x}^n + \Delta t
        \left(\frac{\mathbf{v}^n + \mathbf{v}^{n+1}}{2}\right).
\end{align*}
The role of this RK2-average scheme is important in the weak-boundary
analysis. Proposition 1 of \cite{Nabil2024} shows that it conserves
discrete total linear momentum up to $O(h^{m+1})$, where the only defect
comes from the weakly imposed condition
$\boldsymbol{v}\cdot\boldsymbol{n}=0$.
Proposition 2 shows that the same scheme conserves discrete total energy
exactly.
Since our tests omit body forces, source terms, and heat fluxes, this exact
energy-conservation result applies directly to the Lagrange phase used in this paper.


\section{Mesh Optimization Phase}
\label{sec:mesh_opt}

During the Lagrange phase, the mesh follows the physical motion of the material,
which can gradually degrade element quality and thereby force smaller stable
time steps.
Therefore, a mesh optimization step is performed before remapping the solution.
Our remesh phase starts from TMOP \cite{TMOP2019, Knupp2012}, where local
element geometry is optimized relative to prescribed target quality.
TMOP has a broader set of capabilities, including high-order $r$-adaptivity
\cite{TMOP2020, TMOP2026}, $h$-adaptivity \cite{TMOP2022, TMOP2025},
$p$-adaptivity \cite{TMOP2024}, displacement limiting and fitting,
and GPU-oriented implementations \cite{TMOP2023}.
For the present boundary treatment, we need only the core ingredients:
high-order mesh maps, targets, quality metrics, the variational objective,
and the nonlinear solve.


\subsection{Basic TMOP formulation}

We briefly fix the notation used later for the parameter derivatives.
Let $\mathcal{M}$ denote the computational mesh, and let
$E \in \mathcal{M}$ be a physical element with reference element $\bar E$ and
target element $E_t$. The geometric map of $E$ is
\begin{equation}
    \Phi_E(\bar{\boldsymbol{x}})
        = \boldsymbol{x}(\bar{\boldsymbol{x}})
        = \sum_{i=1}^{N_k} \boldsymbol{x}_i
        w_i(\bar{\boldsymbol{x}}),
\end{equation}
where $k$ is the mesh order, $N_k$ is the number of geometric degrees of
freedom inside an element, $\boldsymbol{x}_i$ are the nodal positions,
and $w_i$ are the geometric basis functions. At each quadrature point
$\bar{\boldsymbol{x}}_j$, $j=1,\ldots,N_q$, let
\begin{equation}
    A_j =
    \frac{\partial \Phi_E}{\partial \bar{\boldsymbol{x}}}
    (\bar{\boldsymbol{x}}_j),
    \qquad
    W_j =
    \frac{\partial \Phi_{E_t}}{\partial \bar{\boldsymbol{x}}}
    (\bar{\boldsymbol{x}}_j),
    \qquad
    T_j = A_j W_j^{-1}.
\end{equation}
Here $A_j$ is the physical Jacobian, $W_j$ is the target Jacobian, and
$T_j$ measures the physical element relative to its target. We use the
following shape metrics to penalize skewness and aspect ratio:
\begin{equation}
    \mu_2(T) = \frac{1}{2}\frac{|T|^2}{\det T} - 1
    \quad \text{in 2D},
    \qquad
    \mu_{302}(T) = \frac{|T|^2 |T^{-1}|^2}{9} - 1
    \quad \text{in 3D},
\end{equation}
where $|T|$ is the Frobenius norm. Many other metric choices are possible,
depending on whether the mesh motion should control shape, size, orientation,
or combinations of these quantities \cite{Knupp2020, Knupp2023}. We restrict
the examples to shape metrics, which are sufficient for demonstrating the
curved-boundary treatment and the resulting tangential boundary motion.

The mesh objective used in the computations is
\begin{equation}
\label{eq_tmop_F}
  F(\boldsymbol{x}) = \sum_{E \in \mathcal{M}}
                      \int_{E_t} \omega(x_t)\mu(T(x_t))\,d x_t +
                      \sum_{E \in \mathcal{M}}
                     \int_{E_t} \xi(x-x_0,\delta(x_0))\,d x_t,
\end{equation}
where $x=x(\bar{x})=x(\bar{x}(x_t))$ is determined by the mesh degrees of
freedom, $\omega$ is a spatial weight, and the second integral is the
displacement-limiting term used in the implementation. The reference position
is $x_0=x_0(\bar{x})=x_0(\bar{x}(x_t))$, and $\delta(x_0)$ is the local
limiting distance. In our tests
$\xi(y,\delta)=|y|^2/(2\delta^2)$. The quadrature and Newton iteration are
otherwise the standard TMOP ones: node positions are updated until
$\partial F(\boldsymbol{x})/\partial \boldsymbol{x}=0$ is satisfied.


\subsection{Boundary parametrization}

Standard TMOP treats each mesh coordinate as an independent unknown. On an
analytically curved wall this is not appropriate for boundary degrees of
freedom, because a free Cartesian update can move them off the prescribed
geometry. Our approach starts from an analytic parametrization of each
constrained boundary component:
\begin{itemize}
    \item $x(t),y(t)$ for two-dimensional curves.
    \item $x(t),y(t),z(t)$ for three-dimensional curves.
    \item $x(u,v),y(u,v),z(u,v)$ for three-dimensional surfaces.
\end{itemize}
We assume that, over the part of the boundary used by the mesh, this
parametrization gives a one-to-one correspondence between Cartesian
boundary positions and parameter values. This correspondence is needed in
both directions. Given the Cartesian position of a boundary node after the
Lagrange phase, we determine its parameter value and use that value to
initialize the optimization. Conversely, during TMOP, each trial parameter
value is mapped back to Cartesian space so that the element maps, element
Jacobians, and mesh quality metrics can be evaluated.

The optimization unknowns for boundary nodes are then the analytic
parameters rather than the full set of Cartesian coordinates. In 2D, for
example, a boundary node is represented as
\begin{equation}
    \boldsymbol{x}_i = S(t_i) =
    \begin{pmatrix} x(t_i) \\ y(t_i) \end{pmatrix}.
\end{equation}
The Newton solve optimizes the scalar parameter $t_i$ directly instead of
the two Cartesian coordinates $x_i$ and $y_i$.
In 3D, nodes on boundary edges use a curve parametrization
$\boldsymbol{x}_i = S_e(t_i)$, while nodes in the interior of boundary
faces use a surface parametrization
$\boldsymbol{x}_i = S_f(u_i,v_i)$. Thus edge nodes have one optimization
unknown and face nodes have two. Nodes not constrained to the analytic
boundary retain their usual Cartesian coordinates.

The result is still one global TMOP optimization problem: the boundary
parameters and the unconstrained interior Cartesian degrees of freedom are
solved for simultaneously by Newton's method. The only change is the
representation of selected boundary unknowns. The TMOP derivatives are
therefore modified by the chain rule and include derivatives of the
boundary parametrization, as described next.


\subsection{Derivatives with parameters}

The TMOP kernels still evaluate the same element maps and Jacobians in
Cartesian space. The difference is that some Cartesian coordinates are no
longer independent optimization unknowns.
Let $x_a$ denote the finite element mesh coordinate field in
direction $a$, and let $x_{a i}$ be the coefficient of
basis function $w_i$ in that coordinate field.
For an unconstrained node, $x_{a i}$ remains an active degree of freedom.
For a constrained boundary node,
the Cartesian position is reconstructed from the active parameter vector
$\boldsymbol{p}_i$, e.g., $\boldsymbol{p}_i=(u_i,v_i)$ for a 3D surface node
and $\boldsymbol{p}_i=(t_i)$ for a curve node:
\begin{equation}
  x_{a i}=S_a(\boldsymbol{p}_i),
  \qquad a=1,\ldots,d.
\end{equation}
Thus the element map and Jacobian are still
\begin{equation}
  x_a(\bar{\boldsymbol{x}})
  = \sum_i x_{a i}w_i(\bar{\boldsymbol{x}}),
  \qquad
  A_{ab}
  = \sum_i x_{a i}
  \frac{\partial w_i}{\partial \bar{x}_b},
  \qquad
  T = A W^{-1},
\end{equation}
but the derivatives used by Newton must be taken with respect to the active
unknowns, not necessarily the Cartesian coordinates.

For an unconstrained Cartesian degree of freedom, the standard TMOP derivative
is unchanged:
\begin{equation}
\label{eq:cartesian-A-derivative}
  \frac{\partial A_{cb}}{\partial x_{a i}}
  = \delta_{ac}
  \frac{\partial w_i}{\partial \bar{x}_b},
  \qquad a,b,c=1,\ldots,d.
\end{equation}
For a parameterized boundary node, the chain rule gives
\begin{equation}
  \frac{\partial A_{ab}}{\partial p_{r i}}
  =
  \frac{\partial S_a}{\partial p_r}(\boldsymbol{p}_i)
  \frac{\partial w_i}{\partial \bar{x}_b},
\end{equation}
where $p_{r i}$ is one component of the parameter vector $\boldsymbol{p}_i$.
If two derivatives are taken with respect to parameters of the same node,
the curvature of the parametrization also contributes:
\begin{equation}
\label{eq:param-second-A}
  \frac{\partial^2 A_{ab}}{\partial p_{s i}\partial p_{r i}} =
  \frac{\partial^2 S_a}{\partial p_s\partial p_r}
  (\boldsymbol{p}_i) \frac{\partial w_i}{\partial \bar{x}_b}.
\end{equation}
For distinct nodes, $A$ has no second derivative with respect to the
parameters themselves; the coupling comes only through the nonlinearity of the
metric $\mu(T)$.

Let $G = \partial\mu/\partial T$ and $H=\partial^2\mu/\partial T^2$ at a
quadrature point.
Define
\begin{equation}
  D_{r i} =
  \frac{\partial A}{\partial p_{r i}}W^{-1}.
\end{equation}
Then the metric derivative with respect to a boundary parameter is
\begin{equation}
  \frac{\partial \mu}{\partial p_{r i}}
  = G : D_{r i}.
  \label{eq:tmop-param-grad}
\end{equation}
For two parameterized nodes $i$ and $j$, the off-diagonal Hessian contribution
is
\begin{equation}
  \frac{\partial^2 \mu}
  {\partial p_{s j}\partial p_{r i}}
  =
  H[D_{r i},D_{s j}],
  \qquad i\neq j.
\end{equation}
For two parameters of the same boundary node, the additional curvature term is
\begin{equation}
  \frac{\partial^2 \mu}
  {\partial p_{s i}\partial p_{r i}}
  =
  H[D_{r i},D_{s i}]
  +
  G :
  \left(
  \frac{\partial^2 A}{\partial p_{s i}\partial p_{r i}}W^{-1}
  \right).
  \label{eq:tmop-param-hess}
\end{equation}
Mixed Hessian blocks between a parameterized boundary node and an
unconstrained Cartesian coordinate are handled in the same way. If
\[
  D_{a j} =
  \frac{\partial A}{\partial x_{a j}}W^{-1},
\]
then
\begin{equation}
  \frac{\partial^2 \mu}{\partial x_{a j}\partial p_{r i}}
  =
  H[D_{r i},D_{a j}].
\end{equation}
These element-level contributions are assembled into the same global Newton
system as in standard TMOP, with parameter unknowns replacing the eliminated
boundary Cartesian unknowns. In 3D, nodes on the intersection of two boundary
surfaces can still move along the edge when that edge is represented by its own
curve parametrization. Nodes that are multiply marked without a unique curve or
surface parametrization are completely fixed.


\subsection{Accounting for Lagrange-phase displacement}

The weak boundary enforcement in the Lagrange phase allows tangential boundary
motion, but it does not keep the boundary nodes exactly on the analytic curve
or surface. Consequently, after the Lagrange step a boundary node may have a
small displacement away from the analytic boundary. The remesh phase does not
project this displacement back to the boundary. Instead, it keeps the
Lagrange-phase offset fixed and optimizes the tangential parameter.

There are two reasons for preserving this offset. First, a projection step is a
direct geometric manipulation outside the TMOP nonlinear solve. Even if the
boundary node is placed exactly on the analytic geometry, the surrounding
volume nodes are not adjusted by the optimization at the same time, and the
resulting mesh can contain inverted volumetric quadrature points. Second, the
remap formulation used in this work assumes that the mesh displacement has no
normal component on the boundary. Introducing a normal remesh displacement
would create boundary fluxes and would require a different remap weak form,
with corresponding changes to the conservation argument. We therefore move
boundary nodes tangentially with respect to the analytic parametrization while
preserving their stored offset from the boundary.

To do this, after the Lagrange phase we determine the boundary parameter
$\boldsymbol{p}_i$ associated with each boundary node, for example by the
closest point on the analytic curve or surface, and store the remaining
displacement as a scalar normal offset. In all cases, the offset is taken in
the normal direction of the analytic boundary. Thus the stored data are a
scalar offset $d_i$ and a parameter $\boldsymbol{p}_i$, rather than an
arbitrary Cartesian displacement, and
$\boldsymbol{d}_i=d_i\boldsymbol{n}(\boldsymbol{p}_i)$. This normal-offset
representation is the preferred choice because it separates the normal error
produced by the Lagrange phase from the tangential motion used in remeshing.
During Newton, $d_i$ is not an optimization unknown. Only the parameter
changes, and the Cartesian position is reconstructed through the offset
parametrization
\begin{equation}
  \boldsymbol{x}_i^{\mathrm{new}}
  =
  S(\boldsymbol{p}_i+\Delta\boldsymbol{p}_i)
  +
  d_i\boldsymbol{n}(\boldsymbol{p}_i+\Delta\boldsymbol{p}_i).
\end{equation}
Thus an off-boundary point remains on the same offset curve or surface while
moving tangentially with the optimized parameter. In the limit
$d_i=0$, the update reduces to the exact boundary map
$S(\boldsymbol{p}_i+\Delta\boldsymbol{p}_i)$.


\section{Remap Phase}
\label{sec:remap-phase}

After the mesh optimization phase, the state variables must be transferred
from the Lagrangian mesh to the optimized mesh. We use the standard
advection-based ALE remap \cite{Hirt1974, Dobrev2015, ALE2018}, in which the
mesh moves over a pseudo-time interval $\tau \in [0,1]$ according to
\begin{equation}
  \boldsymbol{x}(\tau)
  =
  \boldsymbol{x}^{\mathrm{old}}
  +
  \tau\boldsymbol{u},
  \qquad
  \boldsymbol{u}
  =
  \boldsymbol{x}^{\mathrm{new}}
  -
  \boldsymbol{x}^{\mathrm{old}}.
\end{equation}
The remap is conservative because density, internal energy density, and
momentum are transported by the corresponding advection equations:
\begin{align}
  \frac{\partial \rho}{\partial \tau}
  &=
  \boldsymbol{u}\cdot\nabla\rho, \\
  \frac{\partial(\rho e)}{\partial \tau}
  &=
  \boldsymbol{u}\cdot\nabla(\rho e), \\
  \frac{\partial(\rho\boldsymbol{v})}{\partial \tau}
  &=
  \boldsymbol{u}\cdot\nabla(\rho\boldsymbol{v}).
\end{align}
The boundary treatment in this weak form relies on the condition
\begin{equation}
  \boldsymbol{u}\cdot\boldsymbol{n}=0
  \qquad \text{on } \partial\Omega(\tau).
\end{equation}
This is one of the main reasons why the mesh optimization phase moves boundary
nodes tangentially instead of projecting them back to the analytic boundary.
If $\boldsymbol{u}$ had a normal component on the boundary, the remap would
introduce boundary fluxes, and the present weak form would no longer give the
same conservation statement without additional boundary terms or a modified
remap formulation.

For density and specific internal energy we use the finite element spaces and
moving bases from the Lagrange phase. With density basis functions
$\{\psi_i\}$ and energy basis functions $\{\phi_i\}$, the semi-discrete scalar
remap for density and specific internal energy has the form:
\begin{equation}
\label{eq:remap_semi_rho_e}
  M_\rho\frac{d\rho}{d\tau} = K_\rho \rho, \qquad
  M_e\frac{de}{d\tau} = K_e e.
\end{equation}
Let $\kappa \in \Tt(\tau)$ denote a mesh element at pseudo-time $\tau$.
The matrices are
\begin{align}
  (M_\rho)_{ij} &= \int_{\Omega(\tau)}\psi_i\psi_j, \quad
  (K_\rho)_{ij}  = \sum_{\kappa\in\Tt(\tau)} \int_\kappa
                   \psi_i(\boldsymbol{u}\cdot\nabla\psi_j) -
                   \sum_F \int_F (\boldsymbol{u}\cdot\boldsymbol{n}_F)
                          \psi_j^{\mathrm{up}}\llbracket\psi_i\rrbracket, \\
  (M_e)_{ij} &= \int_{\Omega(\tau)}\rho\phi_i\phi_j, \quad
  (K_e)_{ij}  = \sum_{\kappa\in\Tt(\tau)} \int_\kappa
                \rho\phi_i(\boldsymbol{u}\cdot\nabla\phi_j) -
                \sum_F \int_F \rho^{\mathrm{up}}
                (\boldsymbol{u}\cdot\boldsymbol{n}_F)
                \phi_j^{\mathrm{up}}\llbracket\phi_i\rrbracket.
\end{align}
Here $\llbracket\eta\rrbracket=\eta^- - \eta^+$ denotes the jump across an
interior face $F$, and the upwind trace is selected using the sign of
$\boldsymbol{u}\cdot\boldsymbol{n}_F$. On boundary faces, the normal flux
vanishes because $\boldsymbol{u}\cdot\boldsymbol{n}=0$.
These scalar remap equations are solved using the limited FCT algorithms
developed in earlier remap work. We do not repeat the limiter construction
here; details of the FCT remap algorithms can be found in
\cite{Dobrev2015, Manuel2017}.

\subsection{Momentum remap}
The momentum remap is performed in the continuous $H^1$ finite element space.
Using the density remap equation, the momentum advection equation can be
written componentwise in non-conservative form as
\begin{equation}
  \rho\frac{\partial\boldsymbol{v}}{\partial\tau} =
  \rho\boldsymbol{u}\cdot\nabla\boldsymbol{v}.
\end{equation}
We perform a change of basis from the kinematic basis functions $\{w_i\}$ to
the Bernstein basis $\{\phi_i\}$ and introduce the following index sets.
Let $\mathcal{N}_i$ and $\mathcal{K}_i$ be the sets of all nodes and elements
in the support of $\phi_i$, respectively.
Furthermore, we denote by $\mathcal{N}^\kappa$ the set of all nodes that
belong to element $\kappa$,
$\mathcal{N}_i^{\kappa} = \mathcal{N}^\kappa\cap\mathcal{N}_i$, and
$\mathcal{N}_i^{\kappa*} = \mathcal{N}_i^\kappa\setminus\{i\}$.
Testing with a Bernstein basis function $\phi_i$ yields, for each component
$1\leq \ell\leq d$,
\begin{equation}
\label{eq:sd_target}
  \sum_{j\in\mathcal{N}_i} m_{ij}
  \frac{\mathrm{d} v_{j,\ell}}{\mathrm{d}\tau}
  =
  \sum_{\kappa\in\mathcal{K}_i}\sum_{j\in\mathcal{N}_i}
  k^{\kappa}_{ij} v_{j,\ell}
  =
  \sum_{\kappa\in\mathcal{K}_i}
  \sum_{j\in\mathcal{N}_i^{\kappa*}}
  k^{\kappa}_{ij} (v_{j,\ell} - v_{i,\ell}),
\end{equation}
where
\begin{equation}
  m_{ij} = \int_{\Omega(\tau)}\rho \phi_i \phi_j, \quad \quad
  k_{ij}^\kappa = \int_{\kappa}
  \rho \phi_i(\boldsymbol{u}\cdot\nabla \phi_j).
\end{equation}
Here and below, the superscript $\kappa$ denotes an element-local contribution.
The second equality in \eqref{eq:sd_target} follows from the
partition of unity property of the Bernstein basis.

To apply the limiting procedure, we first derive a bound-preserving
low-order scheme of local Lax--Friedrichs type
\cite{kuzmin2012b, kuzmin2023}.
Lumping the mass matrix, we approximate
\begin{equation}
    \sum_{j\in\mathcal{N}_i} m_{ij}
    \frac{\mathrm{d} v_{j,\ell}}{\mathrm{d}\tau}
    \approx m_i \frac{\mathrm{d} v_{i,\ell}}{\mathrm{d}\tau},
\end{equation}
where
\begin{equation}
    m_i = \sum_{j\in\mathcal{N}_i} m_{ij} = \int_{\Omega(\tau)} \rho \phi_i.
\end{equation}
Owing to the positivity of the Bernstein functions, we have $m_i> 0$.
Furthermore, we add an element-based graph Laplacian to obtain the
semi-discrete low-order formulation
\begin{equation}\label{eq:LO1}
    m_i \frac{\mathrm{d} v_{i,\ell}}{\mathrm{d}\tau}
    = \sum_{\kappa\in\mathcal{K}_i}
      \sum_{j\in\mathcal{N}_i^{\kappa*}}
      (k^{\kappa}_{ij} + d_{ij}^\kappa) (v_{j,\ell} - v_{i,\ell}),
\end{equation}
where the graph viscosity coefficients are given by
\cite{lohmann2017, kuzmin2023}
\begin{equation}
    d_{ij}^\kappa =
    \begin{cases}
        \max(-k_{ij}^\kappa,\, -k_{ji}^\kappa,\, 0)
            &\text{if } i\in\mathcal{N}^{\kappa}
             \text{ and } j\in\mathcal{N}_i^{\kappa*},\\[1mm]
        -\sum_{l\in\mathcal{N}_i^{\kappa*}} d_{il}^\kappa
            &\text{if } i\in\mathcal{N}^{\kappa}
             \text{ and } j=i,\\[1mm]
        0 &\text{otherwise}.
    \end{cases}
\end{equation}
The coefficients are defined so that
$(D^\kappa)_{ij} = d_{ij}^\kappa$ has zero row sums and is symmetric.

Introducing
$d_i^\kappa=\sum_{j\in\mathcal{N}_i^{\kappa*}} d_{ij}^\kappa$,
we rewrite the low-order scheme~\eqref{eq:LO1} in the element-based
bar-state form
\begin{equation}
    m_i \frac{\mathrm{d} v_{i,\ell}}{\mathrm{d}\tau}
    =
    \sum_{\kappa\in\mathcal{K}_i}
    2d_i^\kappa (\bar v_{i,\ell}^\kappa - v_{i,\ell}),
\end{equation}
where the element bar-state
\begin{equation}
    \bar v_{i,\ell}^\kappa
    =
    \frac{1}{2d_i^\kappa}
    \sum_{j\in\mathcal{N}_i^{\kappa*}}
    2d_{ij}^\kappa \bar v_{ij,\ell}^\kappa
\end{equation}
is a convex combination of the intermediate states
\cite{kuzmin2020, kuzmin2023}
\begin{equation}
    \bar v_{ij,\ell}^\kappa
    =
    \frac{v_{j,\ell} + v_{i,\ell}}{2}
    +
    \frac{k_{ij}^\kappa(v_{j,\ell} - v_{i,\ell})}
         {2d_{ij}^\kappa},
    \quad j\in\mathcal{N}_i^{\kappa*}
\end{equation}
each of which is itself a convex combination of $v_{j,\ell}$ and $v_{i,\ell}$.
\begin{remark}
    In practice, the bar-states $\bar v_{ij,\ell}^\kappa$ and
    $\bar v_{i,\ell}^\kappa$ are
    never assembled.
    Instead, we evaluate the products
    \begin{equation}
        2d_{ij}^\kappa\,\bar v_{ij,\ell}^\kappa
        =
        d_{ij}^\kappa (v_{i,\ell} + v_{j,\ell})
        +
        k_{ij}^\kappa (v_{j,\ell} - v_{i,\ell}),
        \qquad
        2d_{i}^\kappa\,\bar v_{i,\ell}^\kappa
        =
        \sum_{j\in\mathcal{N}_i^{\kappa*}}
        2d_{ij}^\kappa\,\bar v_{ij,\ell}^\kappa
    \end{equation}
    to avoid division by zero.
\end{remark}

Discretizing in time with a forward Euler stage of a high-order strong
stability preserving (SSP) scheme yields
\begin{equation}\label{eq:LO2}
    v_{i,\ell}^{\mathrm{SSP}}
    =
    v_{i,\ell}
    +
    \frac{\Delta \tau}{m_i}
    \sum_{\kappa\in\mathcal{K}_i}
    2d_i^\kappa (\bar v_{i,\ell}^{\kappa} - v_{i,\ell})
    =
    \left(
    1 -
    \frac{\Delta \tau}{m_i}
    \sum_{\kappa\in\mathcal{K}_i} 2d_i^\kappa
    \right) v_{i,\ell}
    +
    \frac{\Delta \tau}{m_i}
    \sum_{\kappa\in\mathcal{K}_i}
    2d_i^\kappa \bar v_{i,\ell}^{\kappa}.
\end{equation}
Since the states $\bar v_{i,\ell}^{\kappa}$,
$\kappa\in\mathcal{K}_i$, are convex combinations of $v_{j,\ell}^n$,
$j\in\mathcal{N}_i^{\kappa}$, the state $v_{i,\ell}^{\mathrm{SSP}}$
is also a convex combination if the time step satisfies the CFL-like condition
\begin{equation}
    \frac{\Delta \tau}{m_i}
    \sum_{\kappa\in\mathcal{K}_i} 2d_i^\kappa
    \leq 1.
\end{equation}
Thus, the low-order step \eqref{eq:LO2} remains within the local bounds
\begin{equation}
    v_{i,\ell}^{\min} = \min_{j\in\mathcal{N}_i} v_{j,\ell}, \qquad v_{i,\ell}^{\max} = \max_{j\in\mathcal{N}_i} v_{j,\ell}.
\end{equation}

Adding the antidiffusive fluxes
\begin{equation}
    f_{ij,\ell}^\kappa
    =
    m_{ij}^\kappa(\dot v_{i,\ell} - \dot v_{j,\ell})
    +
    d_{ij}^\kappa (v_{i,\ell} - v_{j,\ell})
    =
    - f_{ji,\ell}^\kappa,
    \qquad
    m_{ij}^\kappa = \int_{\kappa}\rho\phi_i\phi_j,
\end{equation}
to the semi-discrete low-order scheme \eqref{eq:LO1} recovers the high-order
scheme~\eqref{eq:sd_target}.
We approximate $\dot v_{i,\ell}$ by the low-order time derivatives
\begin{equation}
    \dot v_{i,\ell}
    =
    \frac{1}{m_i}
    \sum_{\kappa\in\mathcal{K}_i}
    2d_i^\kappa (\bar v_{i,\ell}^\kappa - v_{i,\ell}),
\end{equation}
which yields the high-order stabilization~\cite{lohmann2019}.
Thus, the semi-discrete target scheme reads
\begin{equation}
    m_i \frac{\mathrm{d} v_{i,\ell}}{\mathrm{d}\tau}
    =
    \sum_{\kappa\in\mathcal{K}_i}
    \sum_{j\in\mathcal{N}_i^{\kappa*}}
    k^{\kappa}_{ij} (v_{j,\ell} - v_{i,\ell})
    =
    \sum_{\kappa\in\mathcal{K}_i}
    2d_{i}^\kappa (\bar v_{i,\ell}^{\kappa,H} - v_{i,\ell}),
\end{equation}
where the high-order bar-states are given by
\begin{equation}
    \bar v_{i,\ell}^{\kappa,H}
    =
    \bar v_{i,\ell}^\kappa
    +
    \frac{f_{i,\ell}^\kappa}{2d_i^\kappa},
    \qquad
    f_{i,\ell}^\kappa
    =
    \sum_{j\in\mathcal{N}_i^{\kappa*}} f_{ij,\ell}^\kappa.
\end{equation}
Owing to the skew-symmetry of
$f_{ij,\ell}^\kappa = - f_{ji,\ell}^\kappa$, the element contributions
satisfy $\sum_{ i\in\mathcal{N}^\kappa} f_{i,\ell}^\kappa = 0$, so
the flux correction is conservative.

To enforce local bound preservation, we replace the raw antidiffusive
contributions by limited counterparts
$f_{i,\ell}^{\kappa,*}\approx f_{i,\ell}^\kappa$, such that
\begin{equation}
    v_{i,\ell}^{\min}
    \leq
    \bar v_{i,\ell}^{\kappa}
    +
    \frac{f_{i,\ell}^{\kappa,*}}{2d_i^\kappa}
    \leq v_{i,\ell}^{\max}
    \quad \Leftrightarrow\quad
    f_{i,\ell}^{\kappa,\min}
    \leq f_{i,\ell}^{\kappa,*}
    \leq f_{i,\ell}^{\kappa,\max}
\end{equation}
where
\begin{equation}
    f_{i,\ell}^{\kappa,\min}
    :=
    2d_i^\kappa (v_{i,\ell}^{\min} - \bar v_{i,\ell}^{\kappa})
    \leq 0,
    \qquad
    f_{i,\ell}^{\kappa,\max}
    :=
    2d_i^\kappa (v_{i,\ell}^{\max} - \bar v_{i,\ell}^{\kappa})
    \geq 0.
\end{equation}
We apply the two-step \textit{Clip and Scale}
limiter~\cite{Manuel2017, kuzmin2023, lohmann2017}.
In the first step, we clip the raw contributions to their admissible bounds
\begin{equation}\label{eq:int_fluxes}
    \tilde f_{i,\ell}^{\kappa}
    =
    \min\left(
    \max(f_{i,\ell}^{\kappa}, f_{i,\ell}^{\kappa,\min}),
    f_{i,\ell}^{\kappa,\max}
    \right).
\end{equation}
However, the intermediate fluxes \eqref{eq:int_fluxes} generally do not
satisfy the conservation property
$\sum_{i\in \mathcal{N}^\kappa} \tilde f_{i,\ell}^{\kappa} = 0$.
Thus, we calculate the sum of the positive and negative clipped element
contributions
\begin{equation}
    P_\kappa^+
    =
    \sum_{i\in\mathcal{N}^\kappa}
    \max(0, \tilde f_{i,\ell}^\kappa),
    \quad
    P_\kappa^-
    =
    \sum_{i\in\mathcal{N}^\kappa}
    \min(0, \tilde f_{i,\ell}^\kappa),
\end{equation}
and define the limited fluxes by
\begin{equation}
    f_{i,\ell}^{\kappa,*}=
    \begin{cases}
        -\dfrac{P_\kappa^-}{P_\kappa^+}\tilde f_{i,\ell}^\kappa
            &\text{if }\tilde f_{i,\ell}^\kappa>0
             \text{ and } P_\kappa^+ + P_\kappa^- >0,\\[1mm]
        -\dfrac{P_\kappa^+}{P_\kappa^-}\tilde f_{i,\ell}^\kappa
            &\text{if }\tilde f_{i,\ell}^\kappa<0
             \text{ and } P_\kappa^+ + P_\kappa^- <0,\\[1mm]
        \tilde f_{i,\ell}^\kappa& \text{otherwise.}
    \end{cases}
\end{equation}
Replacing the low-order bar states $\bar v_{i,\ell}^{\kappa}$ in the forward
Euler update~\eqref{eq:LO2} by the limited high-order counterparts
$\bar v_{i,\ell}^{\kappa,*}
= \bar v_{i,\ell}^{\kappa}+f_{i,\ell}^{\kappa,*}/2d_{i}^\kappa$
yields a bound-preserving high-order method that suppresses spurious
oscillations in the vicinity of steep gradients.

After the remap step, the velocity is transformed back to the kinematic basis
functions $\{w_i\}$.
At boundary nodes, the fluid velocity may not satisfy the slip-wall
constraint exactly. In the present implementation this final normal component
is removed by an explicit correction step at the end of the remap procedure.


\section{Numerical Results}
\label{sec:results}

In this section, we illustrate the behavior of the method using variants of
the Sedov blast test \cite{Sedov1993} on several nontrivial 2D and 3D
domains. In all cases, the material is an ideal gas with $\gamma=1.4$.
The initial density is set to one and the initial velocity is zero throughout
the domain. A delta-function source of internal energy is deposited at a fixed
location, with total energy $E_{\mathrm{total}}=0.25$.

The main user-specified parameters are the ALE period, which determines the
duration of each Lagrangian phase, the remesh limiting distance $\delta$ in
\eqref{eq_tmop_F}, and the discretization order. The values used in
the test cases below are listed in Table \ref{tab_params}.

All simulations are performed with a customized version of the open-source
Laghos proxy application \cite{laghos}, which is based on the MFEM finite
element library \cite{MFEM2024}. Visualizations are performed with GLVis
\cite{glvis}.

\begin{table}
\centering
\begin{tabular}{ |c|c|c|c|c| }
  \hline
  \textbf{Test Case} & \textbf{Final Time} & \textbf{ALE Period} & \textbf{Limiting Distance} $\delta$ & \textbf{Discretization} \\
  \hline
  2D Linear Boundary      & 2.5   & 0.25  & 0.05  & $Q_3Q_2$  \\
  3D Bilinear Boundary    & 1.0   & 0.1   & 0.2   & $Q_2Q_1$  \\
  2D Sine Boundary        & 0.9   & 0.3   & 0.01  & $Q_3Q_2$  \\
  3D Sine Boundary        & 0.6   & 0.05  & 0.1   & $Q_2Q_1$  \\
  Annulus                 & 7.0   & 0.5   & 0.3   & $Q_3Q_2$  \\
  Solid Torus             & 4.8   & 0.2   & 0.2   & $Q_2Q_1$  \\
  \hline
\end{tabular}
\caption{User-chosen parameters for the six test cases presented.}
\label{tab_params}
\end{table}


\subsection{Linear Boundary}

We begin with simple non-axis-aligned linear boundaries in 2D and bilinear
boundary faces in 3D as a first test of the boundary treatment.

\begin{figure}[pos=htbp]
\centering
\begin{tabular}{c c c c}
  & Density & Velocity & Specific Internal Energy \\

  $t=0.5$
  & \raisebox{-.5\height}{\includegraphics[width=0.20\textwidth]{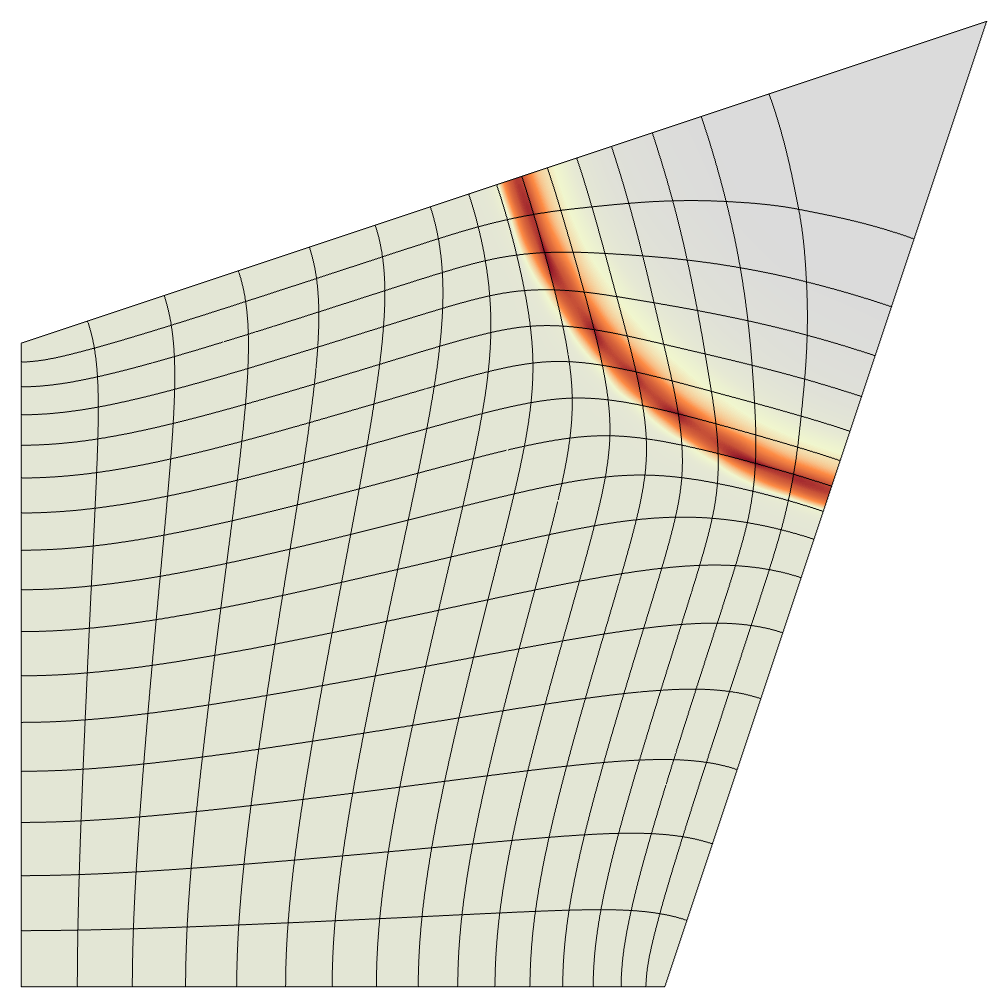}}
  & \raisebox{-.5\height}{\includegraphics[width=0.20\textwidth]{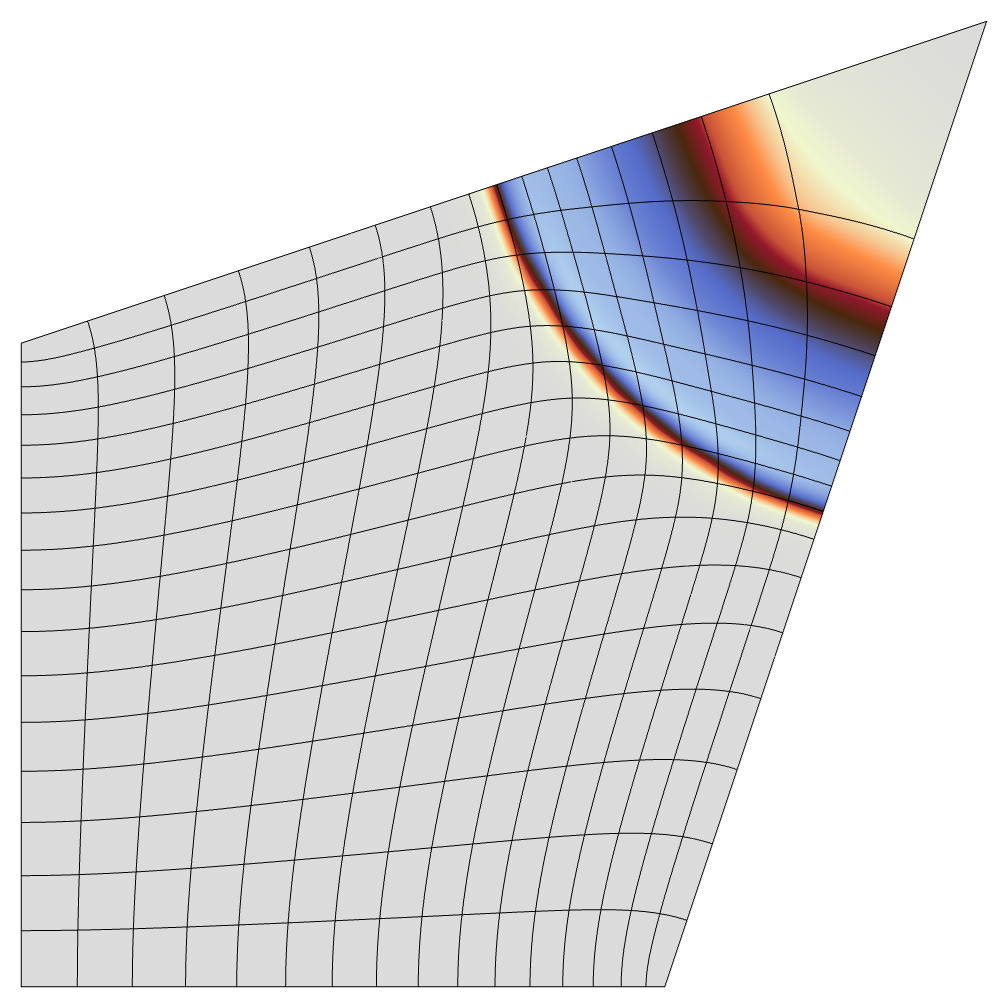}}
  & \raisebox{-.5\height}{\includegraphics[width=0.20\textwidth]{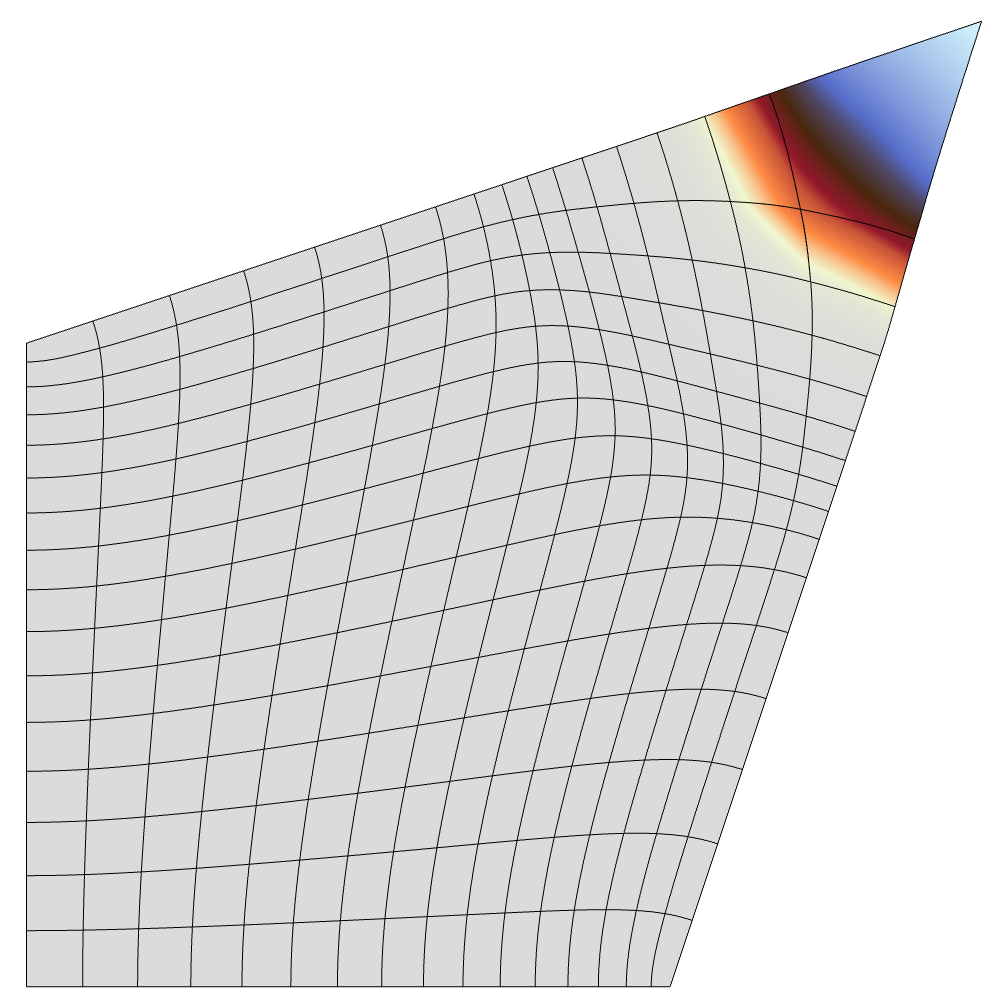}} \\

  $t=1.5$
  & \raisebox{-.5\height}{\includegraphics[width=0.20\textwidth]{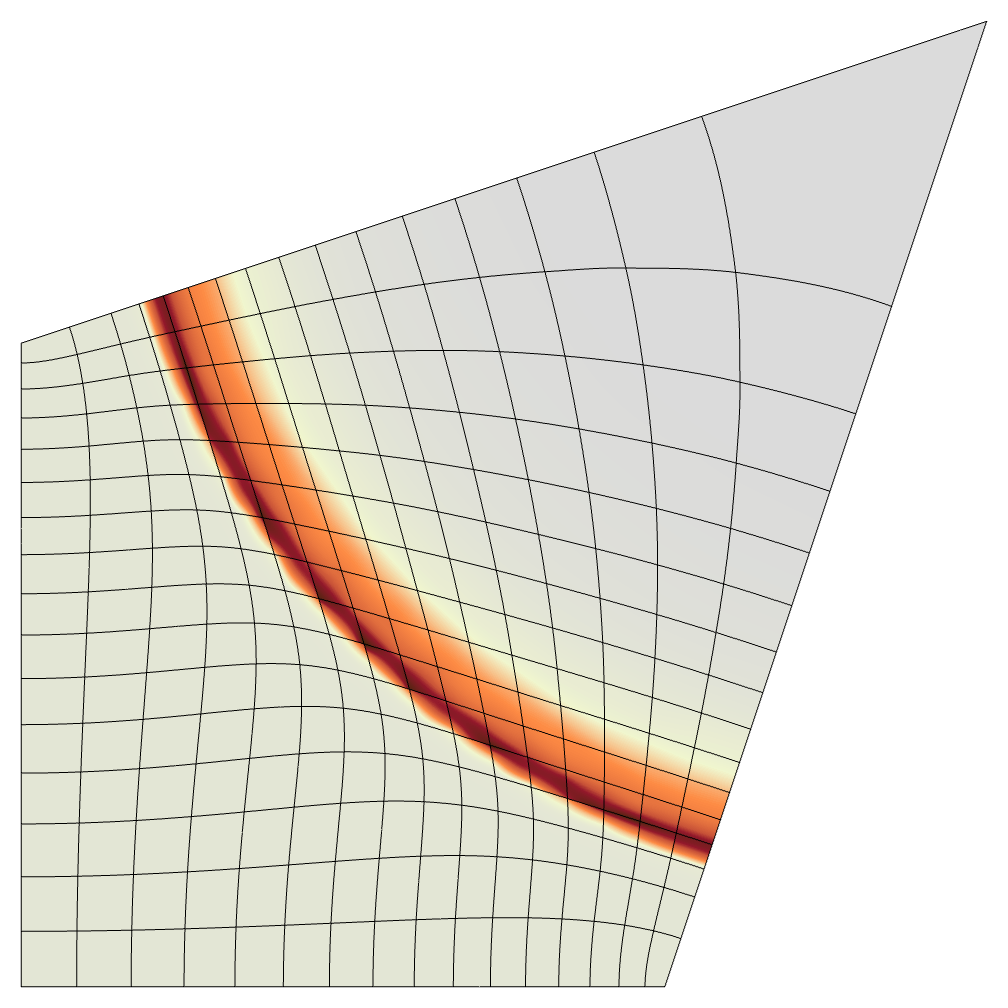}}
  & \raisebox{-.5\height}{\includegraphics[width=0.20\textwidth]{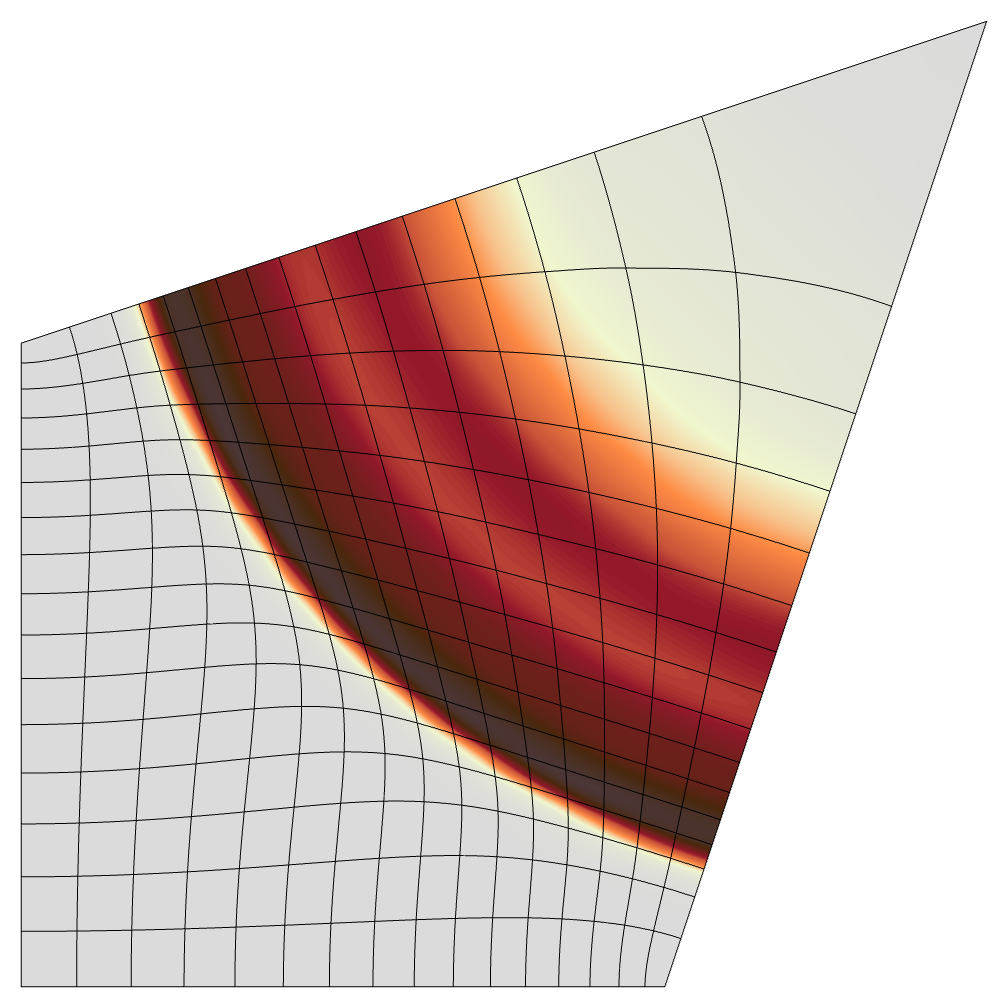}}
  & \raisebox{-.5\height}{\includegraphics[width=0.20\textwidth]{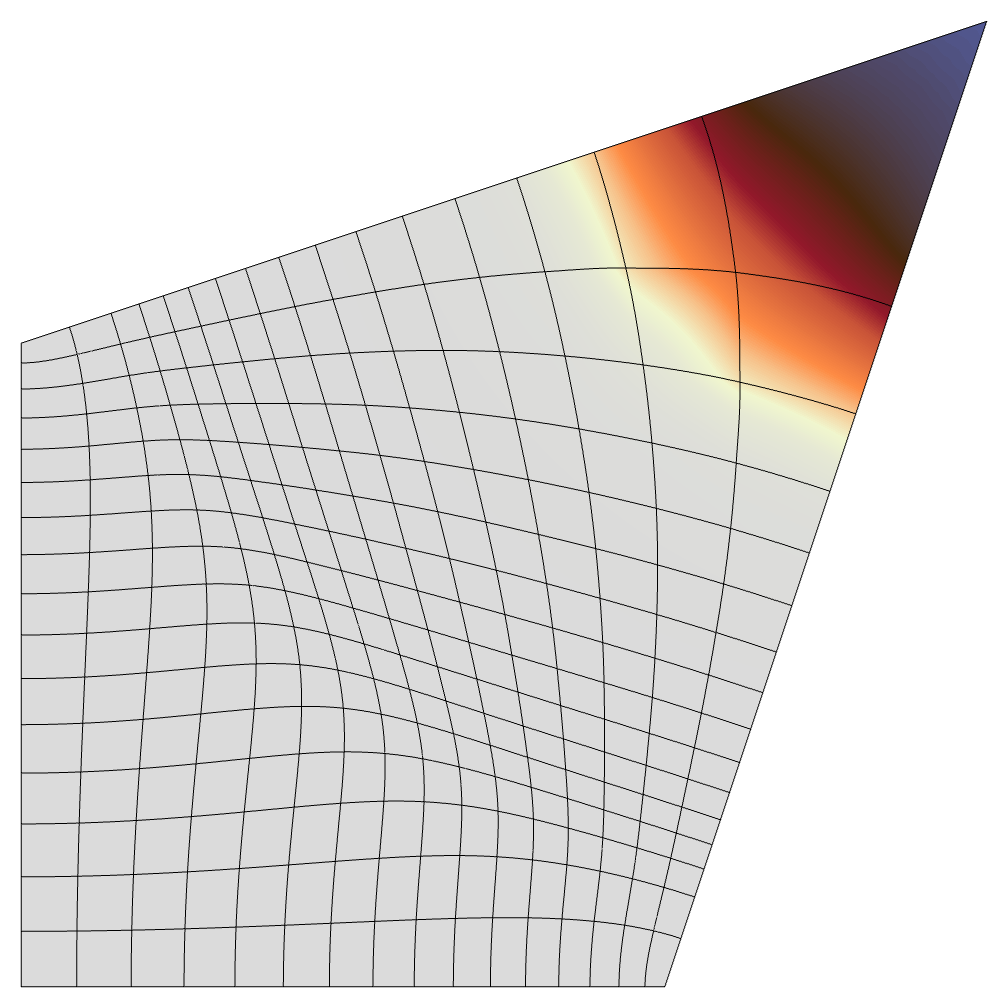}} \\

  $t=2.5$
  & \raisebox{-.5\height}{\includegraphics[width=0.20\textwidth]{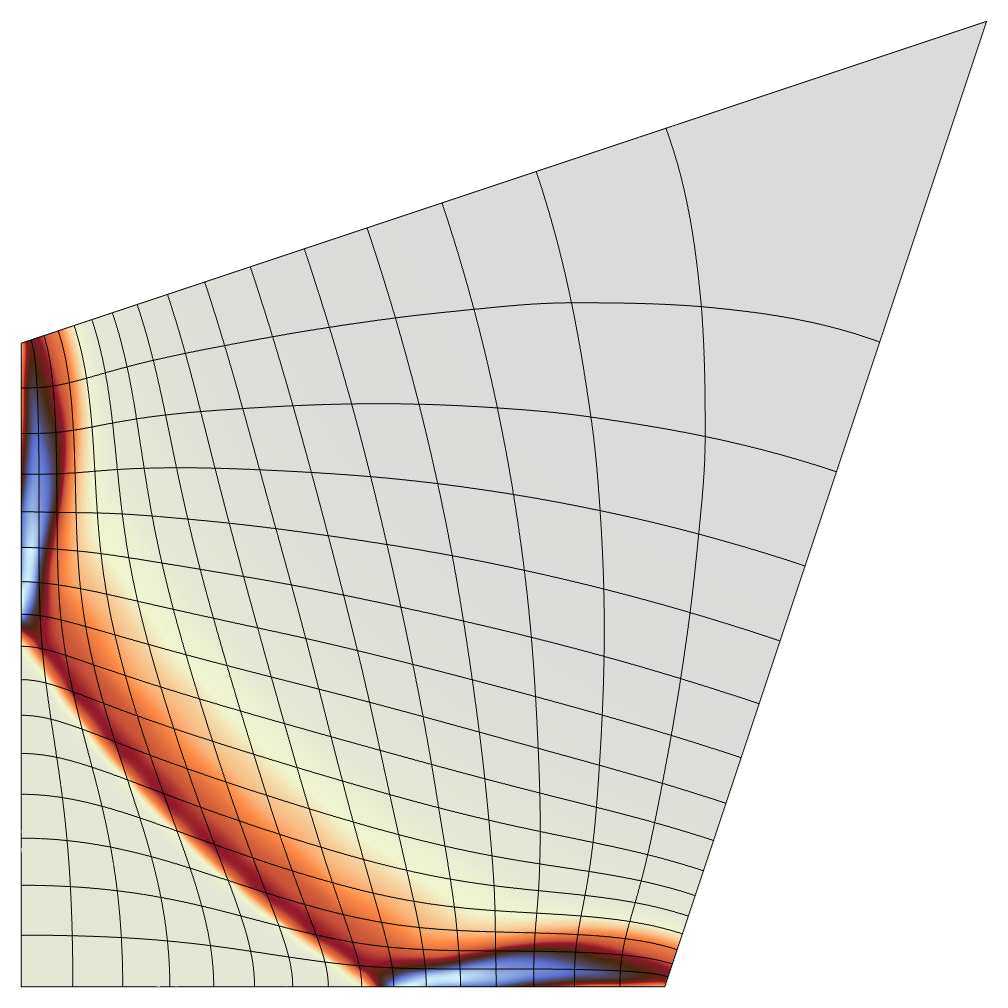}}
  & \raisebox{-.5\height}{\includegraphics[width=0.20\textwidth]{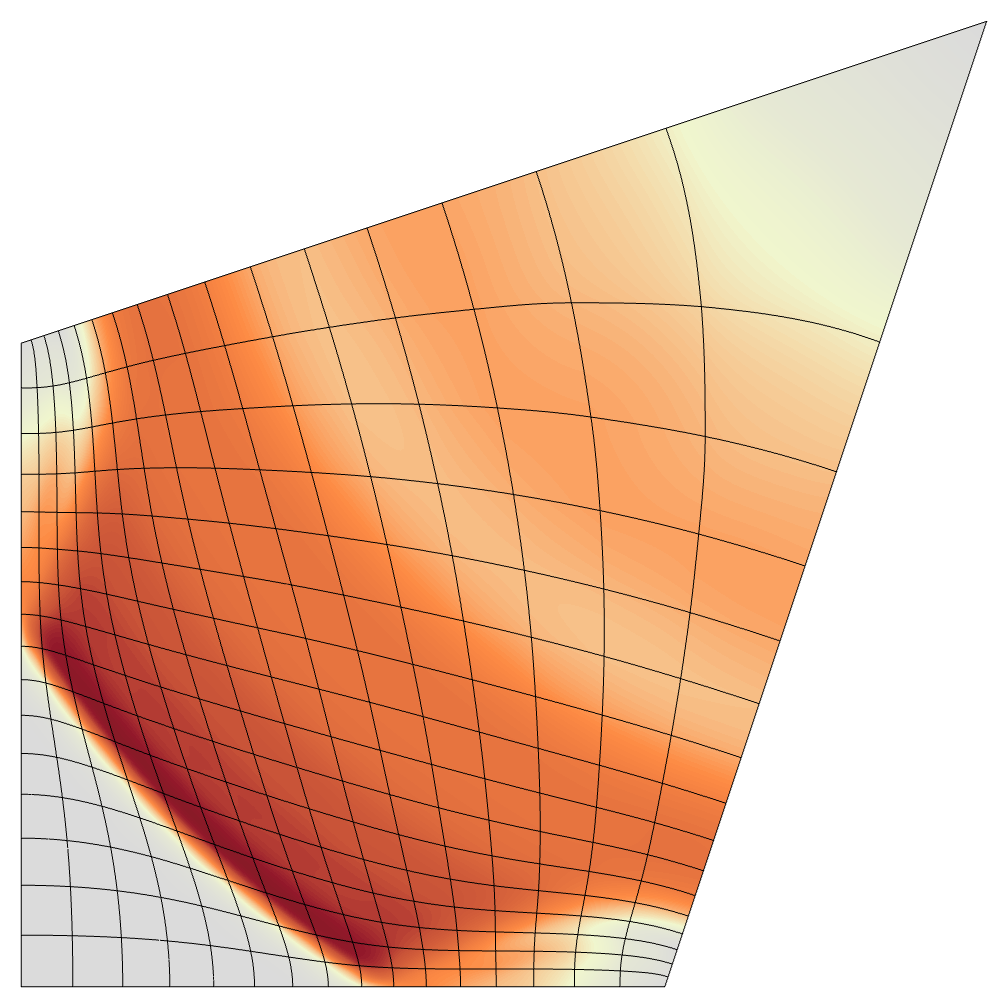}}
  & \raisebox{-.5\height}{\includegraphics[width=0.20\textwidth]{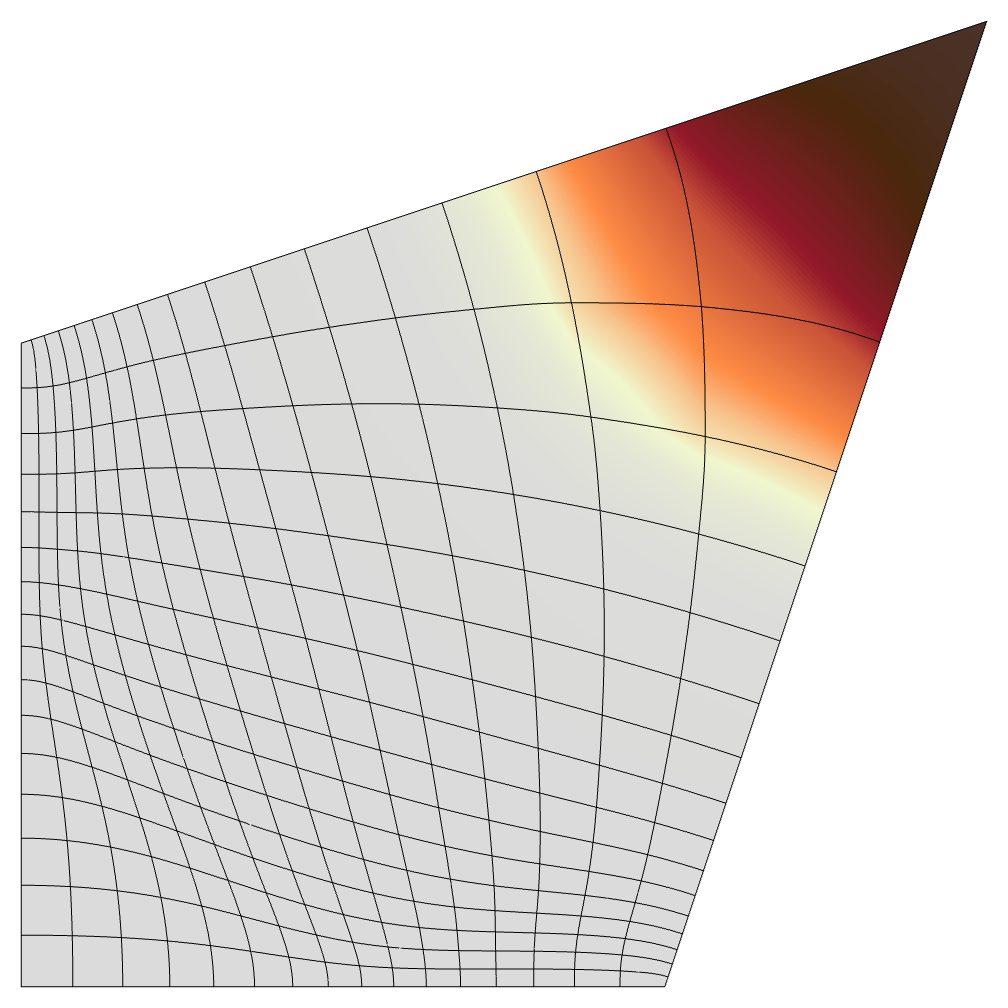}}
\end{tabular}
\caption{Plots of density, velocity, and specific internal energy for
         the Sedov blast on a 2D domain with linear boundaries.}
\label{fig:linear_2D}

\vspace{0.5em}

\centering
\begin{tabular}{c c c c}
  & Density & Velocity & Specific Internal Energy \\

  $t=0.5$
  & \raisebox{-.5\height}{\includegraphics[width=0.20\textwidth]{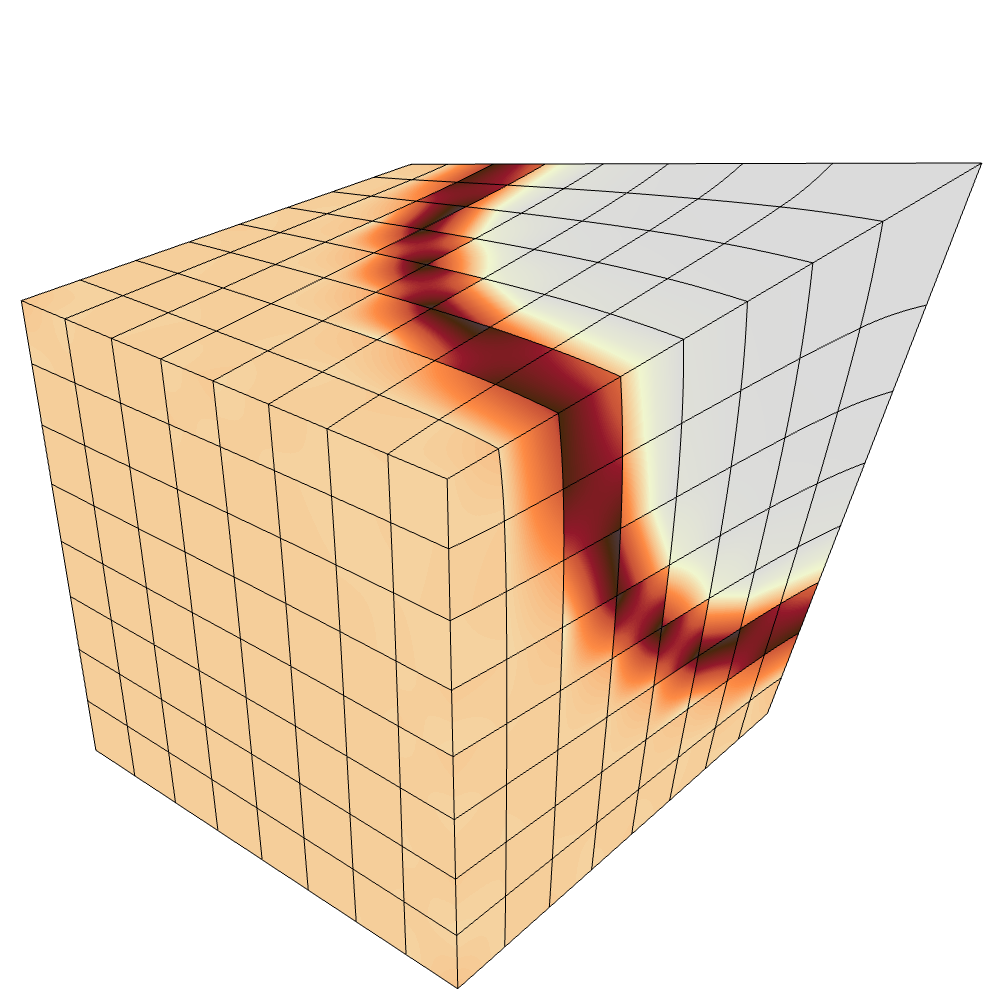}}
  & \raisebox{-.5\height}{\includegraphics[width=0.20\textwidth]{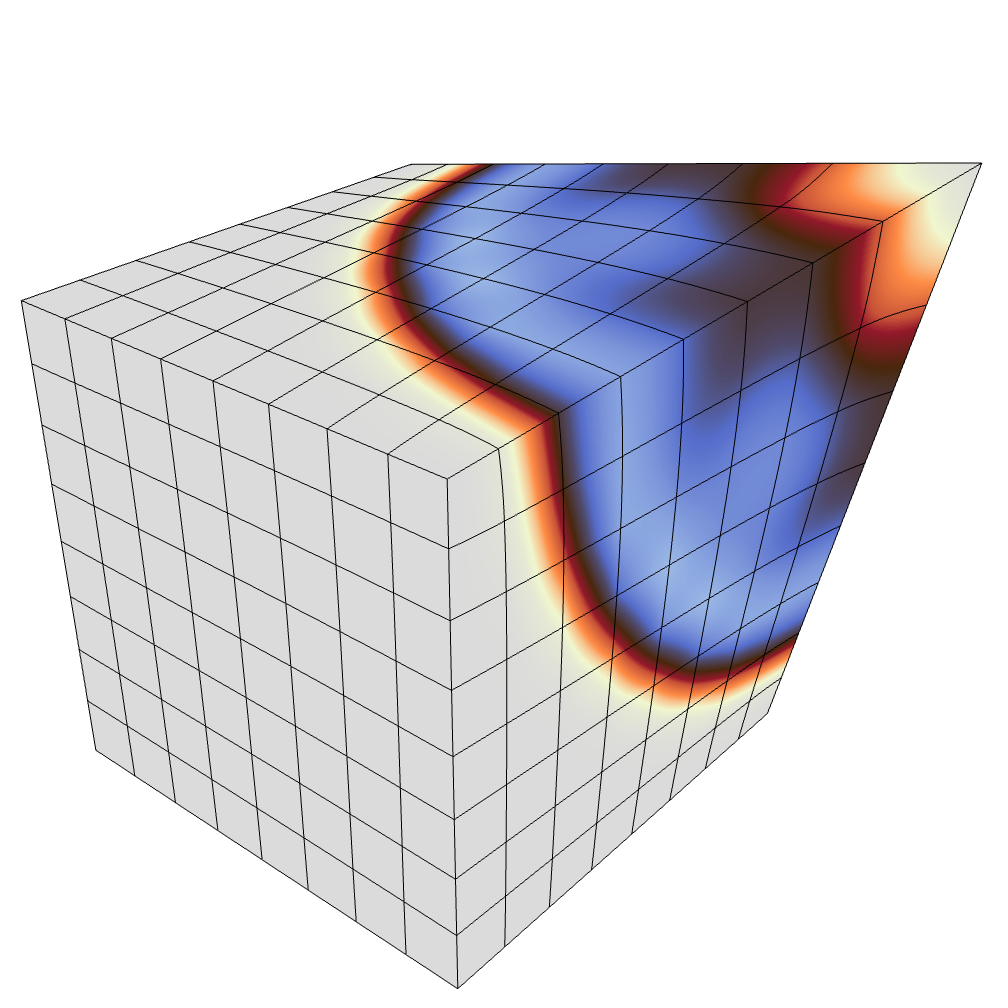}}
  & \raisebox{-.5\height}{\includegraphics[width=0.20\textwidth]{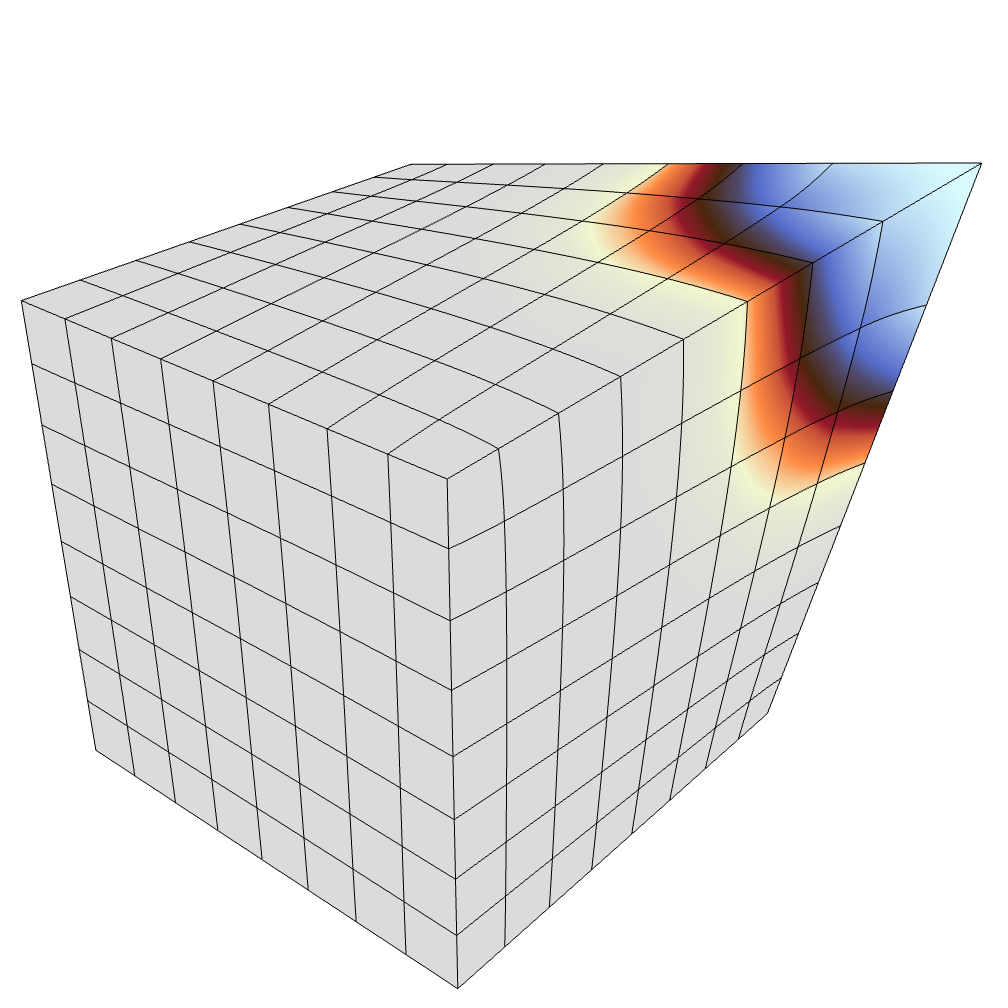}} \\

  $t=1.0$
  & \raisebox{-.5\height}{\includegraphics[width=0.20\textwidth]{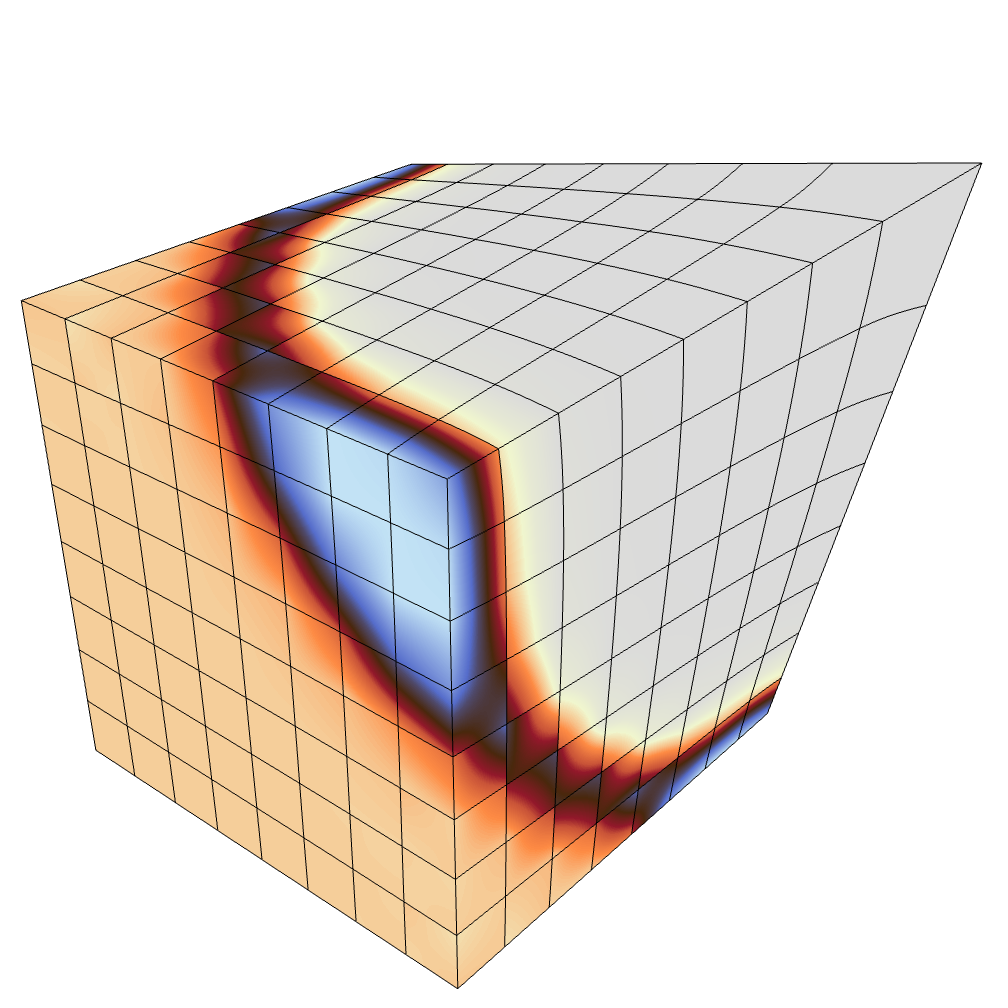}}
  & \raisebox{-.5\height}{\includegraphics[width=0.20\textwidth]{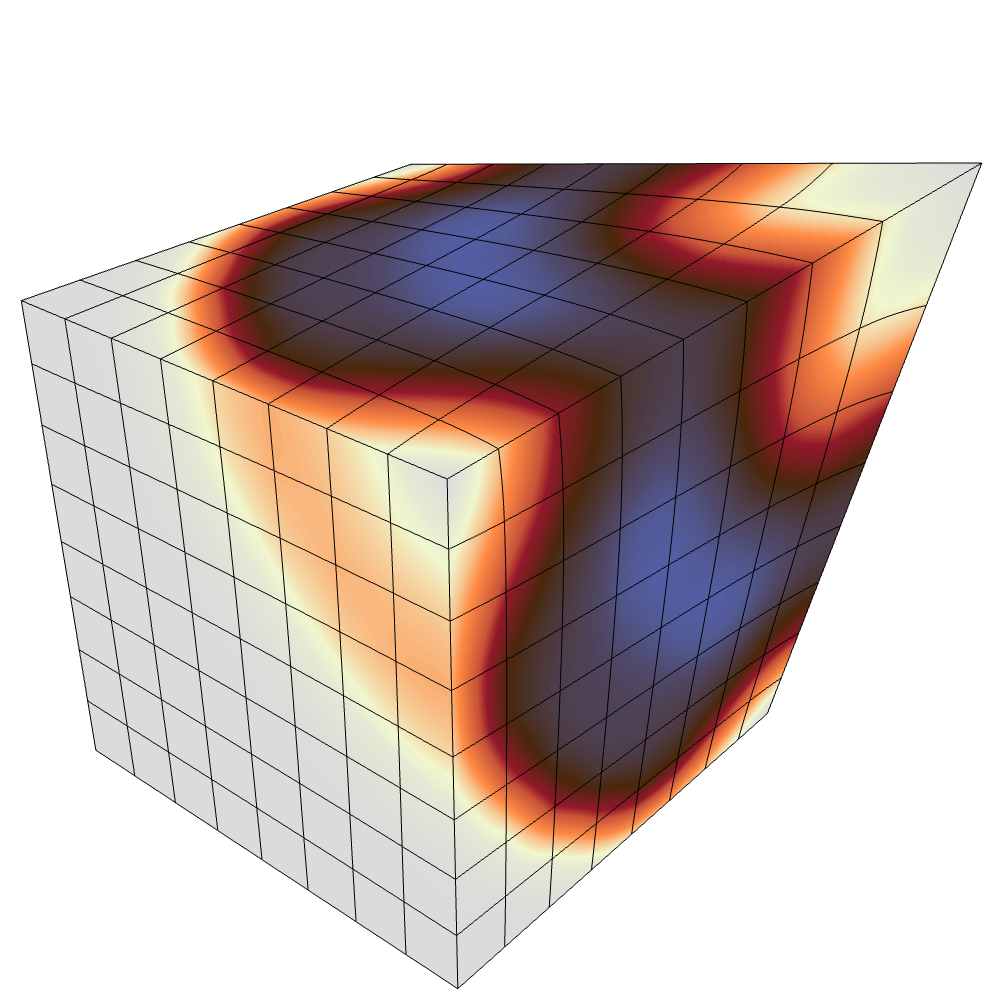}}
  & \raisebox{-.5\height}{\includegraphics[width=0.20\textwidth]{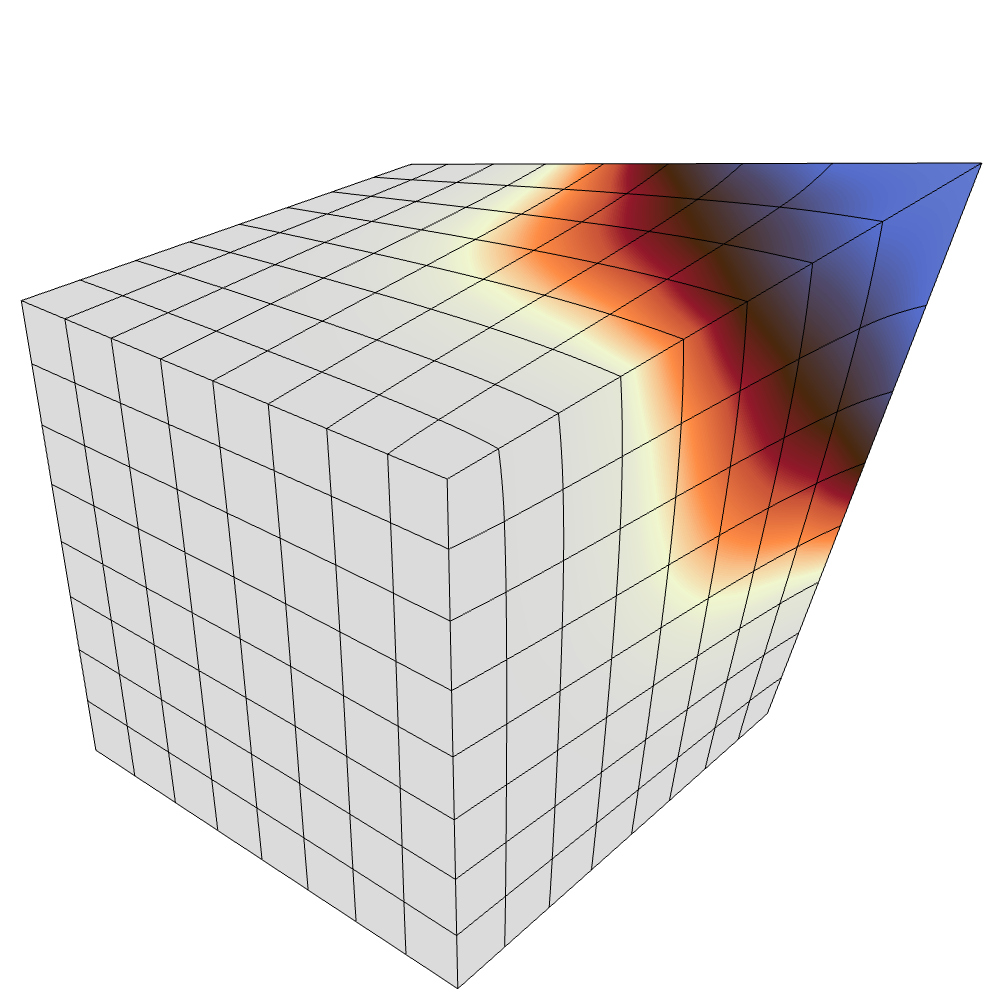}}
\end{tabular}
\caption{Plots of density, velocity, and specific internal energy for
         the Sedov blast on a cube with a ``pulled out'' corner.}
\label{fig:linear_3D}
\end{figure}

The 2D domain is given by the transformation
\[
\begin{bmatrix}
    x \\ y
\end{bmatrix} \mapsto \begin{bmatrix}
    x + \frac{1}{2}xy \\[0.3em]
    y + \frac{1}{2}xy
\end{bmatrix}
\]
applied to the unit square. This gives the following parameterizations
for the non-axis-aligned boundaries with $t \in [0,1]$:
\begin{center}
\begin{tabular}{l l l}
    Top boundary: & & Right boundary: \\
    $x(t) = 1.5t$ & & $x(t) = 1 + 0.5t$ \\
    $y(t) = 1 + 0.5t$ & & $y(t) = 1.5t$
\end{tabular}
\end{center}
The blast is initialized at the upper-right corner of the mesh. The evolution
shown in Figure \ref{fig:linear_2D} exhibits the expected near-wall behavior:
the solution slides along the boundary without visible mesh deformation or
spurious oscillations.

\begin{figure}[pos=htbp]
\centering
\begin{tabular}{c c c}
  \includegraphics[width=0.20\textwidth]{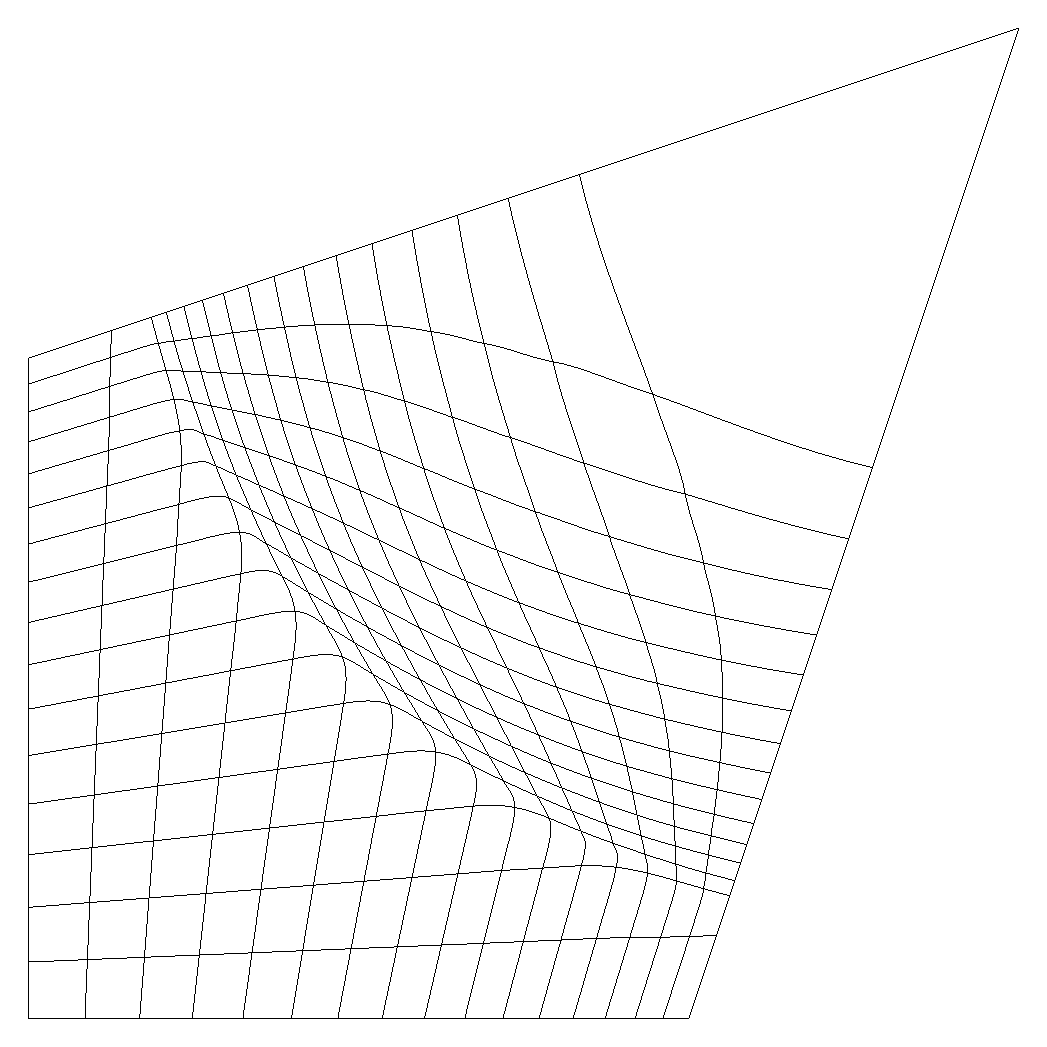} &
  \includegraphics[width=0.20\textwidth]{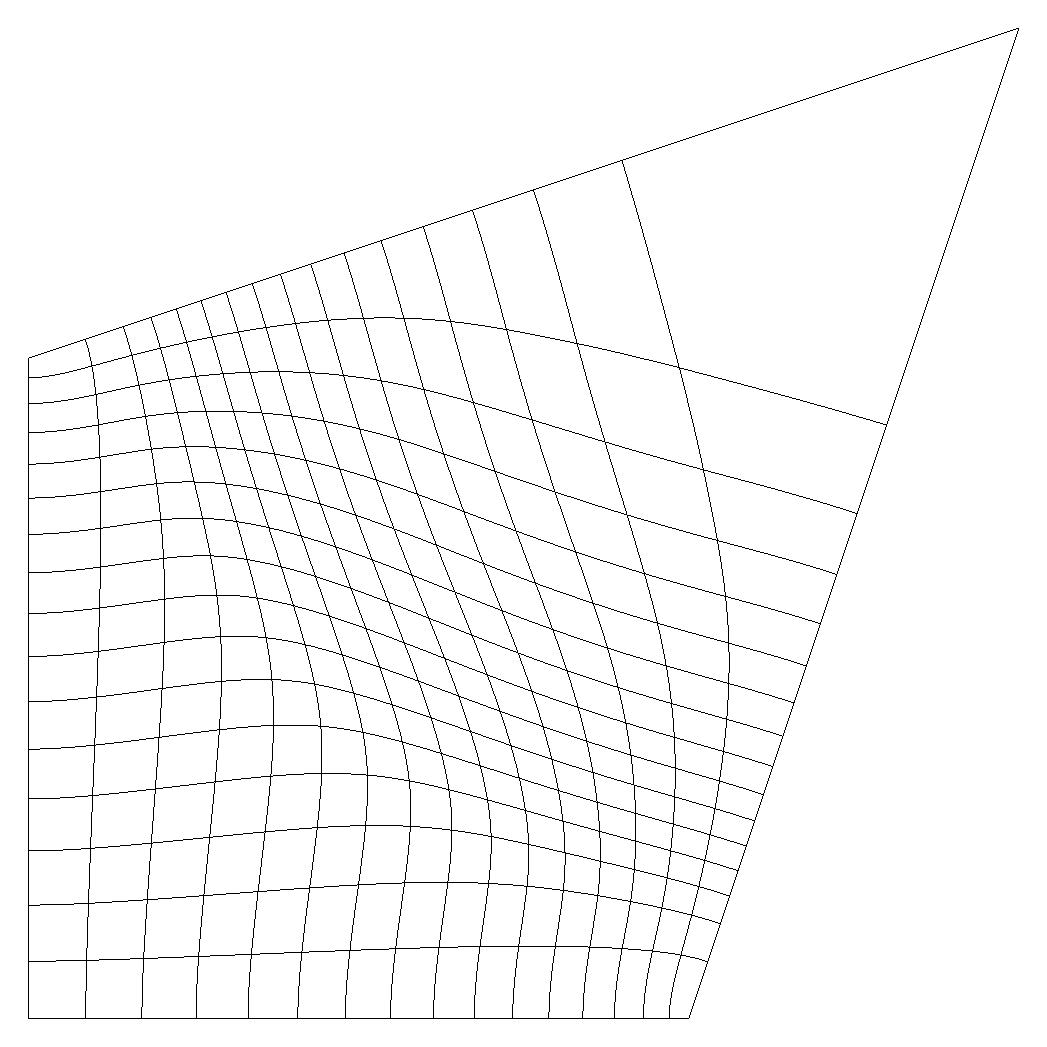} &
  \includegraphics[width=0.20\textwidth]{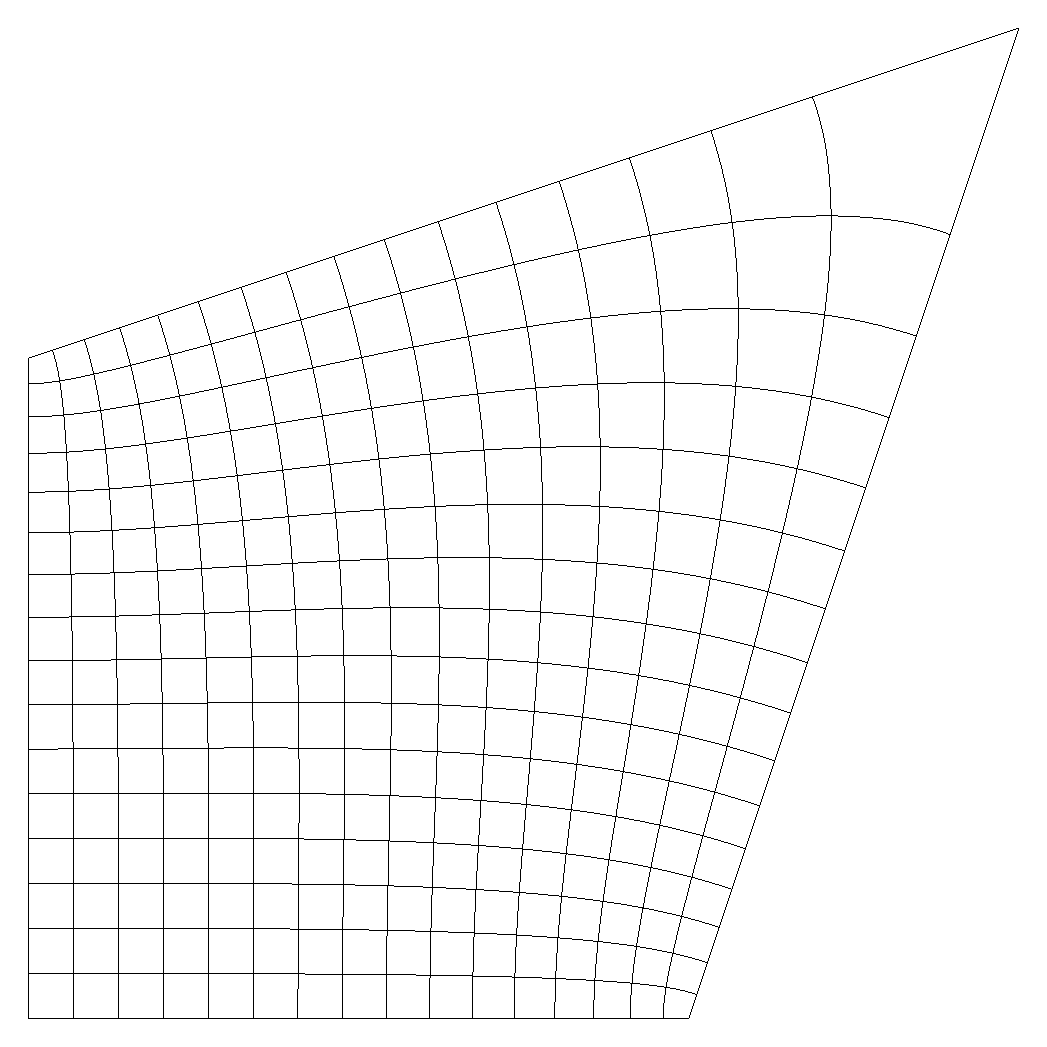}
\end{tabular}
\caption{The 2D mesh with linear boundaries before the ALE remesh phase (left),
         a mesh optimized with limiting distance $\delta = 0.05$ (middle),
         and a fully optimized mesh (right).}
\label{fig:linear_tmop}
\end{figure}

Figure \ref{fig:linear_tmop} shows the behavior of the TMOP-based
mesh optimizer on the 2D linear-boundary case.
The limited optimization improves the mesh while keeping the displacement close
to the Lagrangian mesh, whereas the fully unlimited optimization illustrates
the direction preferred by the quality functional.

The 3D domain with bilinear faces is given by the transformation
\[
    \begin{bmatrix}
        x \\ y \\ z
    \end{bmatrix} \mapsto \begin{bmatrix}
        x + 0.2xyz \\
        y + 0.2xyz \\
        z + 0.2xyz
    \end{bmatrix}
\]
applied to the unit cube. This gives the following parameterizations
with $u,v,t \in [0,1]$:
\begin{center}
\begin{tabular}{l l l}
    Top face: & & Edge between top and front faces: \\
    $x(u,v) = u + 0.2uv$ & & $x(t) = 1.2t$ \\
    $y(u,v) = v + 0.2uv$ & & $y(t) = 1 + 0.2t$ \\
    $z(u,v) = 1 + 0.2uv$ & & $z(t) = 1 + 0.2t$
\end{tabular}
\end{center}
The other nonlinear faces and edges are defined similarly, and the blast is
initialized at the ``pulled out'' corner of the mesh. The 3D results in
Figure \ref{fig:linear_3D} show the same qualitative behavior: the flow follows
the bilinear boundary features with the expected sliding motion, without
visible boundary-induced mesh deformation or solution oscillations.

\subsection{Sine Boundary}

\begin{figure}[pos=htbp]
\centering
\begin{tabular}{c c c c}
  & Density & Velocity & Specific Internal Energy \\

  $t=0.3$
  & \raisebox{-.5\height}{\includegraphics[width=0.20\textwidth]{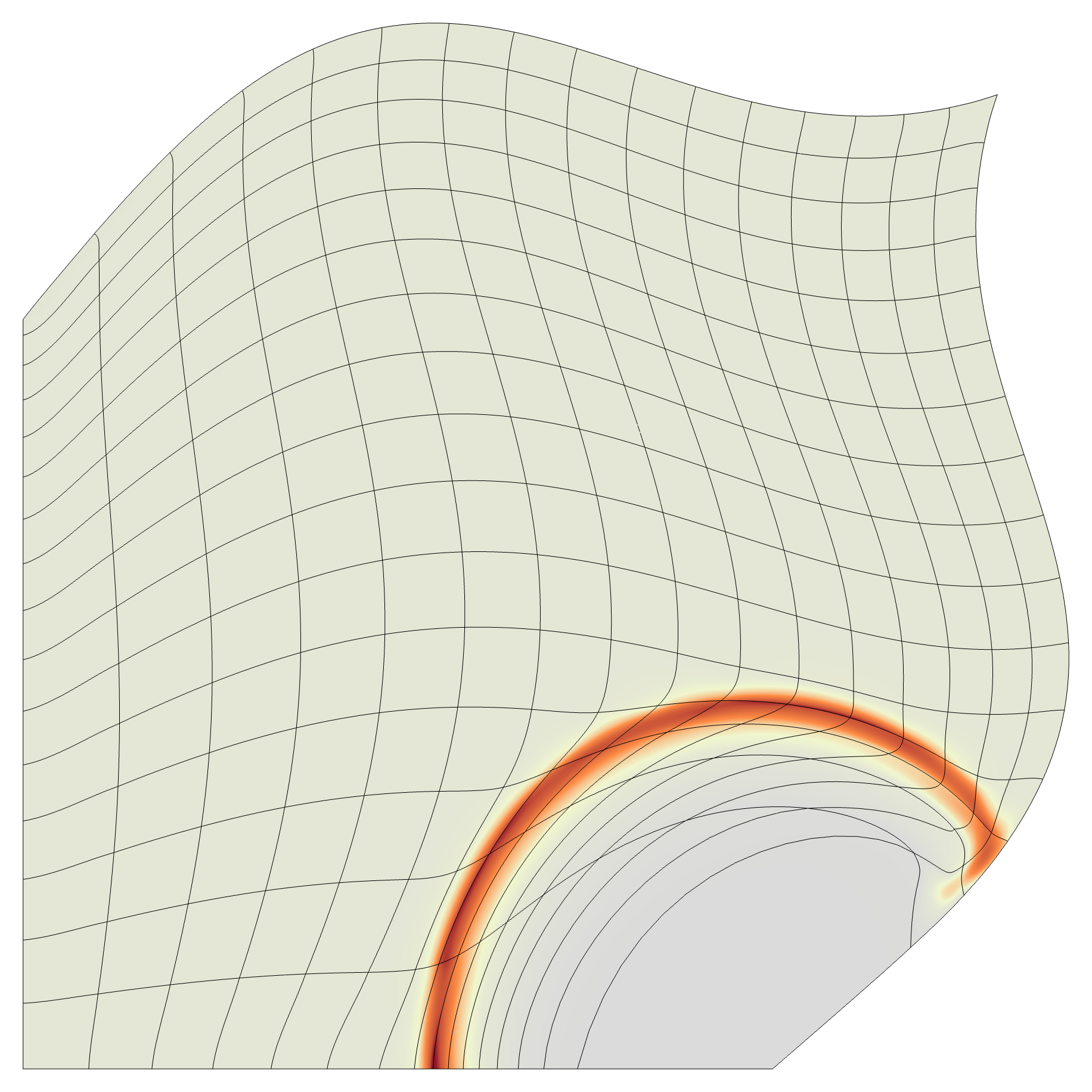}}
  & \raisebox{-.5\height}{\includegraphics[width=0.20\textwidth]{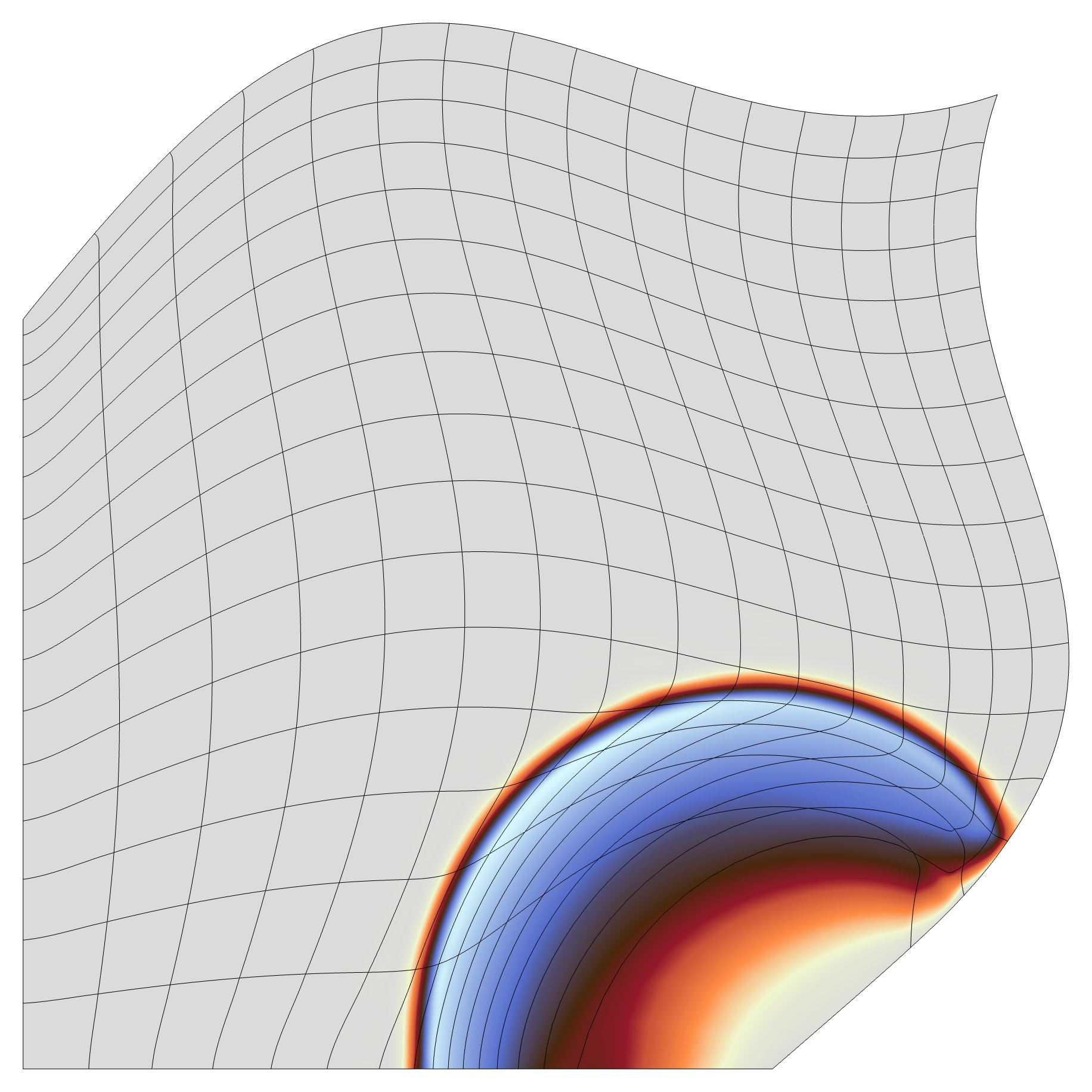}}
  & \raisebox{-.5\height}{\includegraphics[width=0.20\textwidth]{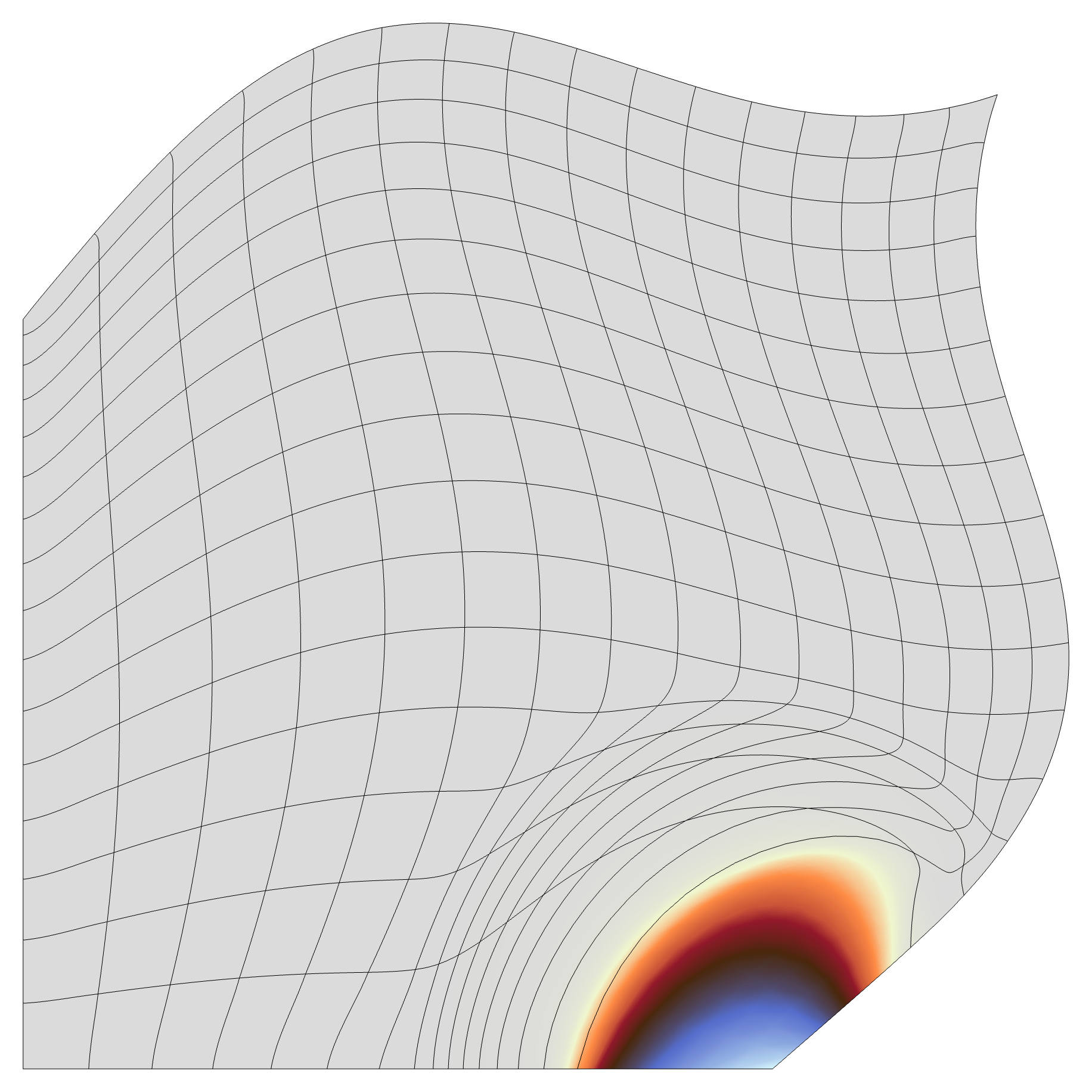}} \\

  $t=0.6$
  & \raisebox{-.5\height}{\includegraphics[width=0.20\textwidth]{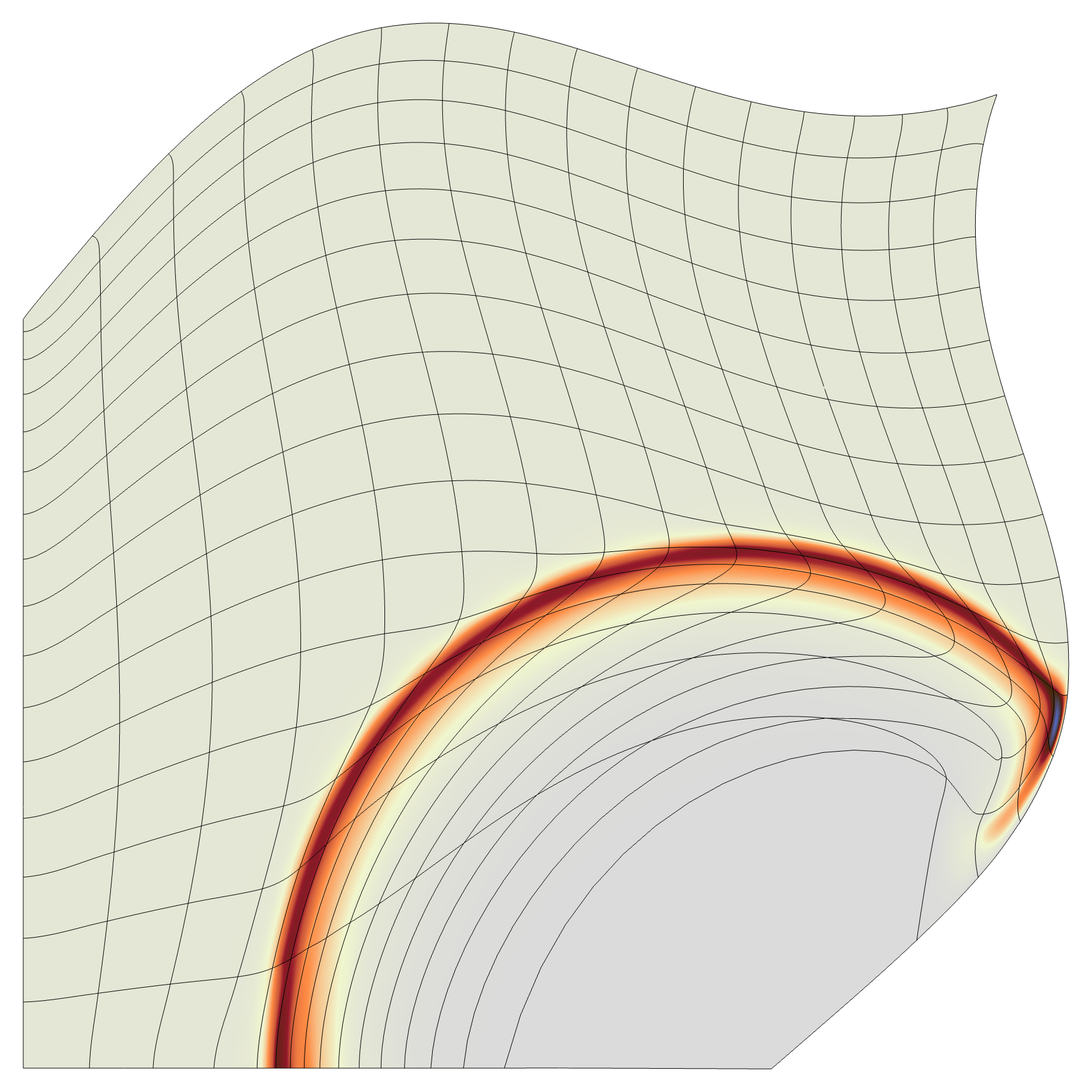}}
  & \raisebox{-.5\height}{\includegraphics[width=0.20\textwidth]{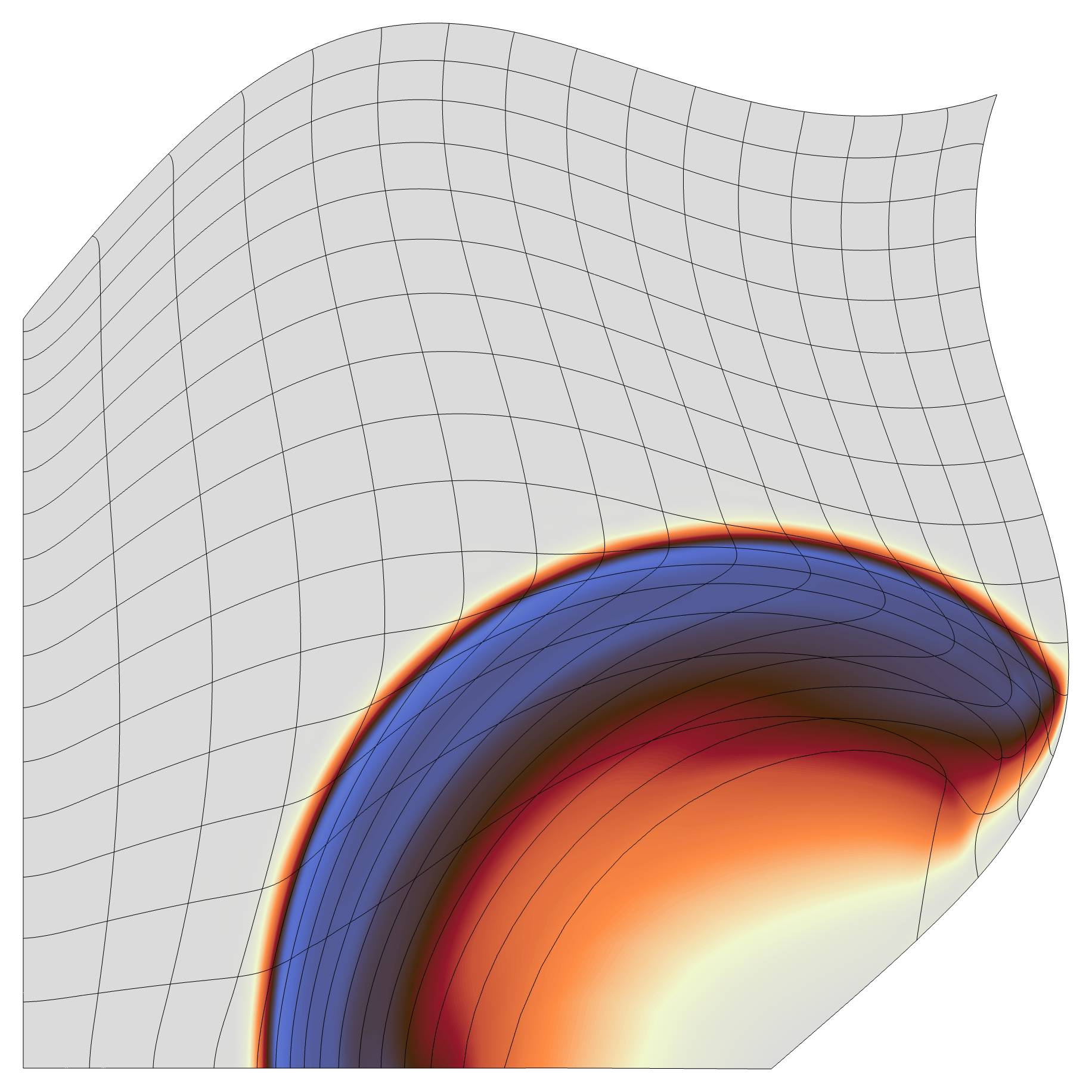}}
  & \raisebox{-.5\height}{\includegraphics[width=0.20\textwidth]{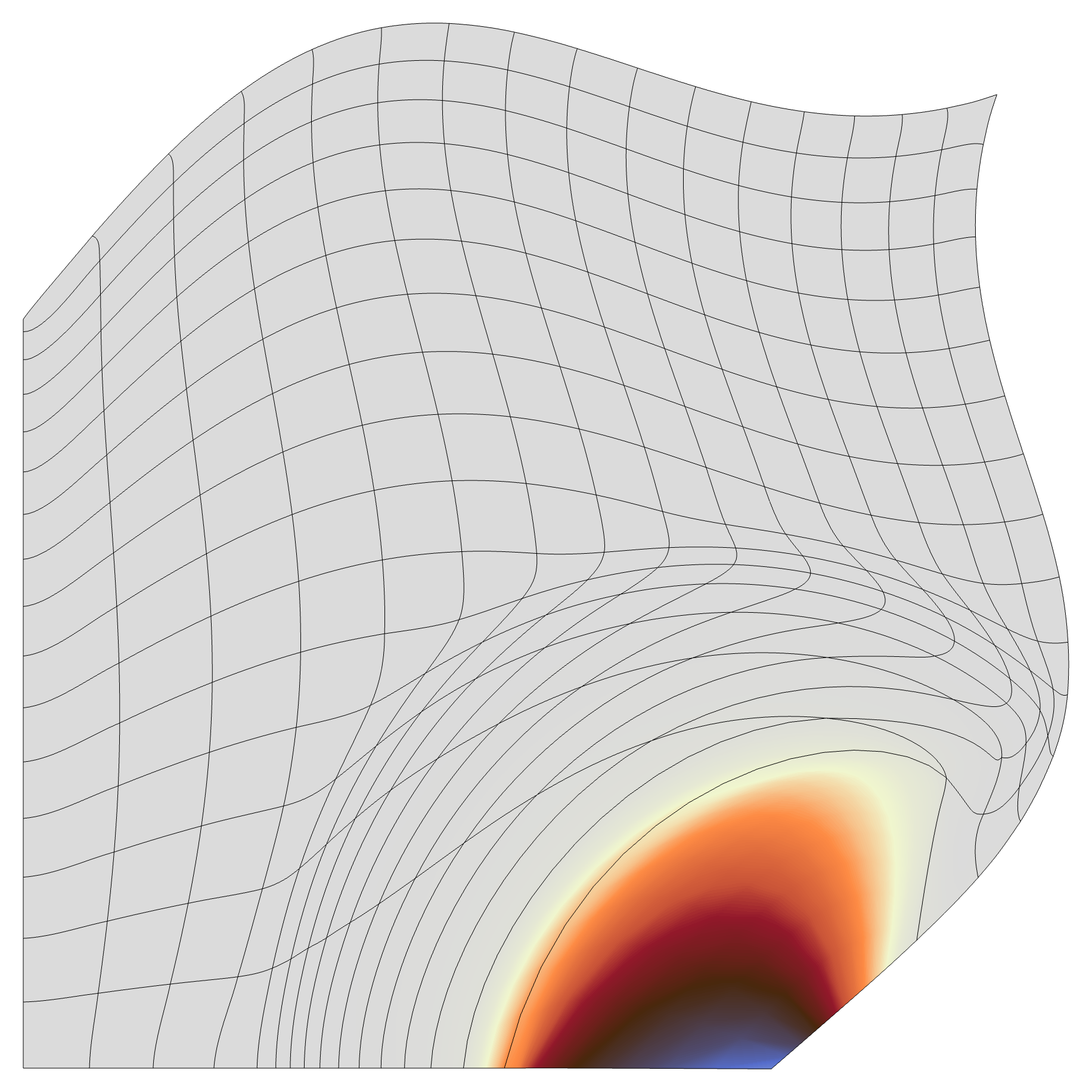}} \\

  $t=0.9$
  & \raisebox{-.5\height}{\includegraphics[width=0.20\textwidth]{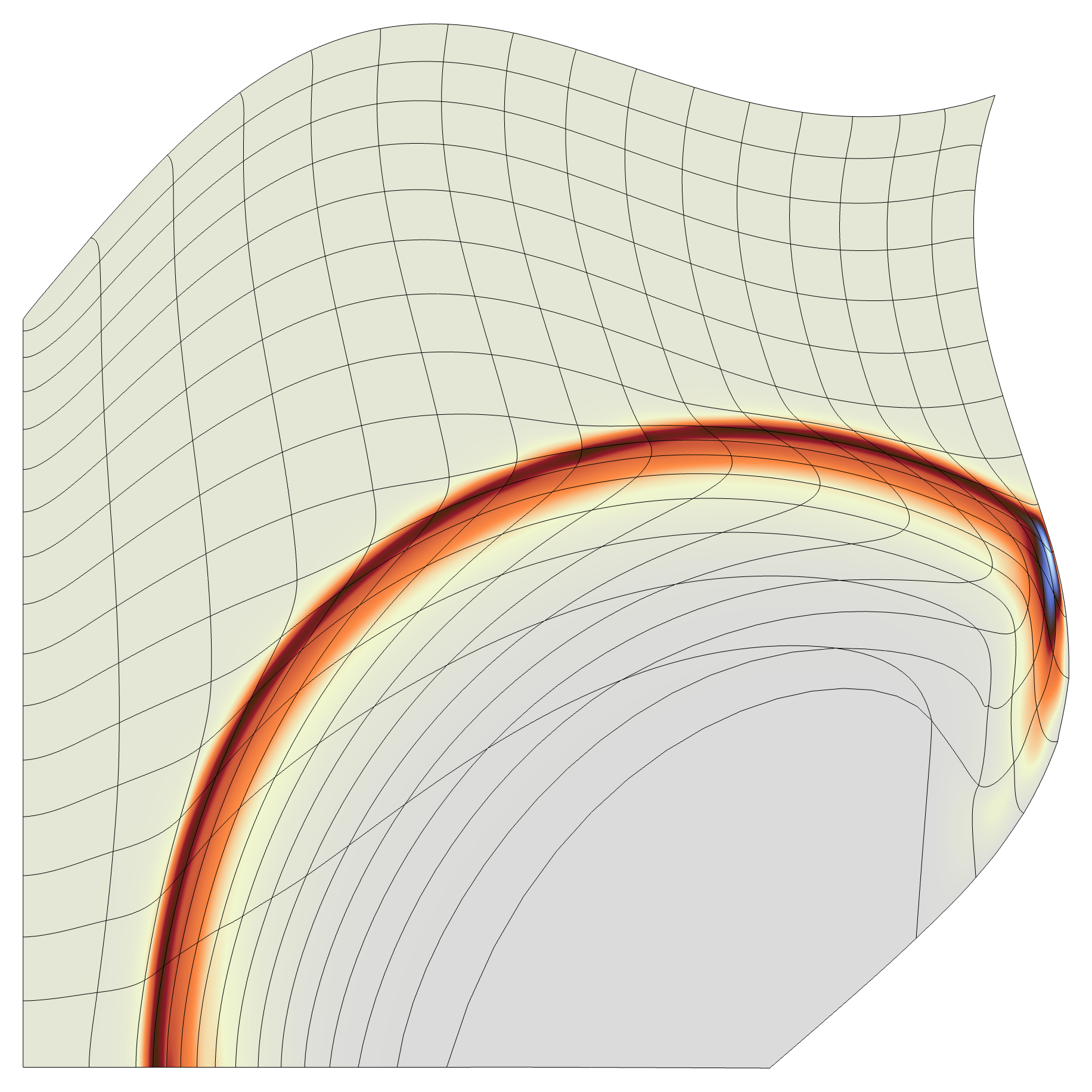}}
  & \raisebox{-.5\height}{\includegraphics[width=0.20\textwidth]{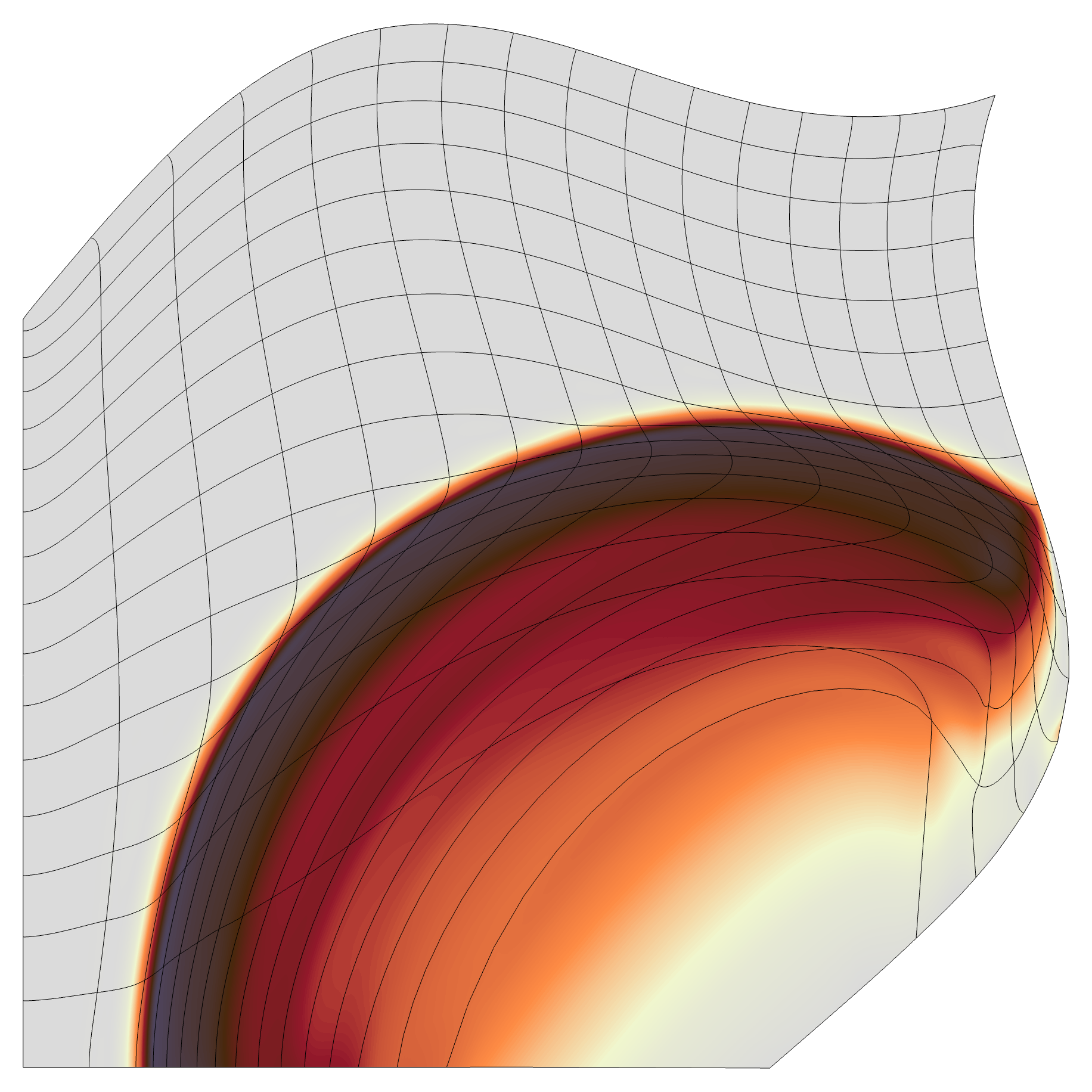}}
  & \raisebox{-.5\height}{\includegraphics[width=0.20\textwidth]{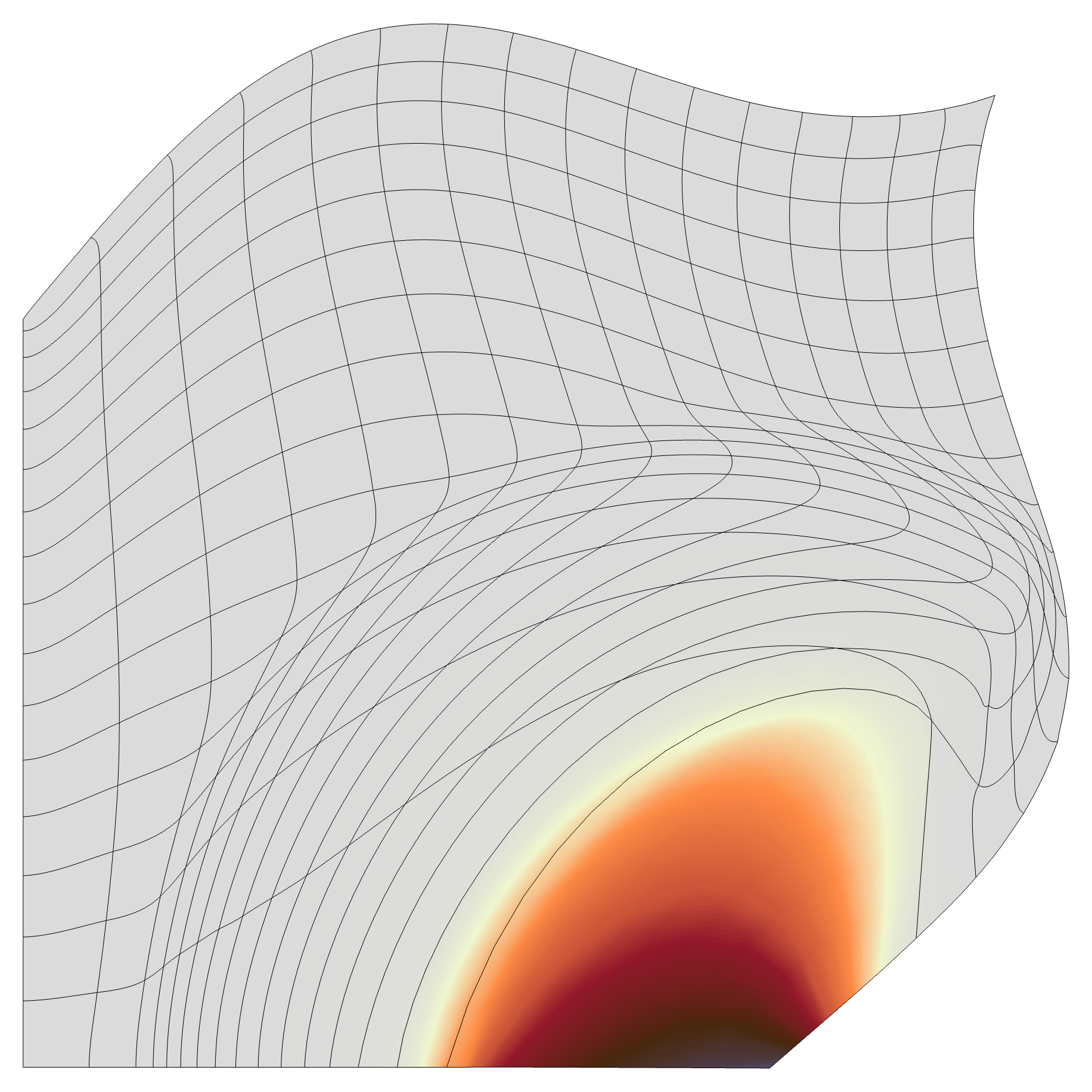}}
\end{tabular}
\caption{Plots of density, velocity, and specific internal energy for
         the Sedov blast on a 2D square-like domain with curved boundaries.}
\label{fig:sine_2D}

\vspace{0.5em}

\centering
\begin{tabular}{c c c c}
  & Density & Velocity & Specific Internal Energy \\

  $t=0.30$
  & \raisebox{-.5\height}{\includegraphics[width=0.20\textwidth]{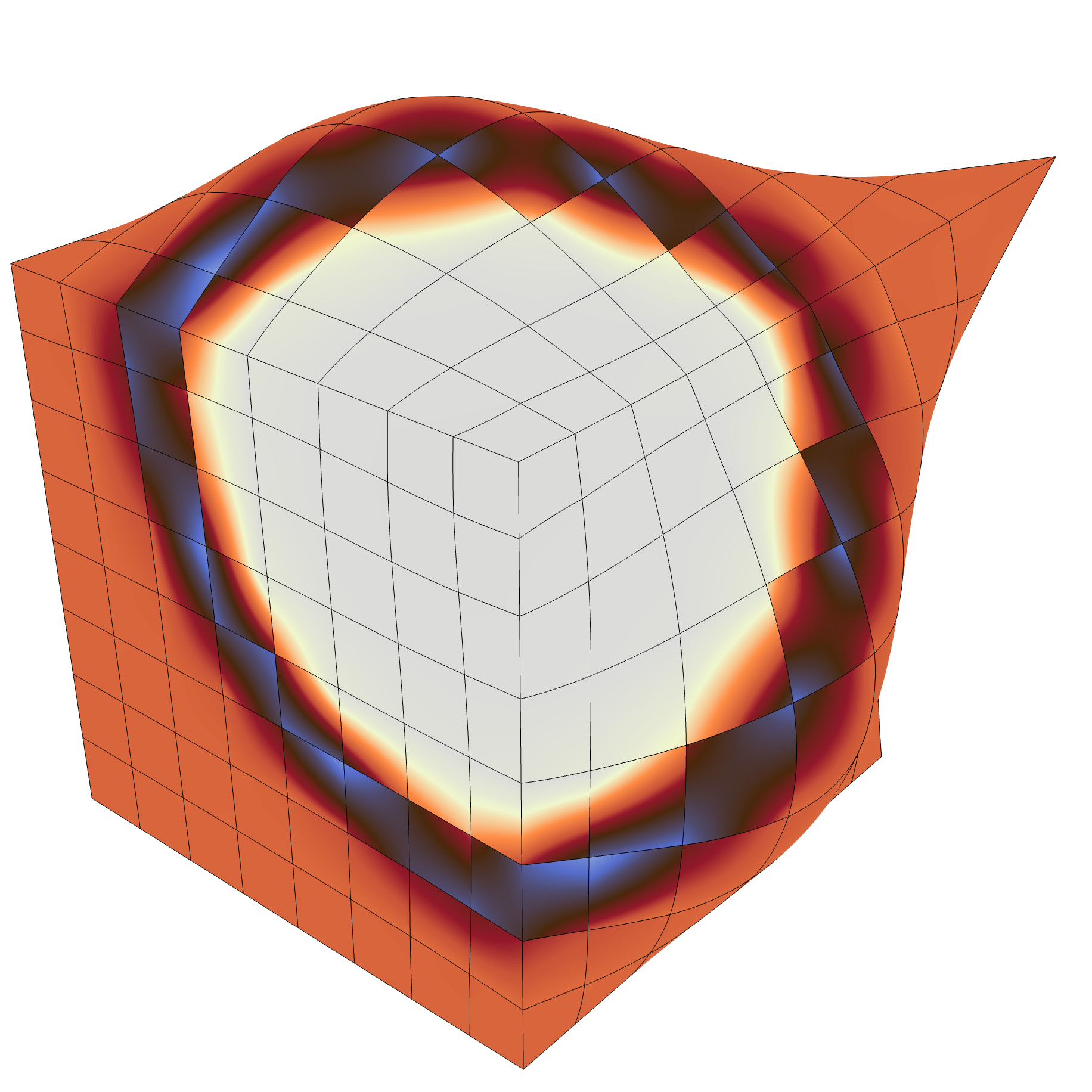}}
  & \raisebox{-.5\height}{\includegraphics[width=0.20\textwidth]{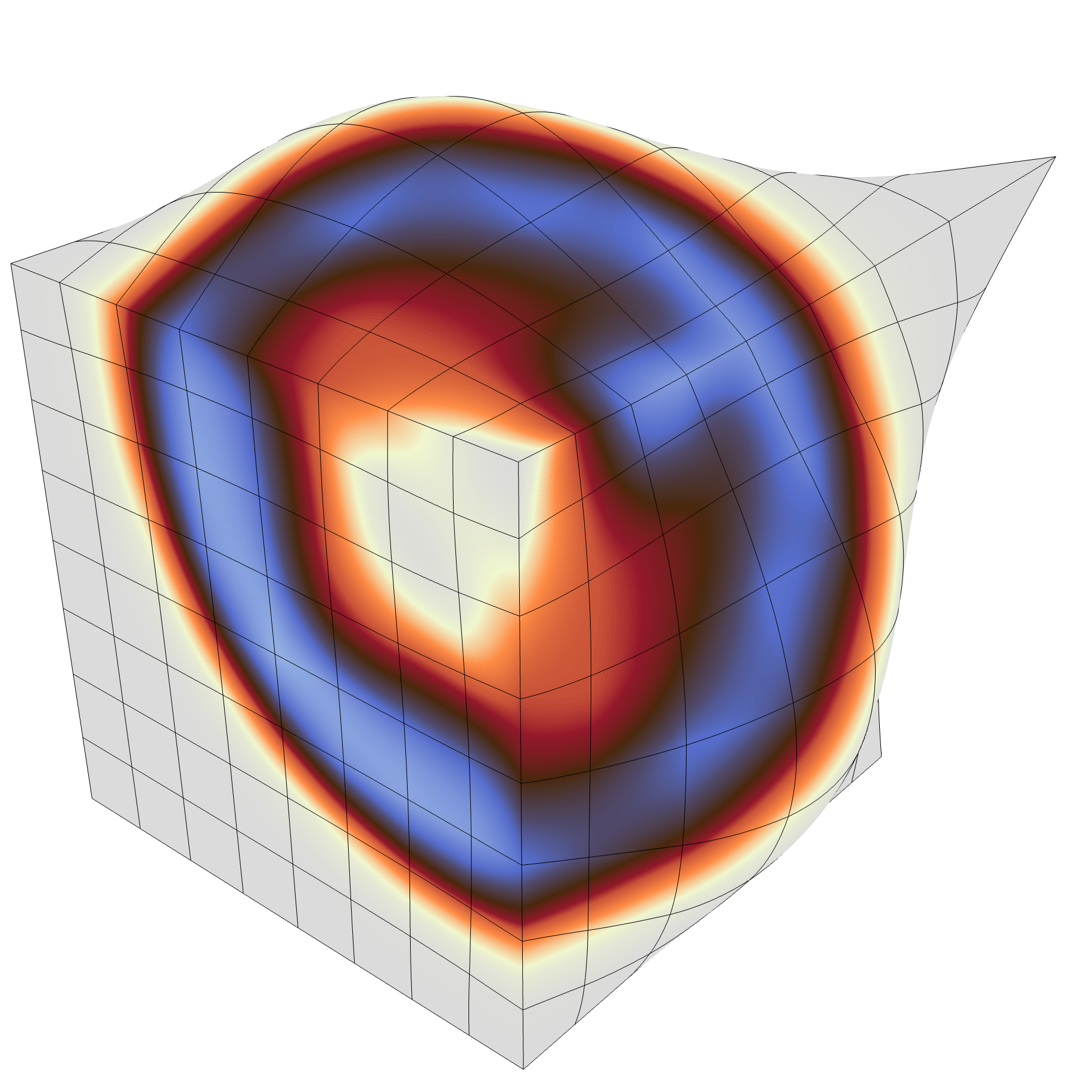}}
  & \raisebox{-.5\height}{\includegraphics[width=0.20\textwidth]{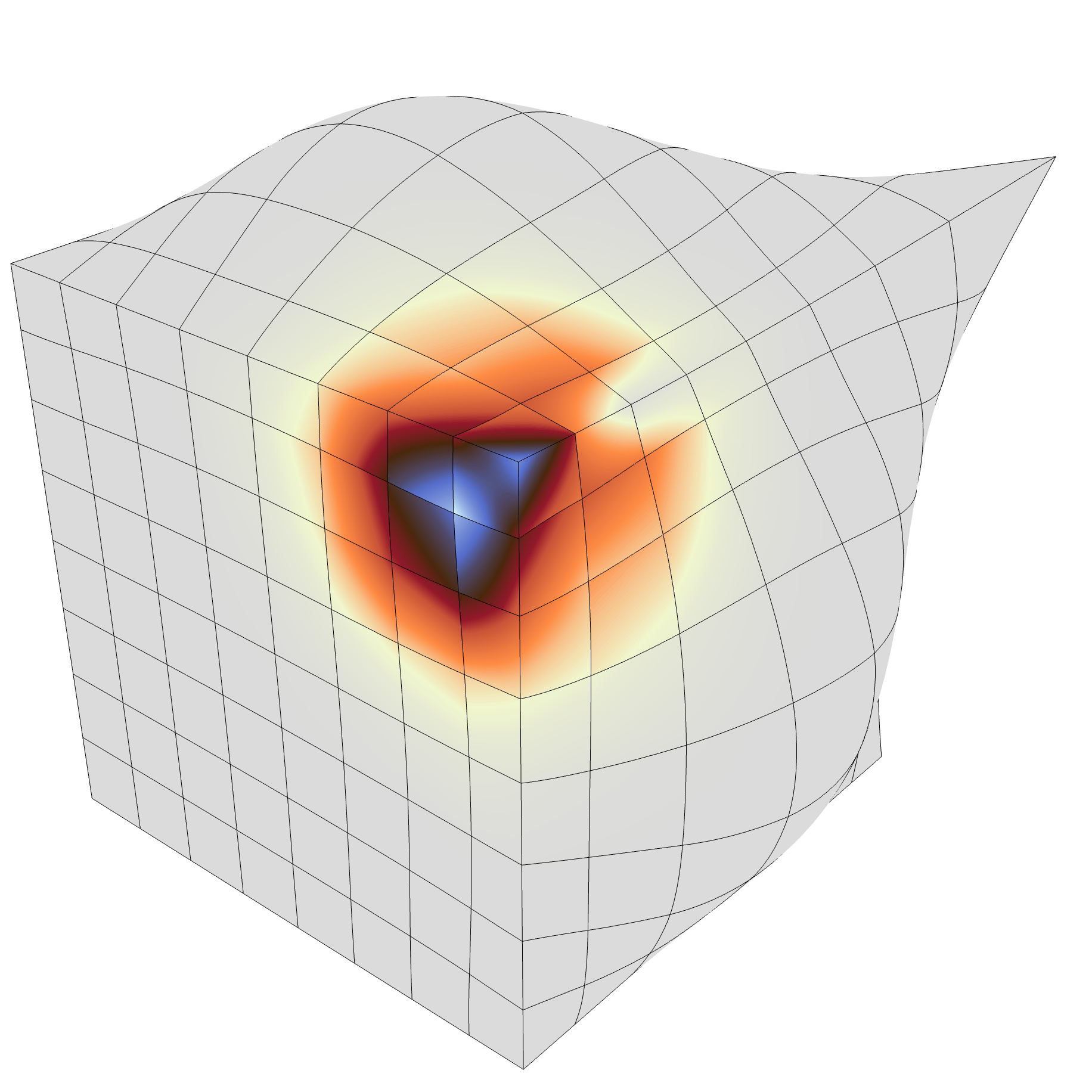}} \\

  $t=0.60$
  & \raisebox{-.5\height}{\includegraphics[width=0.20\textwidth]{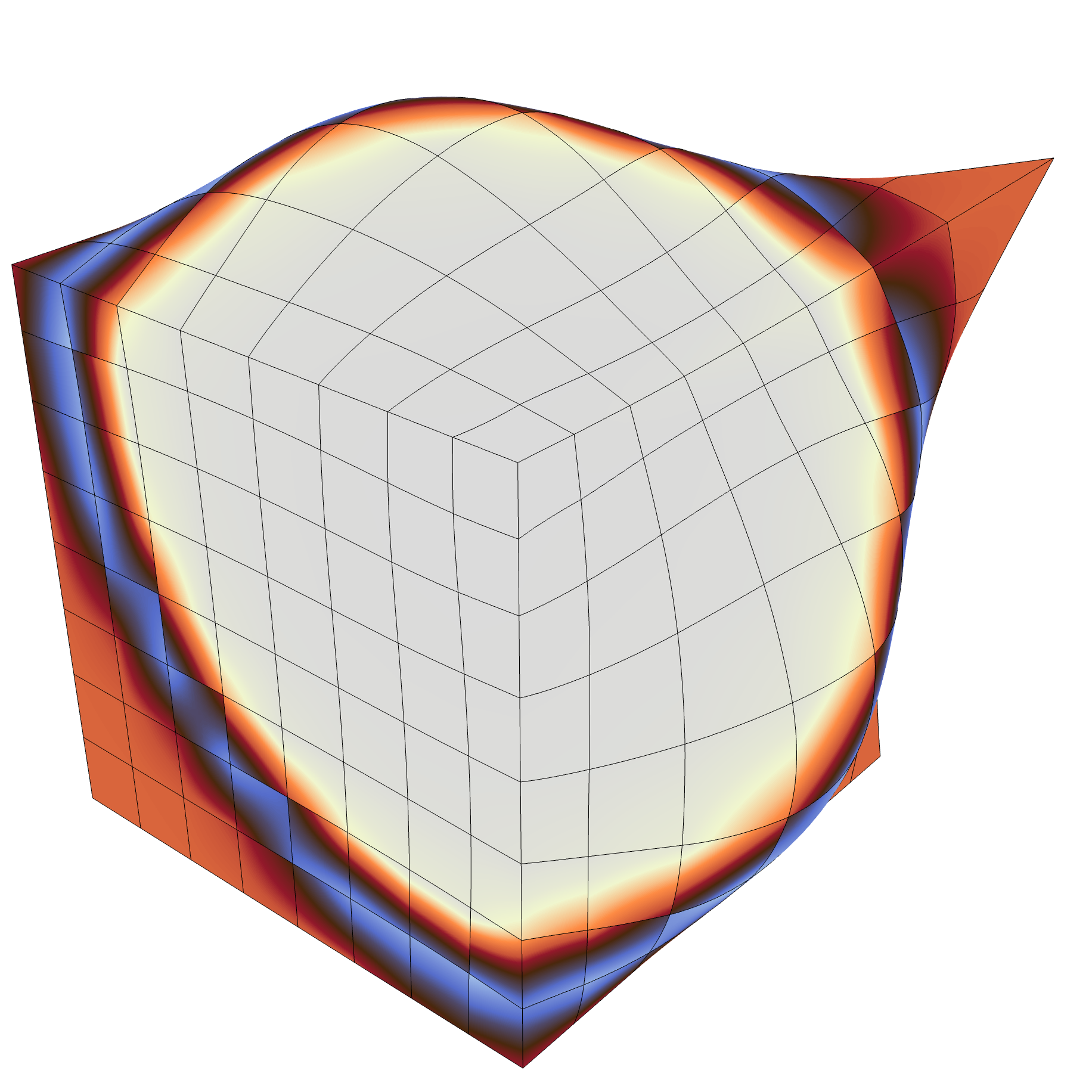}}
  & \raisebox{-.5\height}{\includegraphics[width=0.20\textwidth]{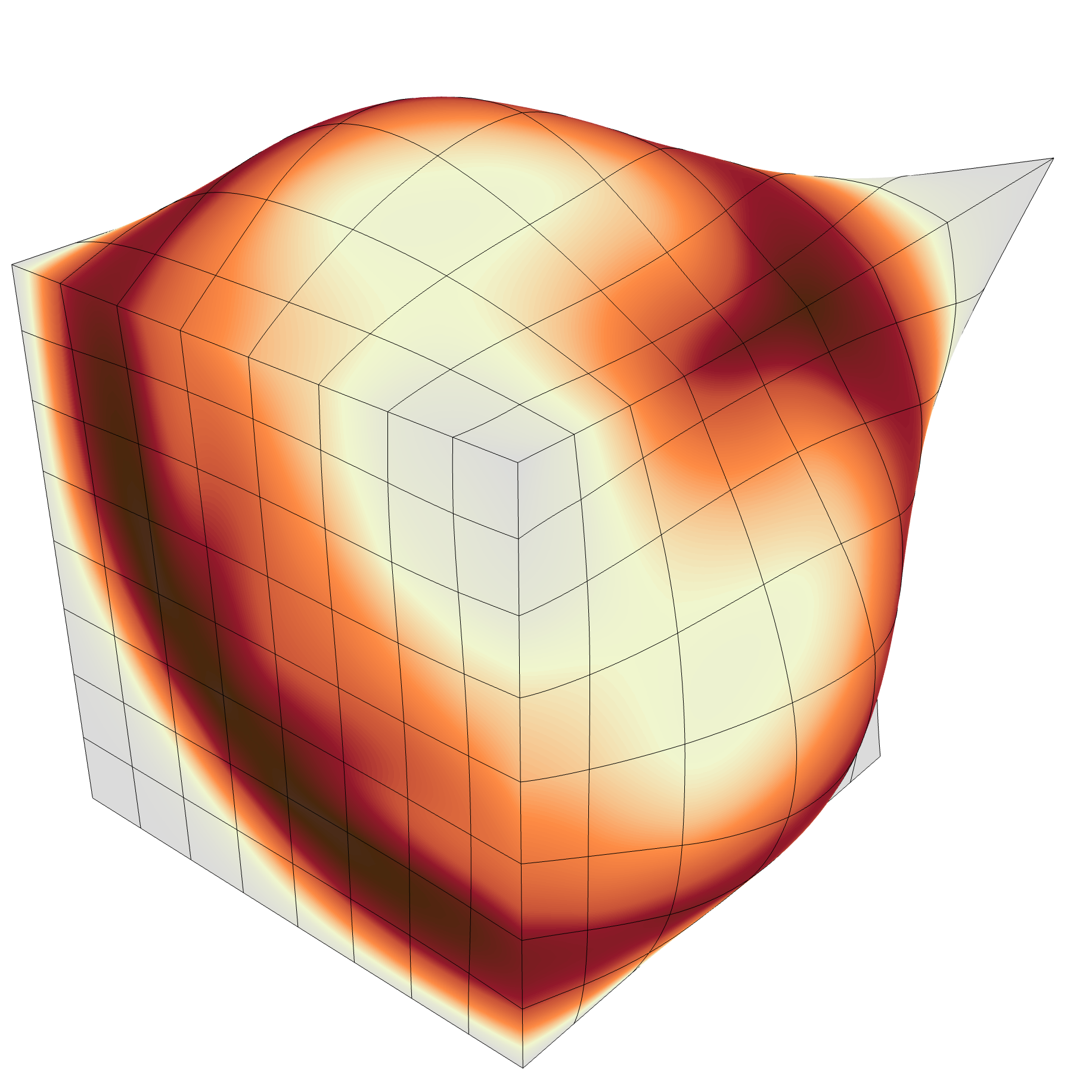}}
  & \raisebox{-.5\height}{\includegraphics[width=0.20\textwidth]{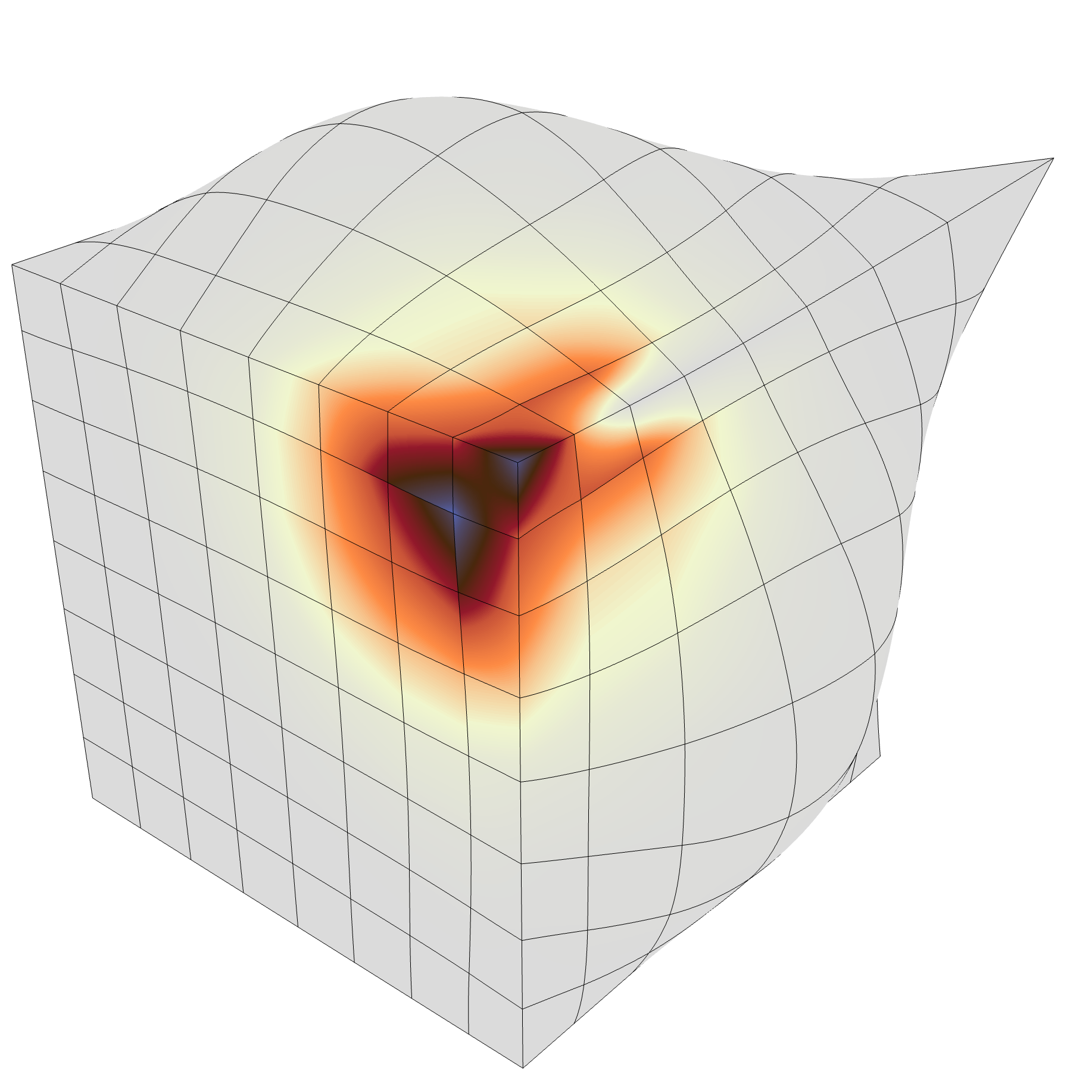}}
\end{tabular}
\caption{Plots of density, velocity, and specific internal energy for
         the Sedov blast on a cube with several curved faces.}
\label{fig:sine_3D}
\end{figure}

We next increase the geometric difficulty by considering domains with
boundaries that have pronounced curvature.
The 2D domain is given by the transformation
\[
\begin{bmatrix}
    x \\ y
\end{bmatrix} \mapsto \begin{bmatrix}
    x + 0.2x \sin\left((\frac{3\pi}{2} y\right) + \frac{1}{2}xy \\[0.3em]
    y + 0.2y \sin\left((\frac{3\pi}{2} x\right) + \frac{1}{2}xy
\end{bmatrix}
\]
applied to the unit square. This gives the following parameterizations
for the nonlinear boundaries with $t \in [0,1]$:
\begin{center}
\begin{tabular}{l l l}
    Top boundary: & & Right boundary: \\
    $x(t) = 1.3t$ & &
    $x(t) = 1 + 0.2 \sin\left(\frac{3\pi}{2} t\right) + 0.5t$ \\
    $y(t) = 1 + 0.2 \sin\left(\frac{3\pi}{2} t\right) + 0.5t$ & &
    $y(t) = 1.3t$
\end{tabular}
\end{center}
The blast is initialized at the bottom-right corner of the mesh. The evolution
shown in Figure \ref{fig:sine_2D} shows the expected behavior as the shock
reaches and reflects from the curved boundary: the motion follows the wall
tangentially, without visible boundary-induced mesh deformation or spurious
solution oscillations.

\begin{figure}[pos=htbp]
\centering
\begin{tabular}{c c c}
  \includegraphics[width=0.2\textwidth]{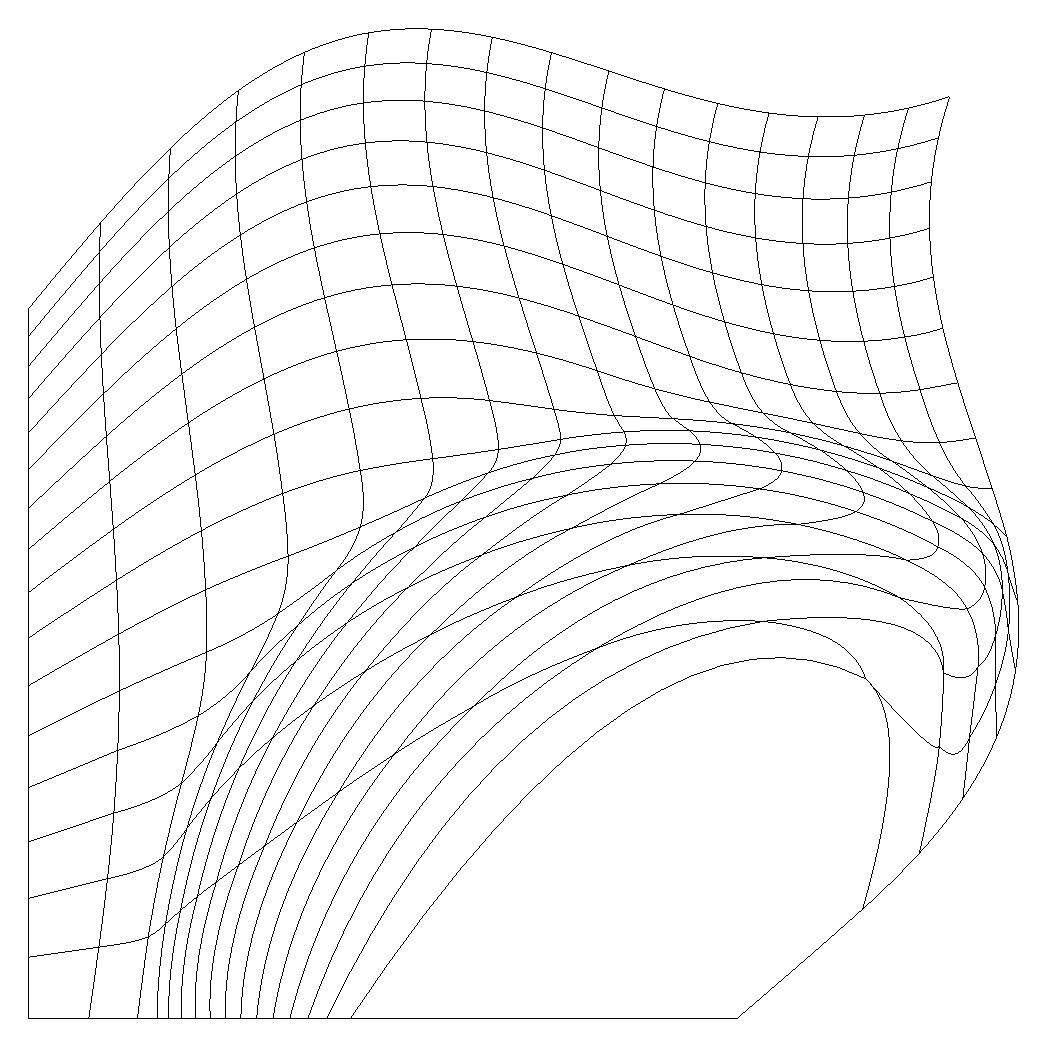} &
  \includegraphics[width=0.2\textwidth]{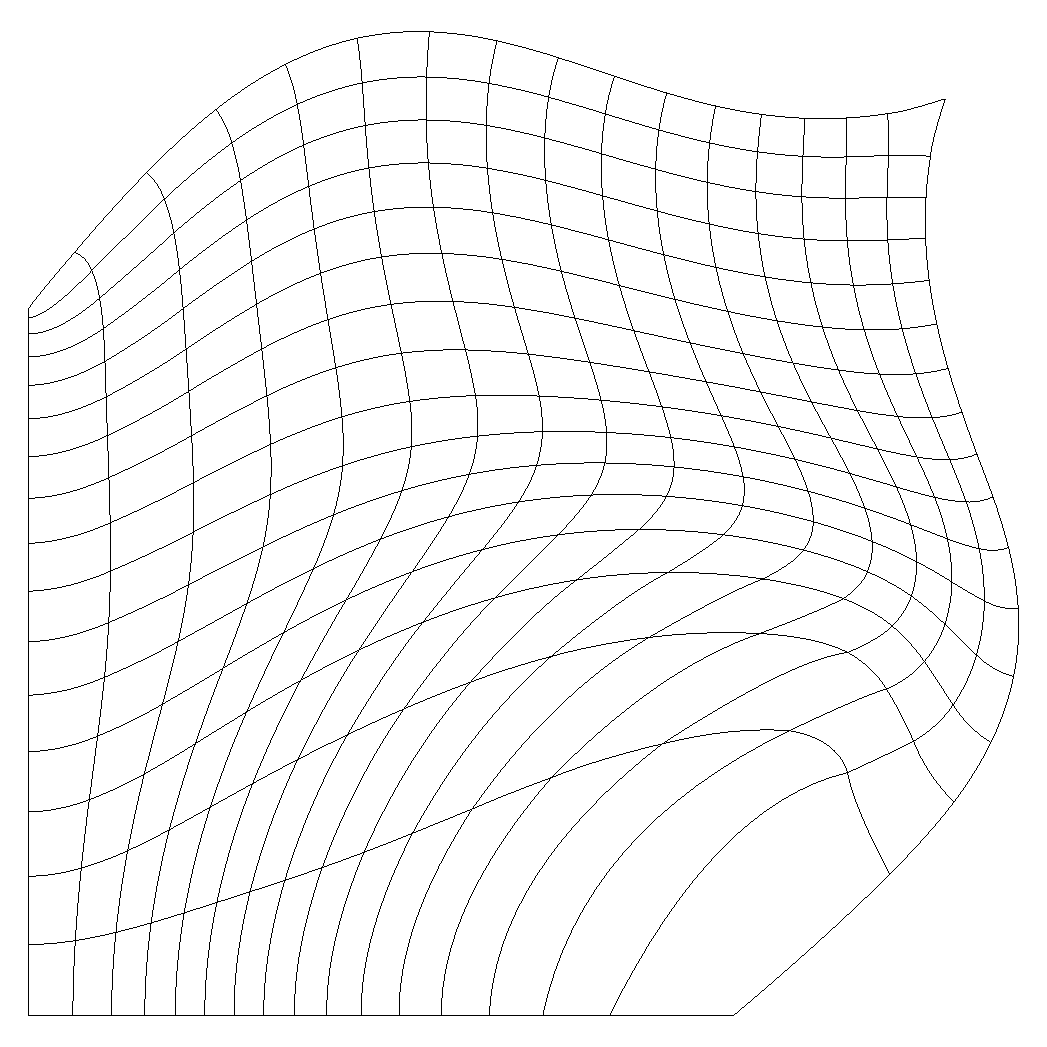} &
  \includegraphics[width=0.2\textwidth]{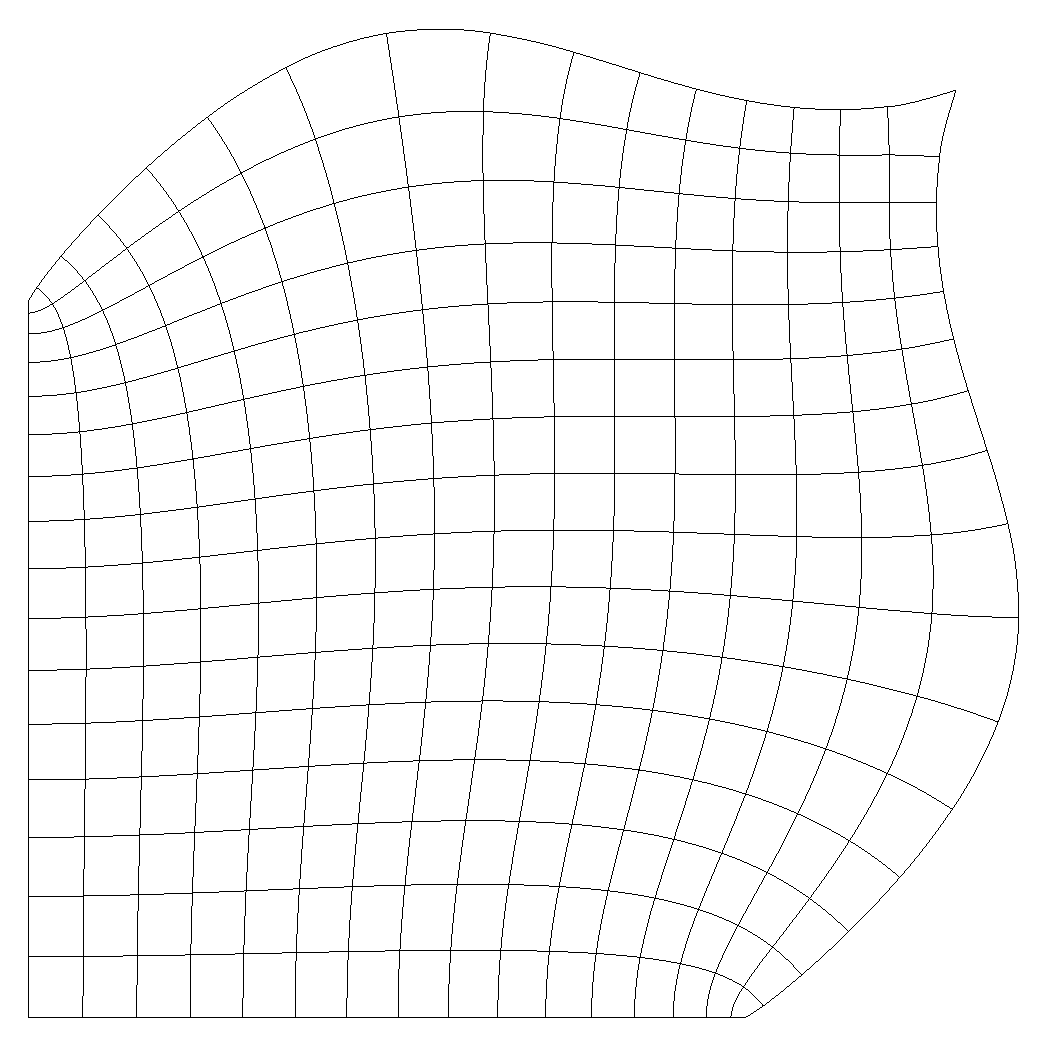}
\end{tabular}
\caption{The 2D mesh with sine boundaries before the ALE remesh phase (left),
         a mesh optimized with limiting distance $\delta = 0.01$ (middle),
         and a fully optimized mesh (right) for a transformed square with
         curved edges.}
\label{fig:sine_tmop}
\end{figure}

Figure \ref{fig:sine_tmop} shows the corresponding behavior of the mesh
optimizer. As in the linear case, the limited optimization improves the mesh
while preserving the Lagrangian displacement near the boundary, and the
unlimited optimization shows the motion preferred by the quality functional.

The 3D domain is given by the transformation
\[
\begin{bmatrix}
    x \\ y \\ z
\end{bmatrix} \mapsto \begin{bmatrix}
    x + 0.2\sin(0.5 \pi x)\sin(1.1 \pi y)\sin(1.1 \pi z)
      + 0.1 xyz \\
    y + 0.2\sin(1.1 \pi x)\sin(0.5 \pi y)\sin(1.1 \pi z)
      + 0.1 xyz \\
    z + 0.2\sin(1.1 \pi x)\sin(1.1 \pi y)\sin(0.5 \pi z)
      + 0.1 xyz
\end{bmatrix}
\]
applied to the unit cube.
This gives the following parameterizations with $u,v,t \in [0,1]$:
\begin{center}
\begin{tabular}{l l l}
    Top face: & & Edge between top and front faces: \\
    $x(u,v) = u + 0.2 \sin(0.5 \pi u)\sin(1.1 \pi v)\sin(1.1 \pi)
      + 0.1uv$ & &
    $x(t) = t + 0.2 \sin^2(1.1 \pi)\sin(0.5 \pi t) + 0.1t$ \\
    $y(u,v) = v + 0.2 \sin(1.1 \pi u)\sin(0.5 \pi v)\sin(1.1 \pi)
      + 0.1uv$ & &
    $y(t) = 1 + 0.2 \sin(1.1 \pi t) + 0.1t$ \\
    $z(u,v) = 1 + 0.2 \sin(1.1 \pi u)\sin(1.1 \pi v) + 0.1uv$ & &
    $z(t) = 1 + 0.2 \sin(1.1 \pi t) + 0.1t$
\end{tabular}
\end{center}
The other nonlinear faces and edges are defined similarly, and the blast is
initialized at a corner of two curved faces and one planar face.
The 3D evolution in Figure \ref{fig:sine_3D} again shows the expected
interaction with the curved walls: after the shock reaches the boundary,
the solution slides along the curved features without visible
mesh deformation or boundary oscillations.

\subsection{Annulus and Solid Torus}

The hardest tests considered here use domains where the shock repeatedly
encounters curved boundaries throughout the simulation. The 2D annulus has an
inner radius $r_0 = 0.4$ and outer radius $r_1 = 1.0$. The boundaries are
thus parameterized as $x(t) = r_i \cos(t)$, $y(t) = r_i \sin(t)$ for
$i \in \{0,1\}$ and $t \in (-\pi, \pi]$. The limiting distance for the
remesh is $\delta = 0.3$, which allows almost unlimited optimization for this
case. The blast origin is the leftmost point of the outer circle.

\begin{figure}[pos=htbp]
\centering
\begin{tabular}{c c c}
  \includegraphics[width=0.2\textwidth]{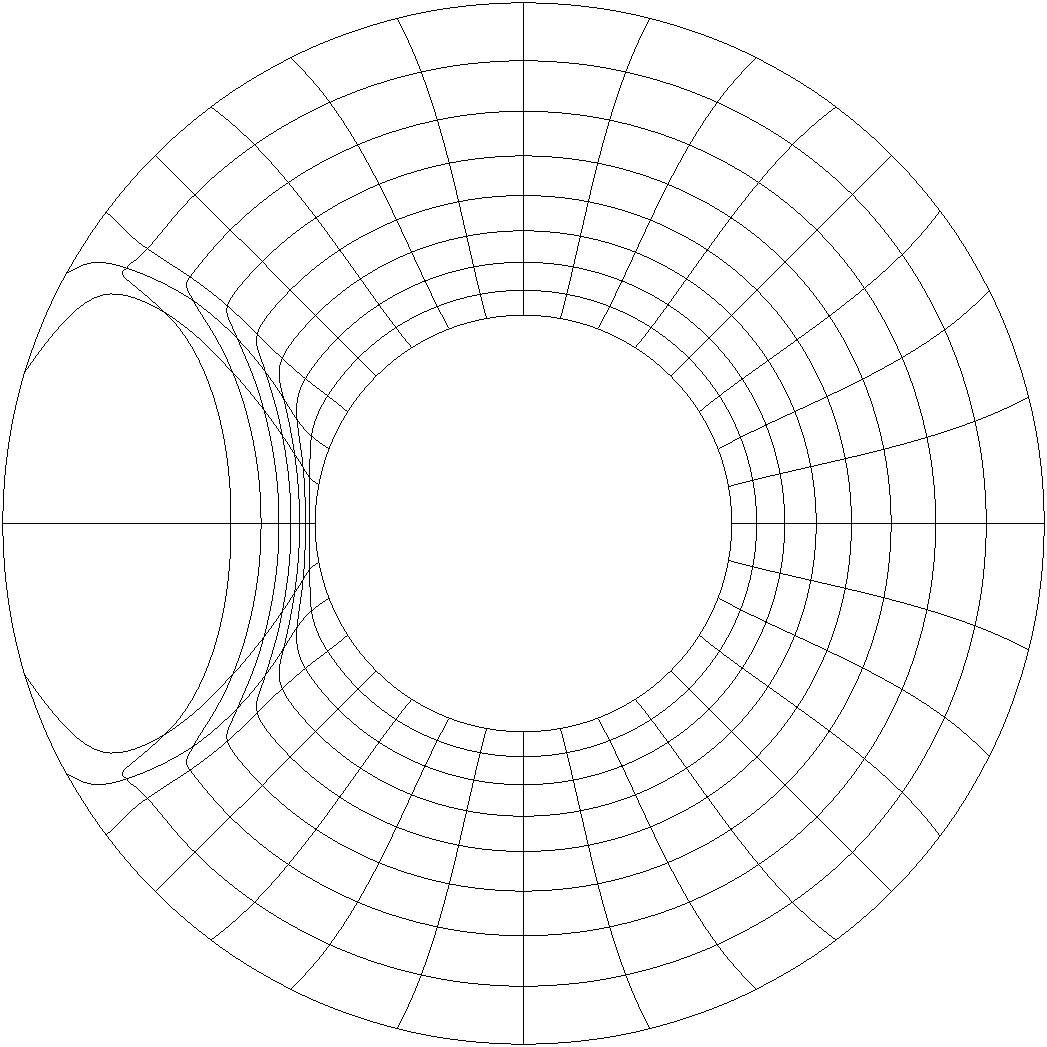} &
  \includegraphics[width=0.2\textwidth]{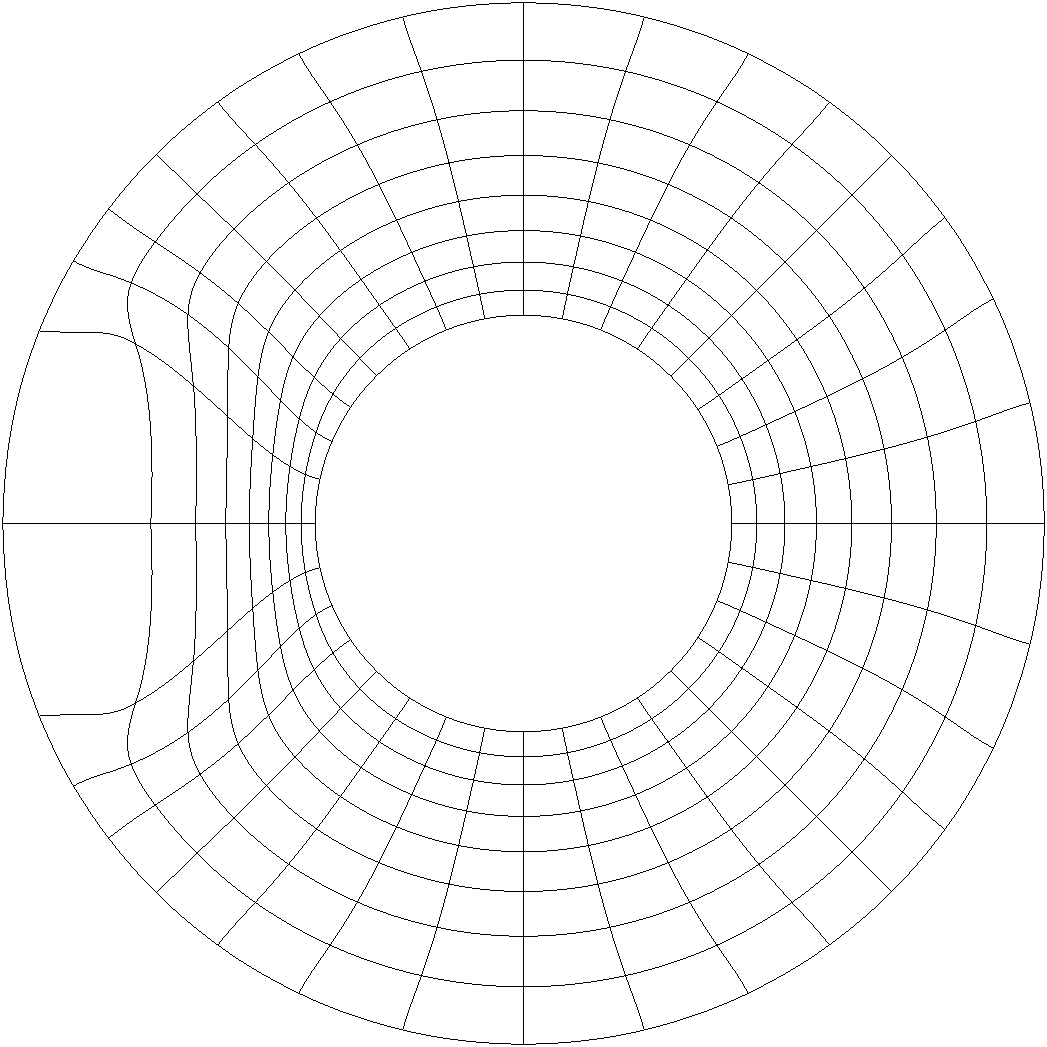} &
  \includegraphics[width=0.2\textwidth]{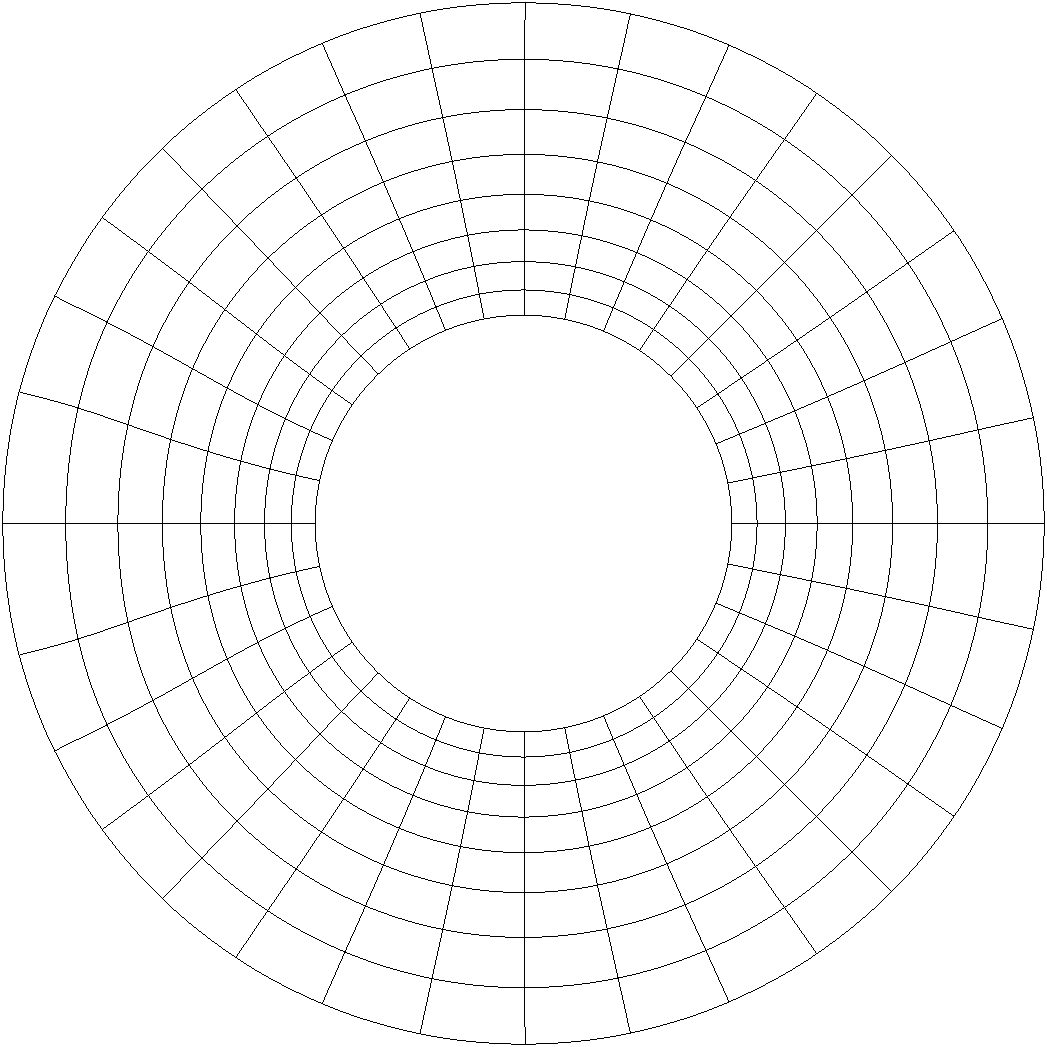}
\end{tabular}
\caption{An annular mesh before the ALE remesh phase (left), a mesh
         optimized with limiting distance $\delta = 0.05$ (middle), and a
         fully optimized mesh (right).}
\label{fig:annulus_tmop}
\end{figure}

Figure \ref{fig:annulus_tmop} shows the behavior of the mesh optimizer for
the annular domain. The limited optimization preserves the displacement near
the curved wall, while the fully optimized mesh shows the quality-driven
motion in the absence of this restriction.

\begin{figure}[pos=htbp]
\centering
\begin{tabular}{c c c c}
  & Density & Velocity & Specific Internal Energy \\

  $t=1.0$
  & \raisebox{-.5\height}{\includegraphics[width=0.20\textwidth]{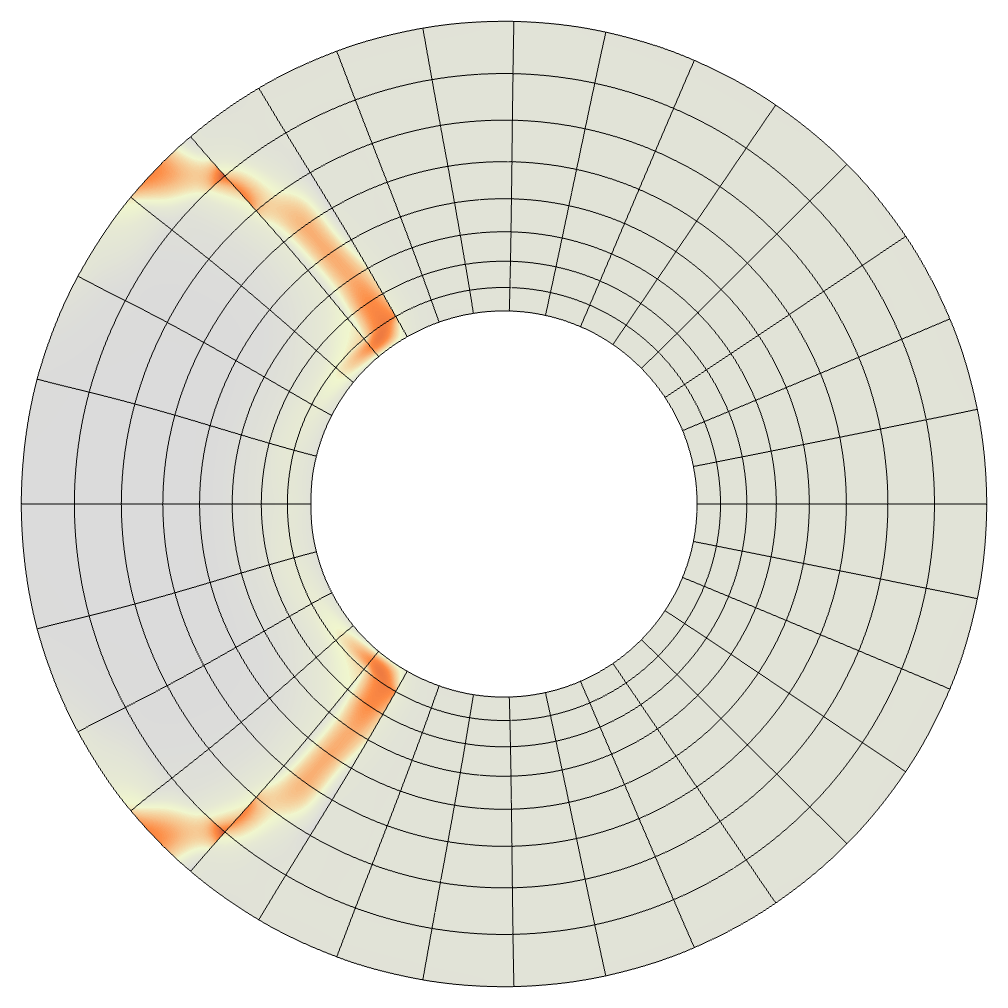}}
  & \raisebox{-.5\height}{\includegraphics[width=0.20\textwidth]{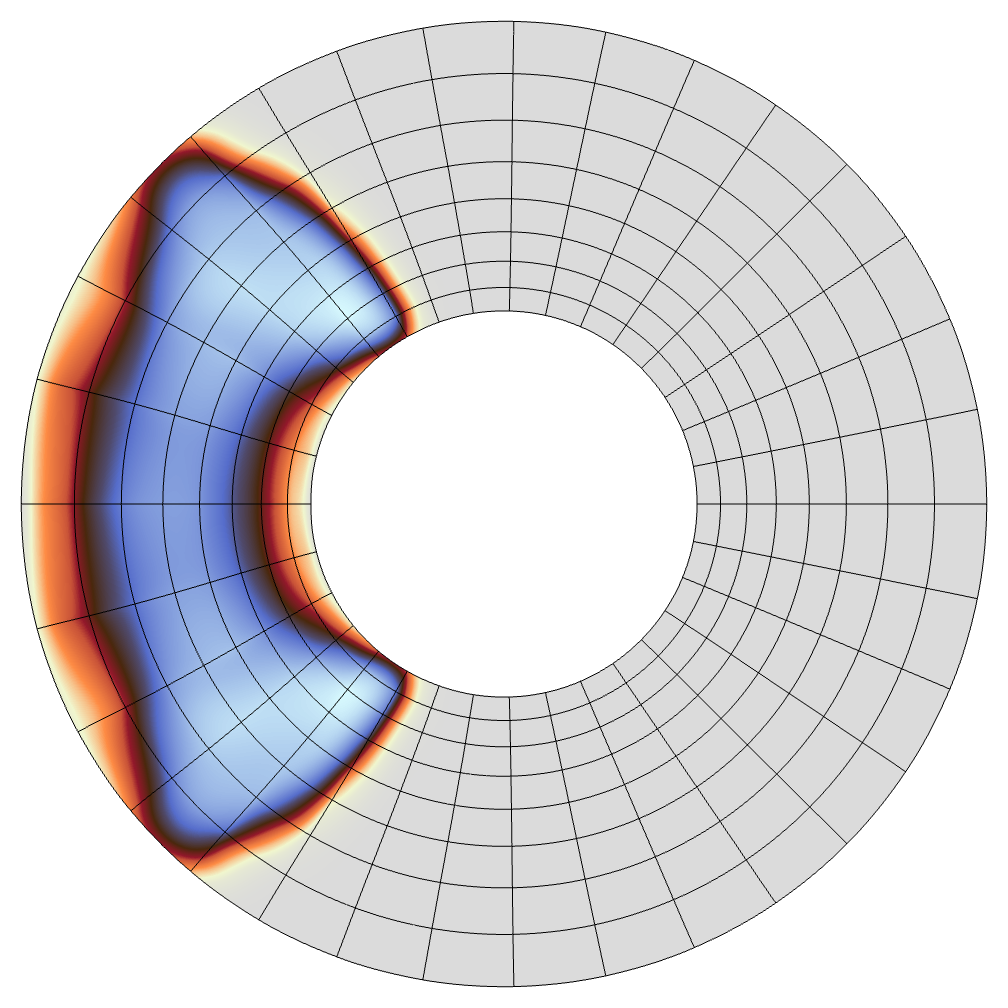}}
  & \raisebox{-.5\height}{\includegraphics[width=0.20\textwidth]{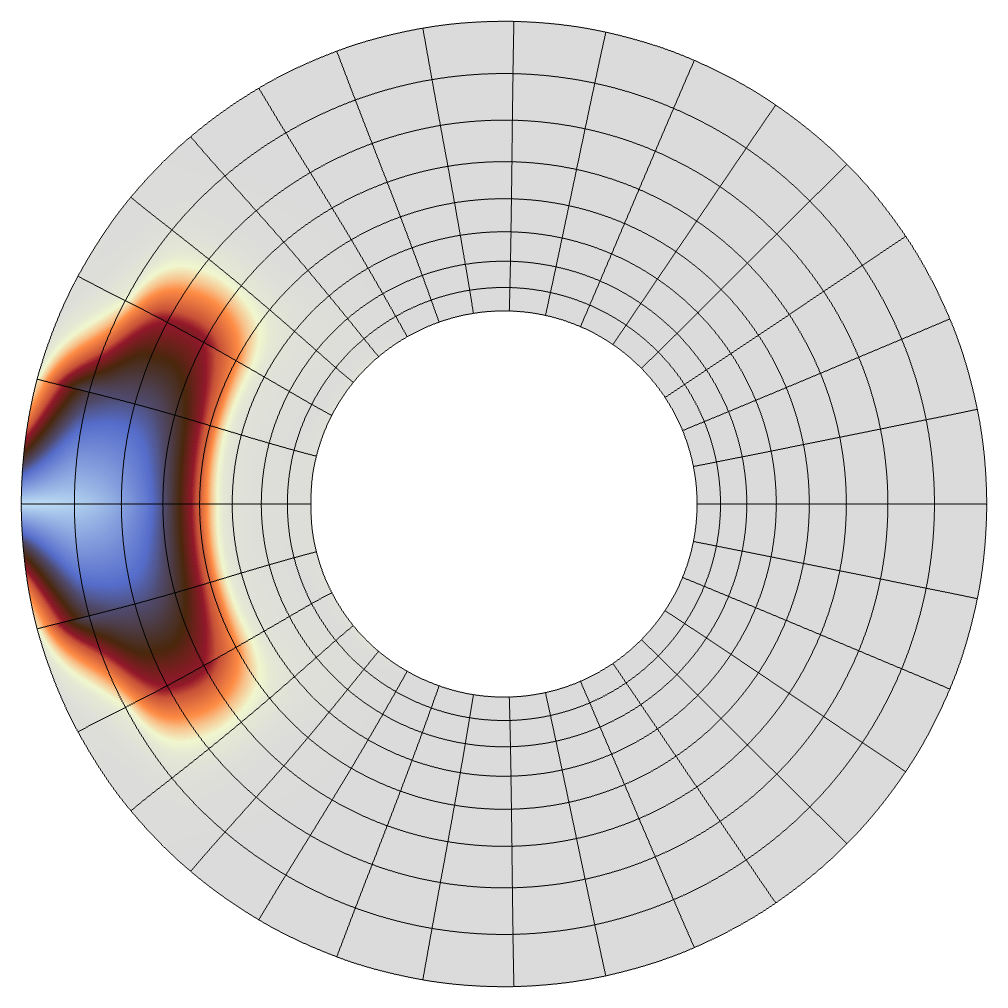}} \\

  $t=4.0$
  & \raisebox{-.5\height}{\includegraphics[width=0.20\textwidth]{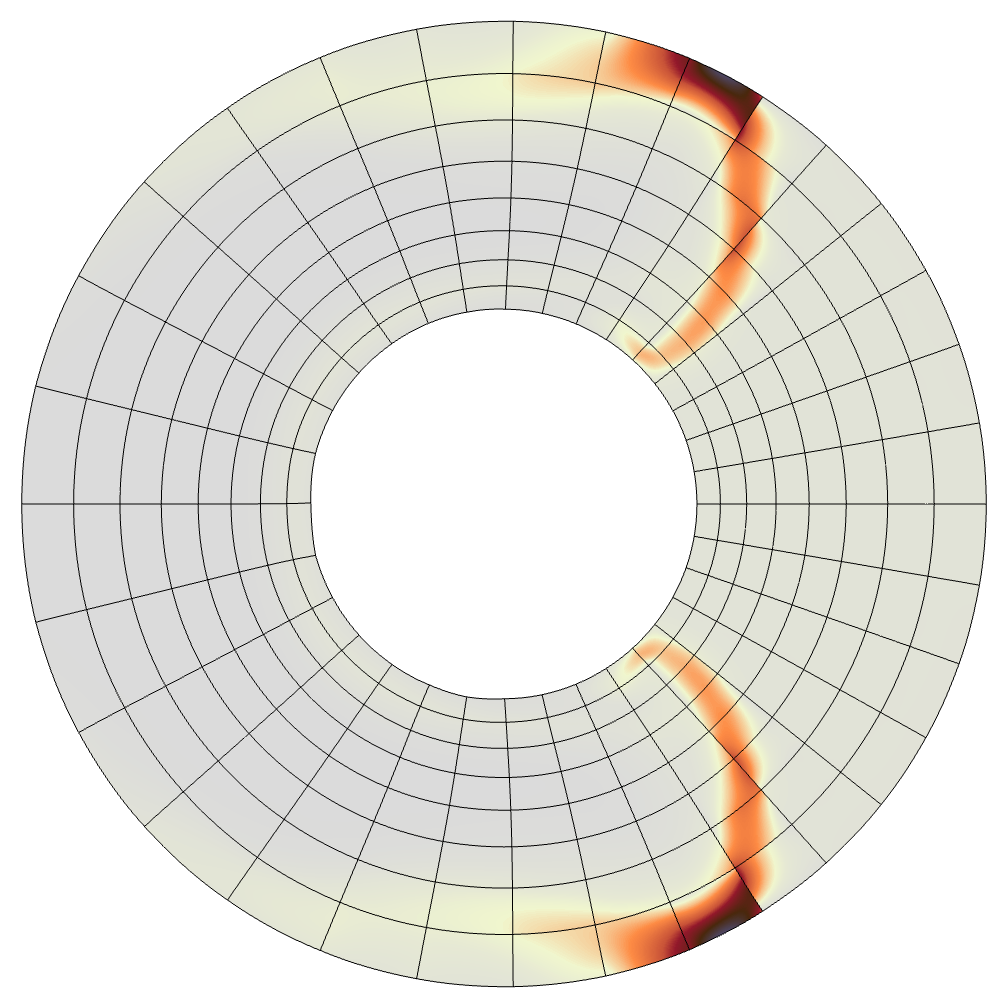}}
  & \raisebox{-.5\height}{\includegraphics[width=0.20\textwidth]{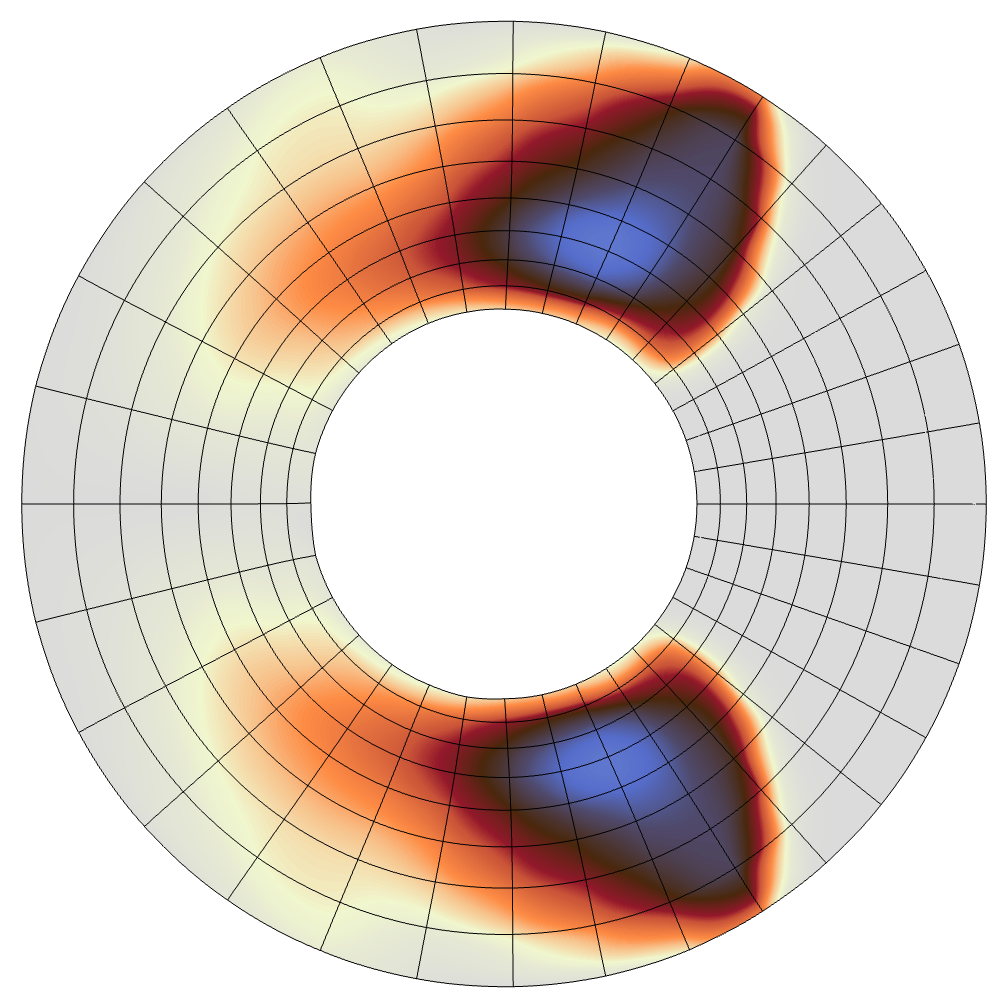}}
  & \raisebox{-.5\height}{\includegraphics[width=0.20\textwidth]{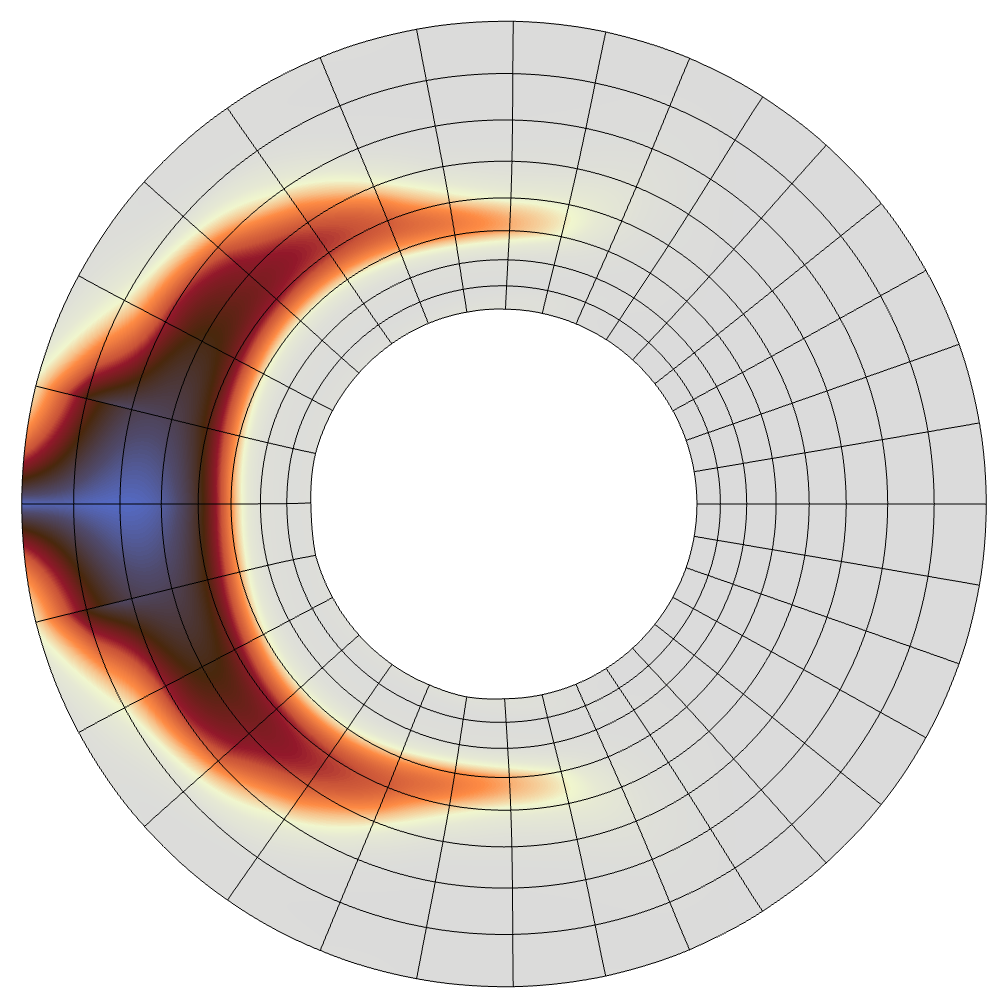}} \\

  $t=7.0$
  & \raisebox{-.5\height}{\includegraphics[width=0.20\textwidth]{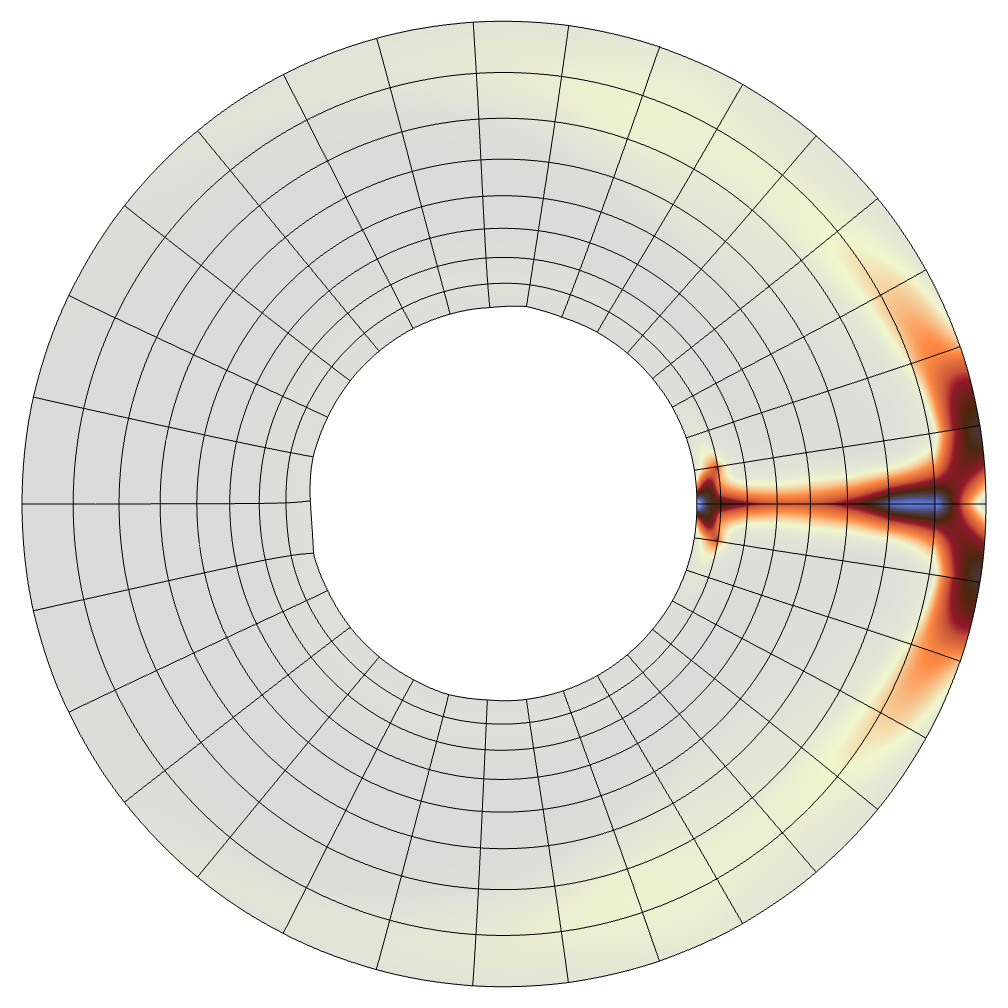}}
  & \raisebox{-.5\height}{\includegraphics[width=0.20\textwidth]{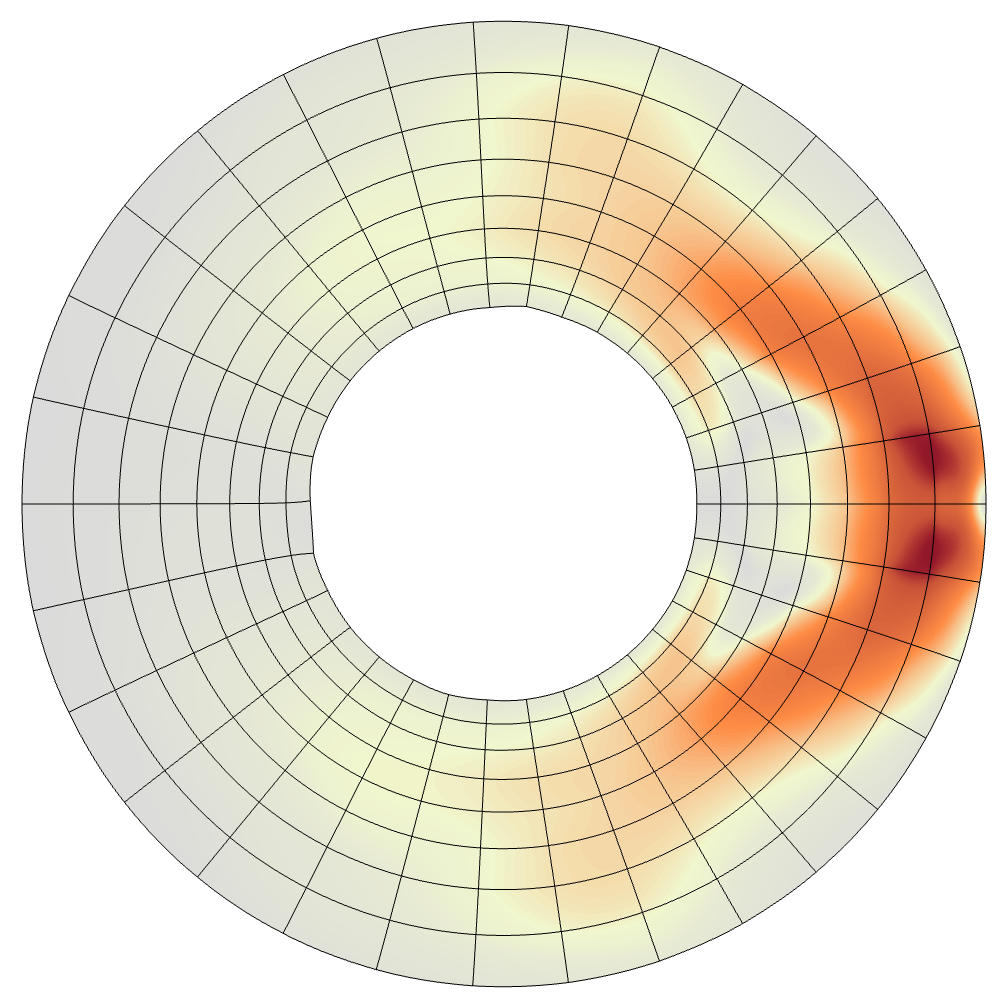}}
  & \raisebox{-.5\height}{\includegraphics[width=0.20\textwidth]{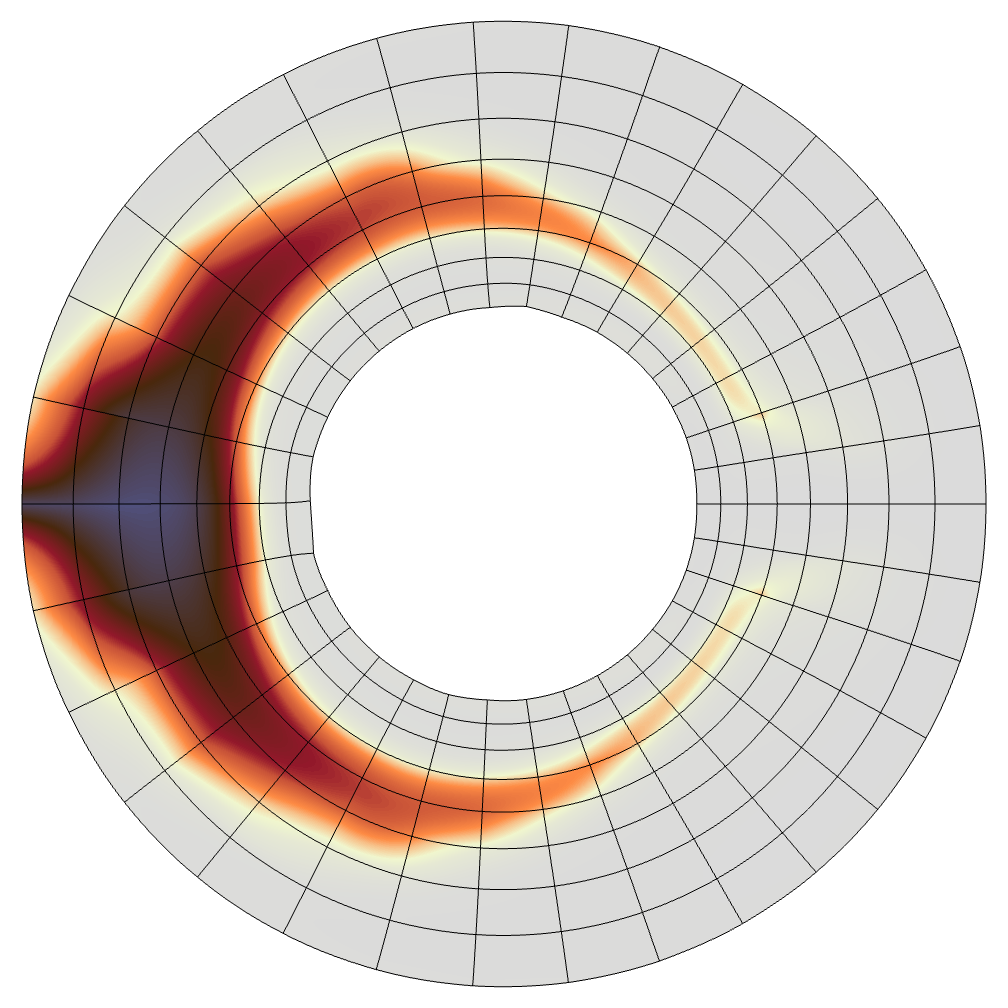}}
\end{tabular}
\caption{Plots of density, velocity, and specific internal energy for
         the Sedov blast on an annular domain with a blast origin on the
         outer boundary.}
\label{fig:annulus_2D}

\vspace{0.5em}

\centering
\begin{tabular}{c c c}
  & Density & Velocity \\

  $t=2.4$
  & \raisebox{-.5\height}{\includegraphics[width=0.30\textwidth]{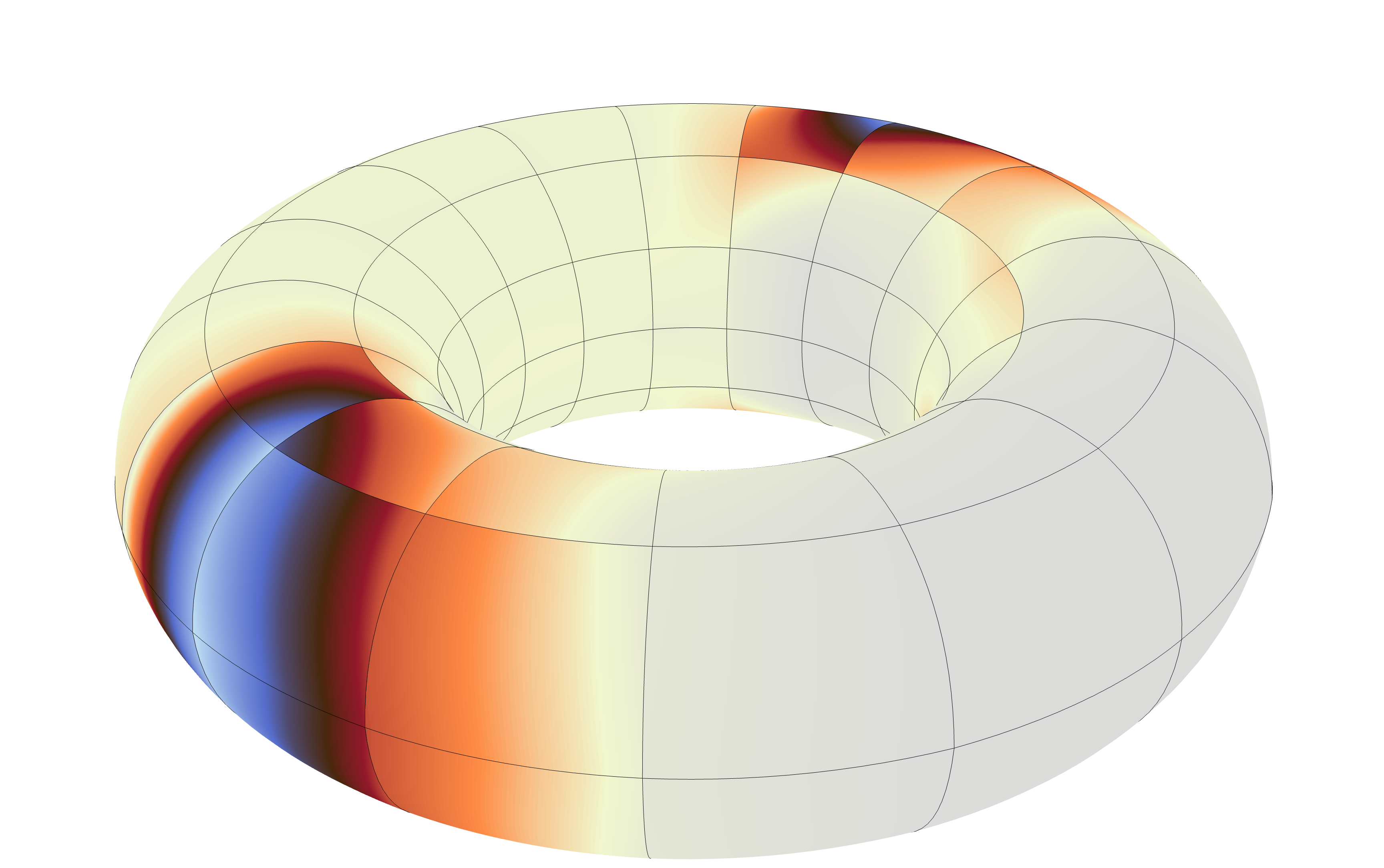}}
  & \raisebox{-.5\height}{\includegraphics[width=0.30\textwidth]{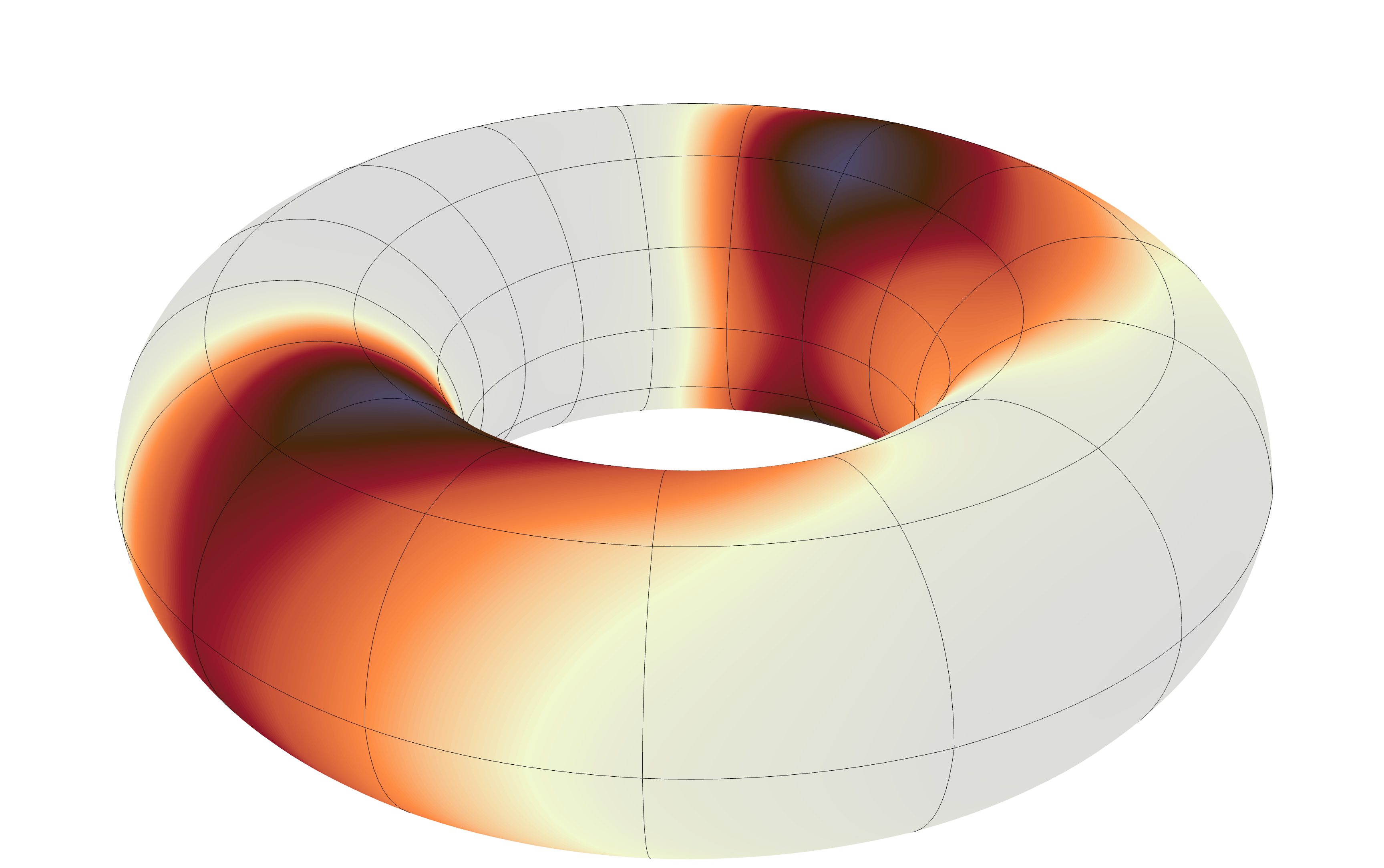}} \\

  $t=4.8$
  & \raisebox{-.5\height}{\includegraphics[width=0.30\textwidth]{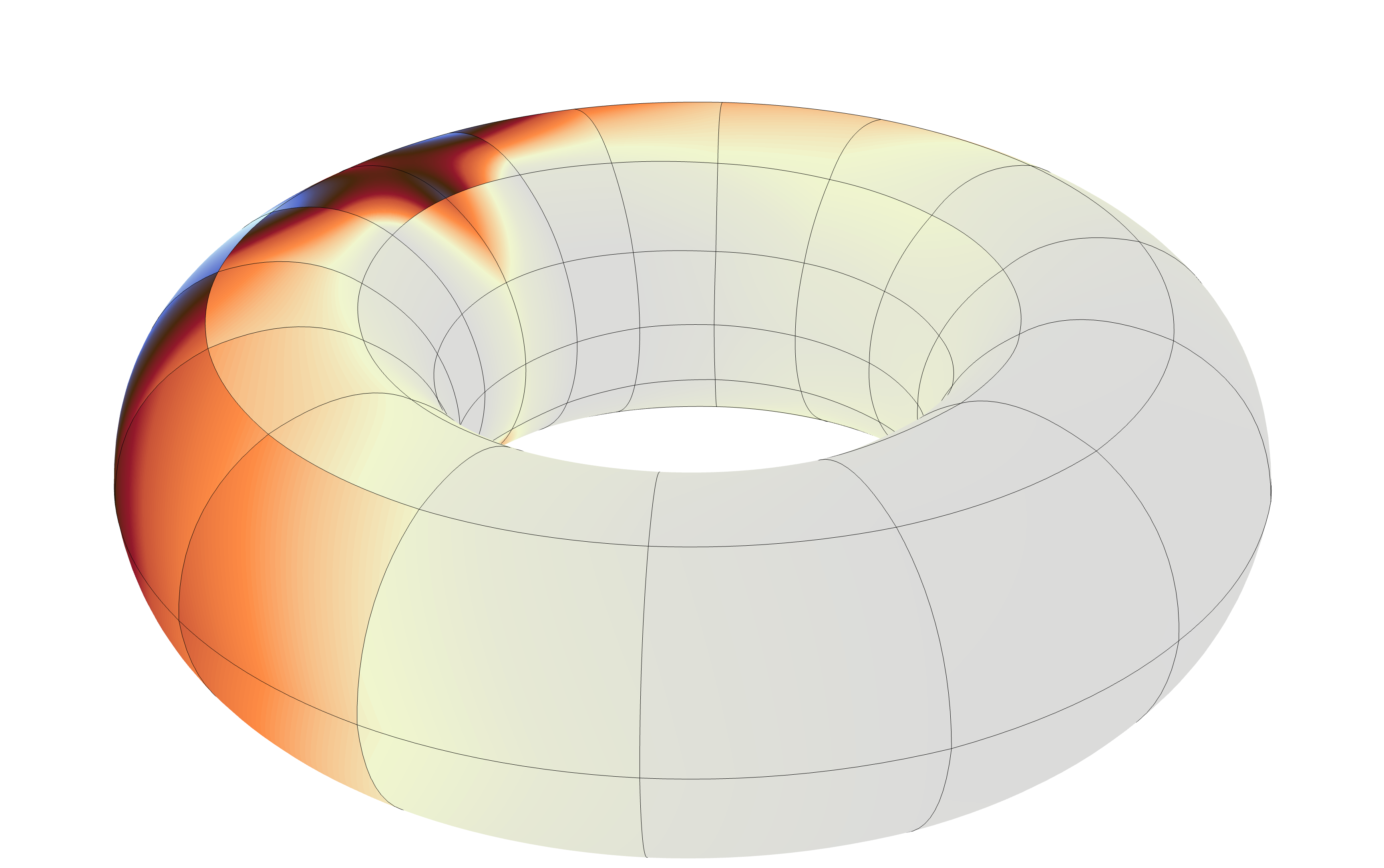}}
  & \raisebox{-.5\height}{\includegraphics[width=0.30\textwidth]{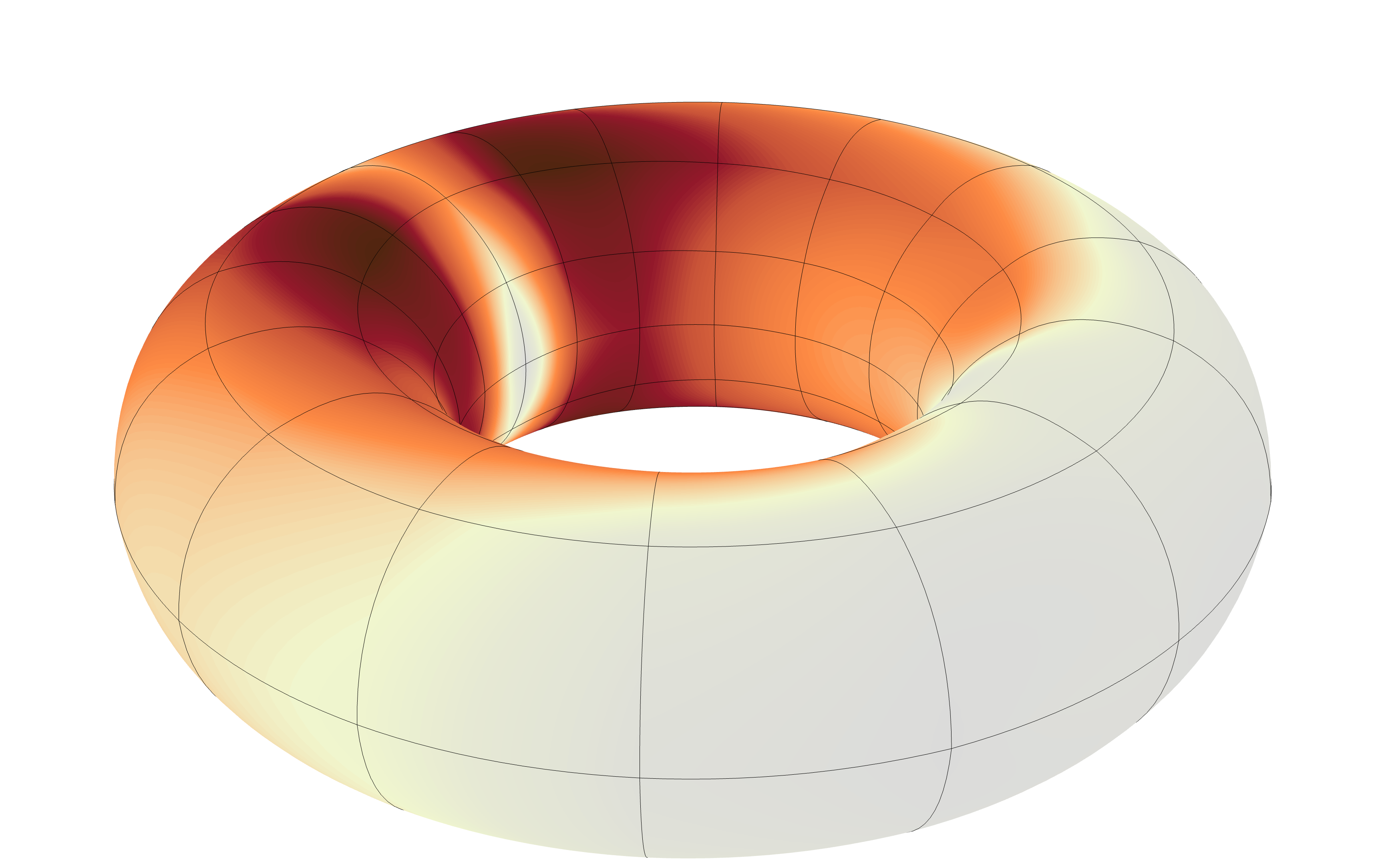}}
\end{tabular}
\caption{Plots of density and velocity for the Sedov blast on a solid torus.}
\label{fig:torus_3D}
\end{figure}

The 3D solid torus has the major radius $R = 0.7$ and the minor radius $r = 0.3$.
Its boundary is parameterized by
\[
    x(u,v) = (R + r \cos(v)) \cos(u), \quad
    y(u,v) = (R + r \cos(v)) \sin(u), \quad
    z(u,v) = r \sin(v),
\]
for $u,v \in [0, 2\pi)$. The blast origin is on the outer edge of the torus.
The annulus results in Figure \ref{fig:annulus_2D} and the solid-torus
results in Figure \ref{fig:torus_3D} show the expected behavior in these more
demanding geometries. The shock repeatedly interacts with curved boundary
segments, and the solution continues to follow the wall tangentially without
visible boundary-induced mesh deformation or solution oscillations.

\section{Conclusion}

We have presented a high-order finite element ALE method for slip-wall
boundary conditions on analytically curved domains, combining weak
Lagrangian boundary enforcement, parameterized TMOP-based mesh optimization,
and conservative remap in a way that preserves curved geometric features and
allows tangential boundary motion.
Future work will combine this approach with our multi-material ALE technology
\cite{ALE2018} and test it on more complex applications, as well as
develop the partial-assembly kernels and performance optimizations
\cite{MARBL2022, MARBL2024} needed for integration in a production code.

\paragraph{License Notice}
This manuscript has been authored by Lawrence Livermore National Security,
LLC under Contract No. DE-AC52-07NA27344 with the U.S. Department of Energy.
The United States Government retains, and the publisher, by accepting the
article for publication, acknowledges that the United States Government retains
a non-exclusive, paid-up, irrevocable, world-wide license to publish or
reproduce the published form of this manuscript, or allow others to do so,
for United States Government purposes.

\bibliographystyle{cas-model2-names}
\bibliography{references}
\end{document}